\documentstyle[12pt]{article}
\newlength{\abstractwidth}
\renewcommand{\thefootnote}{\fnsymbol{footnote}}
\renewcommand{\thanks}[1]{\footnote{#1}} 
\newcommand{\starttext}{ \setcounter{footnote}{0}
\renewcommand{\thefootnote}{\arabic{footnote}}}

\newcommand{\be}{\begin{equation}}
\newcommand{\bea}{\begin{eqnarray}}
\newcommand{\eea}{\end{eqnarray}} \newcommand{\ee}{\end{equation}}
 
\renewcommand{\>}{\rangle} \def\ba{\begin{eqnarray}}
\def\ea{\end{eqnarray}}

\def\cX{{\cal X}}
\def\cL{{\cal L}}
\def\C{{\bf C}}
\def\r{\rho}
\def\Aut{{\rm Aut}\,}

\def\ra{\rightarrow}

\def\o{\omega}

\def\det{{\rm det}}

\def\log{\,{\rm log}\,}
\def\exp{\,{\rm exp}\,}

\def\o{\omega}

\def\al{\alpha}
\def\b{\beta}
\def\g{\gamma}
\def\d{\delta}
\def\e{\varepsilon}

\def\l{\lambda}
\def\m{\mu}

\def\o{\omega}

\def\r{\rho}
\def\si{\sigma}

\def\na{\nabla}

\def\le{\leq}

\def\ti{\tilde}

\def\Z{{\bf Z}}

\def\R{{\bf R}}
\def\C{{\bf C}}
\def\P{{\bf P}}

\def\i{\infty}
\def\I{\int}

\def\s{\sum}

\def\ddb{{\partial\bar\partial}}

\def\sub{\subseteq}
\def\ra{\rightarrow}

\def\cH{{\cal H}}

\def\cH{{\cal H}}

\def\na{{\nabla}}
\def\us{{\underline s}}

\def\cT{{\cal T}}

 \def\v{\vskip .1in}

\def\[{{\bf [}}
\def\]{{\bf ]}}

\def\pl{\partial}

\begin{document}
\starttext \baselineskip=18pt \setcounter{footnote}{0}
\newtheorem{theorem}{Theorem}
\newtheorem{lemma}{Lemma}
\newtheorem{corollary}{Corollary}
\newtheorem{definition}{Definition}
\newtheorem{conjecture}{Conjecture}
\begin{center}
{\Large \bf LECTURES ON STABILITY AND CONSTANT SCALAR
CURVATURE\footnote{Research supported in part by National Science
Foundation grants  DMS-02-45371 and DMS-05-14003.
Contribution to {\it Current Developments in Mathematics 2007},
Harvard University, November 16-17, 2007.}}
\\
\bigskip
\bigskip

{\large D.H. Phong$^*$ and Jacob Sturm$^\dagger$} \\

\bigskip

$^*$ Department of Mathematics\\
Columbia University,
New York, NY 10027\\

\v

$^{\dagger}$ Department of Mathematics \\
Rutgers University, Newark, NJ 07102

\bigskip

\begin{abstract}
{\small An introduction is provided to some current research trends
in stability in geometric invariant theory and the problem of
K\"ahler metrics of constant scalar curvature. Besides classical
notions such as Chow-Mumford stability, the emphasis is on several
new stability conditions, such as K-stability, Donaldson's
infinite-dimensional GIT, and conditions on the closure of orbits of
almost-complex structures under the diffeomorphism group. Related
analytic methods are also discussed, including estimates for energy
functionals, Tian-Yau-Zelditch approximations, estimates for moment
maps, complex Monge-Amp\`ere equations and pluripotential theory,
and the K\"ahler-Ricci flow.}

\end{abstract}

\end{center}

\vfill\break

\baselineskip=15pt
\setcounter{equation}{0}
\setcounter{footnote}{0}

\tableofcontents

\vfill\break

\section{Introduction}
\setcounter{equation}{0}

A central theme in geometry is to try and characterize a given geometric structure by
a metric with ``best'' curvature properties. A classic example is the uniformization
theorem, which says that a complex structure on a compact surface can be characterized
by a metric of constant curvature. More recently, in the case of holomorphic vector
bundles over a compact K\"ahler manifold, the algebraic-geometric notion
of stability in the sense of Mumford-Takemoto has been shown by Donaldson \cite{D87}
and Uhlenbeck-Yau \cite{UY} to be equivalent
to the existence of a Hermitian-Einstein metric.

\medskip

In general, a metric may not be characterized by curvature
properties, since the number of degrees of freedom may not even
match. But in K\"ahler geometry, there is a natural such question:
Given a positive line bundle $L\to X$ over a compact complex
manifold $X$, determine when there exists a K\"ahler form $\o\in
c_1(L)$ with constant scalar curvature \be R=\mu n. \ee Here we have
denoted the dimension of $X$ by $n$ and $\mu$ is a constant. When
$L=K_X^{-\mu}$, it is easy to see that the condition that $R=\mu n$
is equivalent to the K\"ahler-Einstein condition, \be R_{\bar
kj}=\mu g_{\bar kj} \ee so that the search for metrics of constant
scalar curvature also encompasses the famous problem of finding
K\"ahler-Einstein metrics. This problem was solved by Yau \cite{Y78}
and Aubin \cite{A} when $c_1(X)<0$. It was solved by Yau \cite{Y78}
when $c_1(X)=0$, as part of his solution of the Calabi conjecture.
When $c_1(X)>0$, the surface case was treated in \cite{TY87, T90},
but the general case of higher dimensions is still open. More
generally, the main question in the field is the conjecture of Yau
\cite{Y93}, still open at this time, which says that the existence
of a K\"ahler form $\o\in c_1(L)$ with constant scalar curvature
should be equivalent to the stability of $c_1(L)$ in the sense of
geometric invariant theory.

\medskip

The conjecture of Yau can be viewed as a version of the Donaldson-Uhlenbeck-Yau
theorem for manifolds instead of vector bundles. As such, the corresponding
partial differential equations are more non-linear than the equations
for Hermitian-Einstein metrics. But there is a
significant additional difficulty: it is that, unlike the notion of Mumford-Takemoto
stability in the case of vector bundles, it is not yet clear what is the
correct notion of stability, and in fact, finding this correct notion has to
be viewed as a major component of the problem.

\medskip
The present lecture notes are based on introductory talks on the
subject given at the George Kempf Memorial Lectures at Johns Hopkins
University in October 2007, and at the Current Developments in
Mathematics conference at Harvard University in November 2007. An
earlier and shorter version of the lecture notes had been circulated
informally at Harvard at the time of the conference. Their goal is
to provide an accessible introduction to the various notions of
stability which have been introduced in the context of constant
scalar curvature metrics, as well as to some corresponding analytic
results. Besides the classical notions of Chow-Mumford and
Hilbert-Mumford stability, we discuss the notions of analytic and
algebraic $K$-stability due respectively to Tian \cite{T97} and
Donaldson \cite{D02}, Donaldson's infinite-dimensional GIT
\cite{D99}, and stability conditions such as (B) and (S) which arise
naturally in the context of the K\"ahler-Ricci flow
\cite{PS06,PSSW2}. We also discuss several analytic methods,
including inequalities for energy functionals, Tian-Yau-Zelditch
approximations, estimates for moment maps, degenerate complex
Monge-Amp\`ere equations, and the K\"ahler-Ricci flow.

\medskip
Certain aspects of our presentation may be worth mentioning:

\smallskip
(a) It is well-known that
the Mabuchi $K$-energy functional plays a central role in the theory, first in the variational
formulation of the analytic problem, and second as an important link to the notion
of $K$-stability. It does not appear to be as widely appreciated that
the Aubin-Yau functional $F_{\o_0}^0$ functional plays a similar role,
partly because the Euler-Lagrange equation for
$F_{\o_0}^0$ seems unrelated to the constant scalar curvature equation,
and the relevance of $F_{\o_0}^0$ can only be seen upon restriction to the
spaces ${\cal K}_k$ of Bergman metrics. We have thus tried to stress these points,
by grouping in a single chapter \S 5 its key properties: the classic
result of Zhang \cite{Z96} identifying the critical points of $F_{\o_0}^0$
on ${\cal K}_k$ with balanced imbeddings, and thus establishing a bridge to
Chow-Mumford stability; the observation due to Donaldson \cite{D05} that
$R-\mu n$ can indeed be interpreted as the Euler-Lagrange equation for $F_{\o_0}^0$
restricted to ${\cal K}_k$; and the basic relation between $F_{\o_0}^0$
and Monge-Amp\`ere masses for a path in the space ${\cal K}$ of K\"ahler potentials
\cite{PS06}.

\smallskip
(b) In \S 7.2, we have simplified to some extent some technical aspects of Donaldson's original
proof \cite{D01} of the necessity of Chow-Mumford stability. Of course, the
underlying motivation from symplectic geometry remains,
in particular the
key insight of Donaldson,
based on Lu's formula \cite{Lu}, that the constancy of the scalar curvature should be closely related
to the constancy of the density of states.
But the simplification of some important estimates which no longer
requires the formalism of infinite-dimensional moment maps may make the proof more
accessible.

\smallskip
(c) In \S 8.3.1-8.3.3, we have taken the opportunity to discuss certain aspects of the K\"ahler-Ricci
flow which seem not to have been covered fully in the literature. In particular,
we describe inequalities between energy functionals along the K\"ahler-Ricci flow.
Combined with Perelman's recent results \cite{Pe} and subsequent uniform Sobolev estimates
\cite{Ye, Zq}, they readily imply a version of Perelman's unpublished result, namely,
that when ${\rm Aut}^0(X)=0$, the existence of a K\"ahler-Einstein metric implies
the convergence in $C^\infty$ of the K\"ahler-Ricci flow.

\smallskip
(d) In \S 12.3.1, we have extracted from \cite{PS06, PS06a} a more general version
of the Ansatz producing generalized solutions of the degenerate complex Monge-Amp\`ere
equation than was originally stated in these papers. We would like to thank S. Zelditch
for stressing to us that such a version may be of interest.

\smallskip

The field is  vast and developing very
rapidly, and our list of topics is necessarily very incomplete.
For example, there is a very rich literature on geometric constructions
of metrics of constant scalar curvature, which we did not discuss, see e.g. \cite{ACGT, AP, KLP, LS,
RS, Si1}. There is also a related exposition of O. Biquard in the s\'eminaire N. Bourbaki
2004-2005 \cite{B06}. Nevertheless,
we hope that these notes can be useful to students interested in getting
a quick sense for certain trends in the subject.

\newpage

\section{The Conjecture of Yau}
\setcounter{equation}{0}

We begin with some background and notation. Let $X$ be a compact
complex manifold of dimension $n$. A K\"ahler form is a $(1,1)$ form
$\o={i\over 2}g_{\bar kj}\,dz^j\wedge d\bar z^k$ on $X$ which is
closed and strictly positive. The Ricci curvature tensor $R_{\bar
kj}$ and the scalar curvature $R$ of the corresponding metric
$g_{\bar kj}$ are given by
\be
\label{ricci}
R_{\bar
kj}=-\pl_j\pl_{\bar k}\,\log\,\o^n, \qquad R=g^{j\bar k}R_{\bar kj}.
\ee
The Ricci curvature form $Ric(\o)$ is the $(1,1)$ form defined by
\be
Ric(\o)={i\over 2}R_{\bar kj}dz^j\wedge d\bar z^k.
\ee
Since $\o^n$
is a metric on the anti-canonical bundle $K_X^{-1}$, $Ric(\o)=-{i\over 2}\ddb \log\,\o^n$ can
be viewed as the curvature form of $K_X^{-1}$ with respect to the
metric $\o^n$, and, as such, must be in the cohomology
class\footnote{Strictly speaking, $Ric(\o)\in\pi\,c_1(X)$.
We omit such factors of $\pi$ for notational simplicity.}
$c_1(K_X^{-1})\equiv
c_1(X)$.
Sometimes, we also denote the scalar curvature $R$ by
$R=R(\o)$, to emphasize its dependence on the K\"ahler form $\o$.

\subsection{Constant scalar curvature metrics in a given K\"ahler
class}

Let now $L\to X$ be a positive line bundle, that is, a holomorphic line bundle
admitting a metric $h_0$ whose curvature $\o_0\equiv -{i\over 2}\ddb\log\,h_0\in c_1(L)$
is a positive $(1,1)$-form. Then $\o_0$ equips $X$ with a K\"ahler structure. The main
question addressed in the present lecture series is whether there exists a K\"ahler form
$\o\in c_1(L)$ whose corresponding scalar curvature $R(\o)$ is constant.

\medskip
The condition of constant scalar curvature is not particularly rigid for a general Riemannian
metric, since it is a single scalar condition on an object with a much higher number of degrees
of freedom. However, the situation changes drastically with the above additional constraint
that the metric be a K\"ahler metric in a given K\"ahler class. For example, if $L=K_X^{-1}$,
then we have the  equivalence
\be
\label{kaehlereinstein}
R(\o)=n\ \Longleftrightarrow\ Ric(\o)=\o,
\ee
for $\o\in c_1(L)$. Indeed, as we noted above, for any K\"ahler form $\o$,
we have $Ric(\o)\in c_1(K_X^{-1})$. If we require that $\o\in c_1(K_X^{-1})$,
then $Ric(\omega)$ and $\omega$ are cohomologous,
and by the $\ddb$ lemma,
we have
\be
\label{riccipotential}
Ric(\o)-\o={i\over 2}\ddb f,
\ee
for some $f\in C^\infty(X)$ defined uniquely up to an additive constant.
The function $f$ is called the Ricci potential, and will play an important role
in the sequel.
Contracting both sides with $g^{j\bar k}$, we obtain $R(\o)-n=\Delta f$.
If $R(\o)$ is constant, the left hand side of this identity is constant. But the
range of $\Delta$ is orthogonal to constants, and thus this constant must be $0$.
The function $f$ is then constant, and hence $Ric(\o)-\o=0$. This establishes the
claim.

\subsection{The special case of K\"ahler-Einstein metrics}

The problem of finding K\"ahler-Einstein metrics, that is, metrics
satisfying
\be
Ric(\o)=\mu\, \o,
\ee
is one of the most celebrated problems in complex geometry.
Here $\mu$ is a constant that can be normalized to be $-1,0,$ or $1$.
When $\mu\not=0$, the K\"ahler-Einstein problem is the special case of the
constant
scalar curvature problem for metrics in the K\"ahler class $c_1(L)$,
with $L=K_X^{-\mu}$. The case $\mu=0$ is a special case of the Calabi conjecture,
which asserts that for any compact K\"ahler manifold $(X,\o)$, and any given $(1,1)$-form
$T\in c_1(X)$, there exists
within the K\"ahler class of $\o$ a metric with Ricci curvature equal to $T$.

\medskip
The K\"ahler-Einstein problem when $\mu =-1$ was solved independently by Yau \cite{Y78}
and Aubin \cite{A}. The Calabi conjecture was solved by Yau \cite{Y78}. However,
the case $\mu=1$ is still open, and henceforth we refer only to this case when
we speak of the K\"ahler-Einstein problem.

\medskip
In general, there are obstructions to the existence of K\"ahler-Einstein metrics
of positive scalar curvature. A classic obstruction is the theorem of Matsushima
\cite{Mat}, which says that if there exists a constant scalar curvature metric,
then the automorphism group of $X$ would have to be
reductive. Another obstruction, due to Futaki \cite{Fu1}
and also related to the automorphisms of $X$, is the
vanishing of the Futaki invariant defined as follows. Given any K\"ahler
metric $\o\in c_1(X)$, let $f$ be its Ricci potential,
as defined earlier by (\ref{riccipotential}).
Then the Futaki invariant $Fut$ is the character on $H^0(X,T^{1,0})$ defined by
\be
Fut(V)=\int_X (Vf)\,\o^n,\qquad V\in H^0(X,T^{1,0}).
\ee
The key propery of $Fut(V)$ is that it is actually independent of the choice of
$\o$ within $c_1(X)$. Thus, if $c_1(X)$ admits a K\"ahler-Einstein metric,
then $Fut$ must vanish identically.

\medskip
The vanishing of the Futaki invariant was shown to imply the existence of a K\"ahler-Einstein
metric when $dim\,X=2$ \cite{T90} and when $X$ is a toric variety \cite{WZ}. However,
a counterexample was provided by Tian in 1997 \cite{T97} of a compact complex manifold
$X$ with $c_1(X)>0$, no holomorphic vector fields, and yet no K\"ahler-Einstein
metrics (see \S 7.1 for a description of this and some other results in \cite{T97}).

\subsection{The conjecture of Yau}

The existence of K\"ahler-Einstein metrics, and more generally, of constant scalar
curvature metrics in a given K\"ahler class $c_1(L)$, is expected to be related to
deeper properties of the bundle $L\to X$. In fact, the guiding light of much of the
current research in the area has been the conjecture of Yau \cite{Y93},
which says that the existence of $\o\in c_1(L)$ with $R(\o)$ constant
should be equivalent to the stability of $L\to X$ in the sense of geometric
invariant theory.

\medskip
This conjecture has clearly a strong analogy with the case of holomorphic vector bundles.
Let $E\to X$ be a holomorphic vector bundle over a compact K\"ahler manifold $(X,\o)$.
Given a Hermitian metric $H_{\bar\al\beta}$ on $E$,
let $F_{\bar kj}{}^\al{}_\beta=-\pl_{\bar k}(H^{\al\bar\g}\pl_j H_{\bar\g\beta})$
be its curvature. A metric $H_{\bar\al\beta}$ is said to be Hermitian-Einstein
if $g^{j\bar k}F_{\bar kj}{}^\al{}_\beta=\mu\delta^\al{}_\beta$ for some constant $\mu$.
The theorem of Donaldson-Uhlenbeck-Yau \cite{D87, UY} asserts that
$E\to X$ admits a Hermitian-Einstein metric if and only if $E\to X$ is stable
in the sense of Mumford-Takemoto.

\medskip
As in the case of Hermitian-Einstein metrics and Mumford-Takemoto stability, a particularly
striking aspect
of the conjecture of Yau is that it asserts the equivalence between the existence of
a solution to a non-linear partial differential equation and a global, algebraic-geometric
property of the underlying space.

\newpage

\section{The Analytic Problem}
\setcounter{equation}{0}

In this section, we begin by setting up the problem
from the analytic point of view. Analytically, there are
several possible formulations and approaches.

\subsection{Fourth order non-linear PDE and Monge-Amp\`ere equations}

The most direct formulation of the problem is as a non-linear PDE in the potential $\phi$.
More precisely, fixing a metric $h_0$ on $L$ with curvature $\o_0=-{i\over 2}\ddb\log h_0>0$,
we seek another metric $h=h_0e ^{-\phi}$ with curvature $\o=-{i\over 2}\ddb\log h>0$
so that $R(\o)$ is constant. This means that $\phi$
must satisfy the ``$\o_0$-plurisubharmonicity" constraint
\be
\label{plurisubharmonicity}
\o_0+{i\over 2}\ddb\phi>0
\ee
and the partial differential equation
\be
\label{4thorder}
-g^{j\bar k}\pl_j\pl_{\bar k}\log \,(\o_0+{i\over 2}\ddb\phi)^n= \bar R
\ee
where $g_{\bar k j}$ is the metric corresponding to
the K\"ahler form $\o$. The value of the constant $\bar R$
is cohomological: clearly, $\bar R$ must be given by the average of the scalar
curvature, and thus
\be
\label{mu}
\bar R={1\over V}\int_X R\, \o^n
=
{n\over V}\int_X Ric(\o)\wedge \o^{n-1}=
n
{[c_1(X)][c_1(L)]^{n-1}
\over [c_1(L)]^n}
\equiv n\,\mu.
\ee
Here $V$ is the volume of $X$ and $\mu$ is the cohomological constant defined by
\be
\label{volume}
V=\int_X \o^n= [\pi\,c_1(L)]^n,
\qquad
\mu=
{[c_1(X)][c_1(L)]^{n-1}
\over [c_1(L)]^n}.
\ee
Thus the problem can be viewed as a 4th-order non-linear elliptic PDE in
$\phi$.

\medskip
In the special case $L=K_X^{-1}$, we have seen that the constant scalar curvature
condition is  equivalent to the constant Ricci curvature condition. This last
condition is well-known to be equivalent to a complex elliptic Monge-Amp\`ere equation.
Indeed, we always have $Ric(\o_0)\in c_1(K_X^{-1})$, and
if $\o_0\in c_1(L)=c_1(K_X^{-1})$ also, then we can write
$Ric(\o_0)-\o_0={i\over 2}\ddb f_0$, where $f_0$ is the Ricci potential of $\o_0$.
It is now easily seen,
simply by taking $-\ddb\,\log$ of both sides,
that
the following complex Monge-Amp\`ere equation
\be
\label{volumeequation}
(\o_0+{i\over 2}\ddb \phi)^n=e^{f_0-\phi}\o_0^n
\ee
for $\phi$ still satisfying the $\o_0$-plurisubharmonicity constraint
$\o_0+{i\over 2}\ddb\phi>0$, is equivalent to the constant Ricci curvature equation.
Geometrically, this means that, when $L=K_X^{-1}$,
the 4th order constant scalar curvature equation on the potential can been reduced to a
2nd order equation in $\phi$, involving
the volume of the K\"ahler metric $\o=\o_0+{i\over 2}\ddb\phi$.

\subsection{Geometric heat flows}

One approach to finding a solution to a given equation is to interpret it as the fixed point
of a dynamical system. In the present case, this amounts to replacing
the elliptic non-linear PDE by a non-linear parabolic flow. The problem becomes
then that of the long-time existence and convergence of the flow.

\medskip
For general $L$, a parabolic version of the equation (\ref{4thorder})
is the following 4th order parabolic flow,
\be
\label{calabiphi}
\dot\phi= R-\bar R, \qquad \phi(0)=c_0,
\ee
for $\phi$ satisfying the $\o_0$-plurisubharmonicity constraint
(\ref{plurisubharmonicity}).
This is equivalent to the Calabi flow \cite{Calabi}, which is the following
flow for metrics $g_{\bar kj}$,
\be
\label{calabig}
\dot g_{\bar kj}=\pl_j\pl_{\bar k}R, \qquad g_{\bar kj}(0)=g_{\bar kj}^0.
\ee

Just as in the elliptic case, when $L=K_X^{-1}$,
we need only consider a second-order parabolic flow,
namely the K\"ahler-Ricci flow. This is the K\"ahler version
of the Ricci flow introduced by Hamilton \cite{H82} (see also
\cite{CH, ChKnopf} for a detailed treatment, with an extensive list
of references). It can
be written either as
a flow of metrics,
\be
\label{KRmetric}
\dot g_{\bar kj}=- (R_{\bar kj}-\mu g_{\bar kj}),
\qquad g_{\bar kj}(0)=g_{\bar kj}^0
\ee
(with $\mu$ the cohomological value defined by (\ref{mu})),
or as a parabolic Monge-Amp\`ere equation for the potential $\phi$,
\be
\label{KRpotential}
\dot\phi= \log {(\o_0+{i\over 2}\ddb\phi)^n\over\o_0^n}+\mu \phi-f_0,
\qquad\phi(0)=c_0,
\ee
with $f_0$ the Ricci potential (c.f. (\ref{riccipotential}))
for the original K\"ahler form $\o_0$.

\subsection{Variational formulation and energy functionals}

Another approach which plays an important role in the sequel
is the variational approach. Thus we seek a functional whose
Euler-Lagrange equation is precisely the given equation,
so that the problem reduces to determining whether the
functional admits a critical point.

\subsubsection{The Mabuchi K-energy $K_{\o_0}(\phi)$}

A first important fact in the theory, established by Mabuchi \cite{Ma86},
is that the equation $R(\o)-\bar R=0$ for $\o$ in a given K\"ahler class
$c_1(L)$, $L\to X$ positive, is indeed realizable as an
Euler-Lagrange equation. More precisely, there is a functional $K_{\o_0}(\phi)$,
now called the Mabuchi $K$-energy and
defined on the space ${\cal K}$ of K\"ahler potentials,
\be
\label{calK}
{\cal K}
\equiv\{\,\phi\in C^\infty(X)\,;\ \o_\phi\equiv\o_0+{i\over 2}\ddb\phi>0\,\}
\ee
satisfying
\be
\label{Kvariation}
\delta K_{\o_0}(\phi)
=
-{1\over V}\int_X \delta\phi\,(R-\bar R)\o_\phi^n.
\ee
Since the $K$-energy is characterized by its variation, its exact definition depends
on the choice of a reference metric, which we have taken to be $\o_0$.
An explicit expression for $K_{\o_0}(\phi)$ is (see e.g. the derivation in Theorem
\ref{delignepairingstheorem} below)
\be
\label{Kexplicit}
K_{\o_0}(\phi)
=
{1\over V}
\bigg[\int_X(\log {\o_\phi^n\over\o_0^n})\o_\phi^n
-
\phi\sum_{j=0}^{n-1}Ric(\o_0)\o_\phi^j\o_0^{n-1-j}
+
\mu\sum_{j=0}^n\phi \o_\phi^j\o_0^{n-j}\bigg],
\ee
although we shall mostly use only its characterizing variational property
(\ref{Kvariation}). We note that $K_{\o_0}(\phi+c)=K_{\o_0}(\phi)$ for constant $c$,
so that $K_{\o_0}(\phi)$ descends to a functional on the space of K\"ahler metrics
in $c_1(L)$. Also, the variational formula for $K_{\o_0}(\phi)$ implies
the following important cocycle property for $\phi,\phi+\psi\in{\cal K}$,
\be
\label{Kcocycle}
K_{\o_0}(\phi+\psi)=K_{\o_0}(\phi)+K_{\o_\phi}(\psi).
\ee

\subsubsection{The Aubin-Yau functional $F_{\o_0}^0(\phi)$}

There is another functional in the theory, the Aubin-Yau functional
$F_{\o_0}^0(\phi)$, whose role is as important as that of the Mabuchi K-energy,
but which is more subtle.
First, let $J_{\o_0}(\phi)$
be the functional on ${\cal K}$ given by
\be
J_{\o_0}(\phi)
=
{i\over 2V}
\sum_{j=0}^{n-1}{n-j\over n+1}\int_X
\pl\phi\wedge\bar\pl\phi\wedge \o_\phi^{n-1-j}\wedge\o_0^j.
\ee
Then the functional $F_{\o_0}^0(\phi)$ is defined by
\be
\label{aubinyau}
F_{\o_0}^0(\phi)
=
J_{\o_0}(\phi)-{1\over V}\int_X\phi\,\o_0^n.
\ee
Its variation can be readily verified to be
\be
\label{aubinyauvariation}
\delta F_{\o_0}(\phi)= -{1\over V}\int_X \delta\phi\,\o_\phi^n,
\ee
and it also satisfies the cocycle condition
\be
\label{Fcocycle}
F_{\o_0}^0(\phi+\psi)=F_{\o_0}^0(\phi)+F_{\o_\phi}^0(\psi),
\qquad \phi,\phi+\psi\in {\cal K}.
\ee
In the case of $L=K_X^{-1}$, the Monge-Amp\`ere equation (\ref{volumeequation})
for K\"ahler-Einstein metrics can be viewed as
an Euler-Lagrange equation for $F_{\o_0}^0(\phi)$ with constraints. More precisely,
let
\be
F_{\o_0}(\phi)=F_{\o_0}^0(\phi)-\log\bigg({1\over V}\int_X e^{f_0-\phi}\o_0^n\bigg).
\ee
Then we have
\be
{\delta F_{\o_0}\over\delta\phi}=0
\
\Longleftrightarrow
\
(\o_0+{i\over 2}\ddb \phi)^n=e^{f_0-\phi}\o_0^n.
\ee
For general $L$, there does not appear to be any evident relation
between $F_{\o_0}^0(\phi)$ and the constant scalar curvature equation $R(\o)-\bar R=0$.
It comes therefore as a surprise that $F_{\o_0}^0(\phi)$ is actually
intimately related to this equation. However, to see this relation, we have
to restrict $F_{\o_0}^0(\phi)$ to certain finite-dimensional subspaces
of ${\cal K}$, namely the subspaces ${\cal K}_k$
of Bergman metrics. We turn next to a description of these spaces.

\newpage

\section{The spaces ${\cal K}_k$ of Bergman metrics}
\setcounter{equation}{0}

The conjecture of Yau links the existence of a transcendental object,
namely a K\"ahler metric with
constant scalar curvature, to an
algebraic-geometric condition, namely stability in GIT.
This last concept depends fundamentally on the realization of $X$
as a projective variety by Kodaira imbeddings.
Associated to each such imbedding is an ``algebraic'' metric on $L$,
namely the induced Fubini-Study metric from $O(1)$.
The Tian-Yau-Zelditch theorem (TYZ) asserts that any metric on $L$ with
positive curvature can be approximated by such ``algebraic''
metrics. Thus the induced Fubini-Study metrics - or ``Bergman metrics'' -
provide a bridge between
the original analytic problem and the underlying algebraic-geometric structure.
The strategy of approximating transcendental objects by algebraic-geometric ones
has been advocated by Yau \cite{Y93} over the years.
It plays an essential role in many developments described in this paper,
in particular Donaldson's theorem
on the necessity of Chow-Mumford stability in \S 7.2,
the construction of geodesics in section \S 12.3, and Donaldson's lower
bound for the Calabi functional in \S 7.3.
In this section, we describe the spaces of Bergman metrics
and the TYZ theorem. In the next, we return to the functional $F_{\o_0}^0(\phi)$
and show how its significance in connection with the equation $R-\bar R=0$
emerges upon restricting it to the spaces of Bergman metrics.

\subsection{Kodaira imbeddings}

Let $L\to X$ be a positive line bundle over a compact complex manifold $X$.
Then for each basis $\us=\{s_\al(z)\}_{\al=0}^{N_k}$ of $H^0(X,L^k)$,
$N_k={\rm dim}\,H^0(X,L^k)$, and $k$
large enough, the Kodaira imbedding theorem asserts that the map
\be
\label{kodairaimbedding}
\iota_{\us}\ :\ X\in z\ \longrightarrow\ [s_0(z),\cdots,s_{N_k}(z)]\in {\bf CP}^{N_k}
\ee
is an imbedding. Under this imbedding, the bundle $O(1)$ over ${\bf CP}^{N_k}$ pulls-back
to $L^k$. Let $h_{FS}$ be the Fubini-Study metric on $O(1)$, and $\o_{FS}=-{i\over 2}\ddb \log\,h_{FS}$
the Fubini-Study metric on ${\bf CP}^{N_k}$. Then $\iota_\us^*(h_{FS})^{1/k}$ is a metric
of $L$ with positive curvature, and the space ${\cal K}_k$ of Bergman metrics is by
definition the space of all metrics of the form $\iota_\us^*(h_{FS})^{1/k}$
as the basis $\us$ varies,
\be
\label{calKk}
{\cal K}_k=\{\ \iota_\us^*(h_{FS})^{1/k}\ ;\ \us\ {\rm basis\ of}\ H^0(X,L^k)\ \}
\ee
If we fix a reference basis $\hat\us$, then any other basis can be obtained from
$\hat\us$ by an element of $GL(N_k+1)$, and, since $h_{FS}$ is invariant
under $SU(N_k+1)$, we obtain
\be
{\cal K}_k= SL(N_k+1)/SU(N_k+1).
\ee
Thus ${\cal K}_k$ can be viewed as a symmetric space with negative curvature.

\subsection{The Tian-Yau-Zelditch theorem}

Explicitly, we can view ${\bf CP}^N$ as ${\bf CP}^N=\{{\bf C}^{N+1}\setminus 0\}/\sim$
with $x,y\in{\bf C}^{N+1}\setminus 0$ equivalent if $x=\lambda y$.
As such, ${\bf CP}^N$ carries a natural line bundle, namely the universal bundle consisting
of the line $\ell_{[y]}=\{x=\lambda y\}$ over the point $[y]\in {\bf CP}^N$.
The $O(1)$ bundle over ${\bf CP}^N$ is the dual of the universal bundle.
Thus its fiber at each $[y]$ is the space of linear functionals $\ell_{[y]}\to{\bf C}$.
A basis of $H^0({\bf CP}^N,O(1))$ is provided by the linear functionals
$x\to x_\al$, for $0\le\al\leq N$. We can then define
the Fubini-Study metrics on $O(1)$ by
\be
\label{metricsprojective}
h_{FS}={1\over\sum_{\al=0}^N|x_\al|^2}\equiv{1\over |x|^2}.
\ee
Its curvature
$\o_{FS}={i\over 2}\,\ddb \,\log\,|x|^2$ is then the Fubini-Study metric on
${\bf CP}^N$.

\medskip
In the setting  $L\to X$ positive, imbedded into $O(1)\to{\bf CP}^{N_k}$ by
a basis $\us$ of $H^0(X,L^k)$, it follows that the induced metrics
on $L$ and on $X$ are given explicitly by
\bea
\label{induced}
\iota_{\us}^*(h_{FS})^{1/k}
=
{1\over (\sum_{\al=0}^{N_k}|s_\al(z)|^2)^{1/k}},
\qquad
{1\over k}\iota_{\us}^*(\o_{FS})
=
{i\over 2k}\ddb \log \sum_{\al=0}^{N_k}|s_\al(z)|^2.
\eea
Fix now a metric $h$ on $L$ with positive curvature $\o=-{i\over 2}\ddb\log h$.
Let $\us$ be an {\it orthonormal basis} for $H^0(X,L^k)$,
with respect to the $L^2$ metric defined on the sections of $L^k$ by the metric
$h$ and the volume form $\o$. Define the density of states $\rho_k(z)$
by
\be
\label{density}
\rho_k(z)=\sum_{\al=0}^{N_k}|s_\al(z)|^2 h^k(z)
\ee
Clearly, $\rho_k(z)$ is independent of the choice of $\us$ among the
orthonormal bases of $H^0(X,L^k)$. Its integral gives $N_k+1={\rm dim}\,H^0(X,L^k)$,
whence its name.
Introducing $\rho_k(z)$ converts
the expression (\ref{induced}) for
the induced metrics
\be
h(k)\equiv\iota_{\us}^*(h_{FS})^{1/k},
\qquad
\o(k)\equiv{1\over k}\iota_{\us}^*(\o_{FS})
\ee
into the following simple relation with
the original metric $h$ and its curvature $\o$,
\bea
\label{relation}
\log {h(k)\over h} = -{1\over k}\log\,\rho_k(z),
\qquad
\o-\o(k)
=-{i\over 2k}\,\ddb\log \rho_k(z).
\eea
The Tian-Yau-Zelditch theorem can be stated as follows:

\begin{theorem}
\label{tyz}
{\rm \cite{Y87, T90a, Z}}
The density of states $\rho_k(z)$ admits an asymptotic expansion
$\sum_{p=0}^\infty A_p(z)k^{n-p}$ with $A_0(z)=1$, $A_p(z)$ smooth
functions,
in the sense that
\be
\label{densityexpansion}
\|\rho_k(z)-\sum_{p=0}^MA_p(z)k^{n-p}\|_{C^L}\leq C_{LM}\,k^{n-M-1}.
\ee
In particular, we have the following approximations, for any
$C^L$ norm $\|\cdot\|$,
\be
\|\log{h(k)\over h}+n{\log\,k\over k}\|=O({1\over k^2}),
\qquad
\|\omega(k)-\omega\|=O({1\over k^2}).
\ee
\end{theorem}

The proof of (\ref{densityexpansion}) is built on the expansion
of Boutet de Monvel-Sj\"ostrand \cite{BoSj} for the Szeg\"o kernel
on strongly pseudoconvex domains, extending
an earlier expansion along the diagonal due to Fefferman \cite{F}.
Another independent proof, also built on the result of
Boutet de Monvel-Sj\"ostrand, is due to Catlin \cite{C}.
The approach in \cite{T90a} is based on H\"ormander's $L^2$ estimates,
and a generalization to open manifolds with bounds on the Ricci curvature
is also given there.
The TYZ theorem implies that an arbitrary metric $h$ in the space ${\cal K}$
of K\"ahler potentials can be approximated in $C^\infty$ by metrics $h(k)$ in ${\cal K}_k$.
In this sense, we have
\be
\label{heuristic}
{\cal K}={\rm lim}_{k\to\infty}{\cal K}_k.
\ee

The second leading coefficient $A_1(z)$ plays a key role in subsequent
developments (c.f. \S 5.2, \S 7.2 and \S 7.3)
and has been determined by Lu \cite{Lu99}:

\begin{theorem}
\label{lutheorem}
Let the set up be the same as in the previous theorem. Then
\be
\label{lu}
A_1(z)={1\over 2}R(z)
\ee
where $R(z)$ is the scalar curvature of the K\"ahler metric $\o$ on $X$.
\end{theorem}

\newpage
\section{The functional $F_{\o_0}^0(\phi)$ on ${\cal K}_k$}
\setcounter{equation}{0}

We return now to the discussion of the functional $F_{\o_0}^0(\phi)$.
Unlike the $K$-energy $K_{\o_0}(\phi)$, for general positive line bundles $L\to X$, there is no
obvious relation between $F_{\o_0}^0(\phi)$
and the constant scalar curvature equation $R-\bar R=0$.
However, remarkably, upon restriction to the spaces ${\cal K}_k$
of Bergman metrics, $F_{\o_0}^0(\phi)$ turns out to be closely
related both to the constant scalar curvature equation
and to stability conditions. In fact, it is related to
all forms of GIT stability conditions, whether it be Chow-Mumford,
$K$-stability, or Donaldson's infinite-dimensional GIT.
We explain this now in some detail.

\medskip
It is convenient to think of the restriction of $F_{\o_0}^0(\phi)$ to
${\cal K}_k$ as a functional on $SL(N_k+1)/U(N_k+1)$. To this end,
we fix a basis $\hat\us$ of $H^0(X,L^k)$, and identify $X$ with
the projective variety $\iota_{\hat\us}(X)$,
\be
X\ \longleftrightarrow \ \hat X\equiv \iota_{\hat\us}(X)\subset{\bf CP}^{N_k}.
\ee
Now for any basis $\us= \sigma\cdot\hat\us$
of $H^0(X,L^k)$ with $\sigma\in SL(N_k+1)$,
we have
\be
\iota_{\sigma\cdot\hat\us}(X)=\sigma\cdot \iota_{\hat\us}(X)
\ee
and the pull-back metrics $\iota_{\sigma\cdot\hat\us}^*(\o_{FS})$
and $\iota_{\hat\us}^*(\o_{FS})$
are related by
\be
\iota_{\sigma\cdot\hat\us}^*(\o_{FS})
=
{i\over 2}\ddb \log |\sigma x|^2
=
\iota_{\hat\us}^*(\o_{FS})+{i\over 2}\ddb\phi_\sigma
\ee
where $[x]\in \hat X$, and $\phi_\sigma(x)$ is the following potential
\be
\phi_\sigma(x)=\log{|\sigma x|^2\over |x|^2},
\qquad x\in {\bf C}^{N_k+1}.
\ee
Now $\iota_{\us}^*(\o_{FS})\in k\,c_1(L)$,
and it is natural to introduce the functional $F_{\iota_{\hat\us}^*(\o_{FS})}^0(\phi)$,
which is a functional on the space of K\"ahler potentials in $k\,c_1(L)$ rather
than in $c_1(L)$. The use of  $F_{\iota_{\hat\us}^*(\o_{FS})}(\phi)$
is particularly convenient when $X$ is identified with $\hat X$
and the extrinsic geometry of $\hat X$ as a projective variety play an
important role.
Clearly, we have in general
\be
F_{k\o_0}^0(k\phi)=k\, F_{\o_0}^0(\phi),
\qquad \phi\in {\cal K},
\ee
so that this scaling from $c_1(L)$ to $k\, c_1(L)$ is a question of emphasis rather than
one really of substance. Furthermore,
in view of the cocycle property (\ref{Fcocycle}),
any other choice of reference K\"ahler form $\o_0$ in $F_{k\o_0}(k\phi)$
can be brought back to $\iota_{\hat\us}^*(\o_{FS})$ when considering the derivatives
of $F_{k\o_0}^0(k\phi)$. Thus we have the following function of $\sigma$,
\be
SL(N_k+1)\ni \sigma\ \longrightarrow\ F_{\iota_{\hat\us}^*(\o_{FS})}^0(\phi_\sigma)
\ee
which can be viewed as a rescaling of the functional $F_{\o_0}^0(\phi)$
restricted to ${\cal K}_k$.

\subsection{$F_{\o_0}^0$ and balanced imbeddings}

A key feature of ${\cal K}_k=SL(N_k+1)/SU(N_k+1)$ is that any two potentials
there can be connected by a one-parameter subgroup $\sigma(t)=e^{(\delta\sigma) t}$
in $SL(N_k+1)$, $\delta\sigma\in sl(N_k+1)$. The derivatives of
$F_{\iota_{\hat\us}^*(\o_{FS})}^0(\phi_\sigma)$
along such
one-parameter subgroups can be evaluated explicitly:
\bea
\label{derivativesF0}
-{d\over dt}\bigg\vert_{t=0}F_{\iota_{\hat\us}^*(\o_{FS})}^0(\phi_{\sigma(t)})
&=&
{\rm Tr}((\delta\sigma+\delta\sigma^*)\cdot M),
\nonumber\\
-{d^2\over dt^2}
\bigg\vert_{t=0}F_{\iota_{\hat\us}^*(\o_{FS})}^0(\phi_{\sigma(t)})
&=&
\|\pi_{\cal N}V\|^2
\eea
where $M=M_{\bar\al\beta}$ is the following matrix
\be
\label{M}
M_{\bar\al\beta}={1\over V_{FS}}\int_{\hat X}{\bar x_\al x_\beta\over |x|^2}\o_{FS}^n
\ee
$V_{FS}$ is the volume of $\hat X$ with
respect to the Fubini-Study metric,
and $\pi_{\cal N} V$ is the projection on the normal bundle to $\hat X$
of the vector field corresponding to the action of $\delta\sigma$ on ${\bf CP}^{N_k}$.
In fact, the variation of $\phi_\sigma$ with respect to $\sigma$ is given by
\be
\delta\phi_\sigma={x^*\sigma^*(\delta \sigma+\delta\sigma^*)\sigma x\over x^*\sigma^*\sigma x}.
\ee
and thus, in view of the defining equation (\ref{aubinyauvariation})
for $\delta F_\o^0(\phi)$, we obtain
\be
{d\over dt}\bigg\vert_{t=0}F_{\iota_{\hat\us}^*(\o_{FS})}^0(\phi_\sigma)
=
-{1\over V_{FS}}\int_X\dot\phi_{\sigma(t)}\iota_{\sigma(t)\cdot\us}^*(\o_{FS})^n
\bigg\vert_{t=0}
=
-{1\over V_{FS}}\int_{\hat X}
{x^*(\delta\sigma+\delta\sigma^*)x\over x^*x}\o_{FS}^n
\ee
and the desired formula. The proof of the formula for the second derivative is slightly
longer, and can be found in \cite{PS03}. Note that it implies that
$F_{\iota_{\hat\us}^*(\o_{FS})}^0(\phi_\sigma)$
is convex along one-parameter subgroups of $SL(N_k+1)$.
Replacing $\hat X$ by $\sigma\cdot\hat X$, we obtain the following immediate
consequence of the formula for the first derivative:

\begin{theorem}
\label{balancedcritical}
A point $\sigma\in SL(N_k+1)$
is a critical point of the functional $F_{\o_0}^0$ restricted to ${\cal K}_k$
if and only if the image $\sigma\cdot\hat X$ of $X$ by the Kodaira
imbedding defined by the basis $\sigma\cdot\hat\us$
of $H^0(X,L^k)$ is a balanced submanifold of ${\bf CP}^{N_k}$, in the sense that
\be
\label{balanced}
{1\over V_{FS}}\int_{\sigma\cdot\hat X}{\bar x_\al x_\beta\over |x|^2}\o_{FS}^n
=
\lambda \,\delta_{\bar\al\beta}
\ee
for some scalar $\lambda$.
\end{theorem}

The notion of a balanced submanifold was introduced by Bourguignon, Li, and Yau \cite{BLY}.
As we shall see later, by a theorem of Zhang \cite{Z96} (see also Luo \cite{Luo}
and \cite{PS03}), the existence of a unique balanced imbedding is equivalent to
Chow-Mumford stability, so that the above theorem actually provides a first
link between the functional $F_\o^0(\phi)$ and stability in GIT.
A similar version of balanced imbeddings for vector bundles and its equivalence
with Gieseker stability is in X. Wang \cite{XWang}.

\medskip
The convexity of $F_{\iota_{\hat\us}^*(\o_{FS})}^0(\phi_\sigma)$ along one-parameter subgroups
provides yet another link to stability,
this time to the notion of $K$-stability to be described in greater detail in section \S 6.1.
Indeed, it implies that for any $\delta\sigma_k\in sl(N_k+1)$, the limit
\be
\mu_k\equiv -{\rm lim}_{t\to-\infty}\,V_{FS}\,\dot F_{\iota_{\hat\us}^*(\o_{FS})}^0
(\phi_{\sigma_k(t)})
\ee
exists. For the natural family of
$\d\sigma_k$ associated to a test
configuration,
this limit is in fact the Mumford numerical invariant, and we shall see subsequently that
\be
F={\rm lim}_{k\to\infty}{\mu_k\over k^n}
\ee
is the Donaldson-Futaki invariant, whose sign defines the notion of $K$-stability
(see Lemma \ref{FTmu} below).

\subsection{$F_{\o_0}^0$ and the Euler-Lagrange equation $R-\bar R=0$}

We show next that the equation $R-\bar R=0$ is, in a sense, the Euler-Lagrange
equation for $F_\o^0$ when restricted to ${\cal K}_k$.
For this, fix a metric $h$ on $L$ with positive curvature $\o$,
and let $\hat\us$ be an orthonormal
basis for $H^0(X,L^k)$ with respect to the $L^2$ metric defined by
$h$ and the volume form $\o^n$. By definition of the Kodaira imbedding,
$x_\al=s_\al(z)$, $V_{FS}=k^nV$,
and the matrix $M_{\bar\al\beta}$ of (\ref{M})
can be expressed as
\bea
\label{M1}
M_{\bar \al\beta}
=
{1\over V_{FS}}\int_X{\bar s_\al s_\beta\over
\sum_{\g=0}^{N_k}|s_\g|^2}(\iota_{\hat \us}^*(\o_{FS}))^n
=
{1\over V}\int_X{\<s_\beta,s_\al\>\over \rho_k(z)}\o^n\,{({1\over k}\iota_{\hat\us}^*(\o_{FS}))^n
\over\o^n}
\eea
where $\rho_k(z)$ is the density of states defined in (\ref{density}). In view of the Tian-Yau-Zelditch theorem
and Lu's formula described in section \S 4,
we can expand the right hand side in powers of $k$ and obtain
\be
\label{M2}
M_{\bar\al\beta}
=
\kappa\,\delta_{\bar\al\beta}-
{1\over Vk^{n+1}}\,\int_X \<s_\beta,s_\al\>(R-\bar R)\o^n + O({1\over k^{n+2}}).
\ee
for some constant $\kappa$ whose exact value is immaterial.
Thus we have, under a variation $\delta\sigma\in sl(N_k+1)$,
\bea
\label{FandR}
\delta F_{\iota_{\hat\us}^*(\o_{FS})}^0(\phi_\sigma)
&=&
{1\over V k^{n+1}}\,\int_X \delta\Phi\, (R-\bar R)\o^n+O({1\over k^{n+2}}),
\nonumber\\
\delta\Phi
&\equiv&{\rm Tr}((\delta\sigma+\delta\sigma^*)\cdot\<s_\beta,s_\al\>)
\eea
This observation is due to Donaldson \cite{D05}. It is in this sense that
$R-\bar R=0$ can also be viewed as the Euler-Lagrange equation for $F_{\o_0}^0$.

\subsection{$F_{\o_0}^0$ and Monge-Amp\`ere masses}

Another key property of $F_{\o_0}^0$ is its relation to
Monge-Amp\`ere masses, and hence to
infinite-dimensional GIT (see \S 12 below).
Strictly speaking, this is a property of $F_{\o_0}^0$
itself, not necessarily restricted to ${\cal K}_k$, but we include it in this section,
because the main application of this relation does involve restrictions to
${\cal K}_k$, and also because it fits in the theme of relations to stability
which are not apparent from the initial definition of $F_{\o_0}^0$.

\medskip
Let $(-T,0]\ni t\rightarrow \phi(z,t)$ be a smooth path in the space
${\cal K}$ of K\"ahler potentials. To such a path, we associate the following
potential $\Phi(z,w)$ defined on $X\times A$, with
$A=\{w\in{\bf C}; e^{-T}<|w|\leq 1\}$,
\be
\label{Phi}
\Phi(z,w)=\phi(z,\log |w|)
\ee
Let $\Omega_0$ be the curvature $\o_0$ of $h_0$, viewed as a
$(1,1)$-form on $X\times A$. Then we have the following basic identity \cite{PS06a}
\be
\label{MAmasses}
{1\over(n+1)V}\int\int_{X\times A}(\Omega_0+{i\over 2}\ddb \Phi)^{n+1}
=
-\dot F_{\o_0}^0(\phi(\cdot,0))+\dot F_{\o_0}^0(\phi(\cdot,T)).
\ee
In particular, for infinitely extended paths, we may take $T\to\infty$
and obtain, with $D^\times=\{w\in{\bf C};0<|w|\leq 1\}$,
\be
\label{MAmassesrays}
{1\over(n+1)V}\int\int_{X\times D^\times}(\Omega_0+{i\over 2}\ddb \Phi)^{n+1}
=
-\dot F_\o^0(\phi(\cdot,0))+{\rm lim}_{T\to\infty}\dot F_\o^0(\phi(\cdot,T)).
\ee
The identity (\ref{MAmasses}) can be established as follows.
A first important observation due to Semmes \cite{Se}
and Donaldson \cite{D99} is that
\be
\label{MAandgeodesic}
{1\over n+1}
(\Omega_0+{i\over 2}\ddb\Phi)^{n+1}=(\ddot\phi-g_\phi^{j\bar k}\pl_j\dot\phi\pl_{\bar k}\dot\phi)
\o_\phi^n \,dt\,d\theta
\ee
where $(e^t,\theta)$ are polar coordinates for $w$. On the other hand,
\be
-\pl_tF_\o^0(\phi(\cdot,t))
=\pl_t\bigg({1\over V}\int_X\dot\phi\,\o_\phi^n\bigg)
=
{1\over V}\int_X (\ddot\phi+\dot\phi\Delta\dot\phi)\o_\phi^n
=
{1\over V}\int_X(\ddot\phi-g_\phi^{j\bar k}\pl_j\dot\phi\pl_{\bar k}\dot\phi)\o_\phi^n.
\ee
Comparing the two identities and integrating in $t$ gives the desired formula.

\medskip
The significance of the identity (\ref{MAmasses}) is as follows. We shall
see later that the space ${\cal K}$ of K\"ahler potentials carries a natural
metric, with respect to which the geodesic equation for paths $t\to\phi(z,t)$ is precisely
the vanishing of the expression
$\ddot\phi-g_\phi^{j\bar k}\pl_j\dot\phi\pl_{\bar k}\dot\phi$, and hence
of $(\Omega_0+{i\over 2}\ddb\Phi)^{n+1}$
in view of the equation (\ref{MAandgeodesic}).
Just as ${\cal K}$ is the limit of the finite-dimensional symmetric spaces ${\cal K}_k$,
the geodesic rays in ${\cal K}$ are the analogues in infinite-dimensional GIT
of the one-parameter subgroups in ${\cal K}_k$, and thus the identity (\ref{MAmasses}) links this
key concept for stability once again with the $F_{\o_0}^0$ functional.

\newpage
\section{Notions of Stability}
\setcounter{equation}{0}

Stability is a condition on a geometric object that should insure that
the moduli space of stable objects be a well-behaved space, and in particular
Hausdorff. As  stressed earlier, it is  still unclear
what is the correct notion of stability for a positive line
bundle $L\to X$ that should apply in Yau's conjecture,
and be equivalent to the existence of constant scalar curvature
metric in the K\"ahler class $c_1(L)$. A prime candidate is
$K$-stability (or its variants), but other seemingly natural notions of stability
also arise, as suggested by infinite-dimensional GIT and the K\"ahler-Ricci
flow. We present here a brief overview of all these notions.

\subsection{Stability in GIT}

In GIT, one
typically associates to the geometric
object a non-zero ``defining vector'' $V$ in a  vector space $H$
carrying an action of $
SL(N+1,\C)$. The
vector $V$ is not intrinsically defined, as it depends on various
choices, but if $[V]\in \P H$ denotes
the image of $V$ in $\P H=(H\backslash \{0\})/\C^\times$, then the orbit
$SL(N+1)\cdot [V]\sub \P H$ is independent of all choices and it
characterizes the object.
We say that the orbit $SL(N+1)\cdot [V]\sub \P H$ is
GIT stable if  $SL(N+1)\cdot V\sub H$  is closed and the stabilizer of
$V$ is finite. Thus the orbit of $[V]$ characterizes the object
in question while the orbit of $V$ is the one which is relevant in
the definition of GIT stability.
It follows then that
the moduli space of stable orbits is Hausdorff and, in fact,
an algebraic variety \cite{Mu}. For an extensive discussion, see \cite{Th}.

\subsubsection{Closedness of orbits}

In our specific context, the geometric object is a positive line bundle $L\to X$
over a compact complex manifold. The GIT procedure outlined above can
be implemented as follows:

\medskip
As before, for $k$ sufficiently large, to each basis $\us=\{s_\al(z)\}_{\al=0}^{N_k}$
of $H^0(X,L^k)$ we can associate the image of $X$ by the corresponding
Kodaira imbedding $\iota_{\us}(X)$, which is a subvariety of ${\bf CP}^{N_k}$.
Clearly, the subvariety $\iota_{\us}(X)$ is not intrinsic,
since it depends on the choice of the basis $\us$. If $\ti\us$ is another basis, then
$\ti\us=\sigma\cdot\us$ and $\iota_{\ti\us}(X)=\sigma\cdot\iota_\us(X)$ for
some $\sigma\in GL(N_k+1)$.
As discussed above,
in order to define an intrinsic object in GIT
we associate to $X$ the whole orbit
\be
\label{orbit}
X\ \longleftrightarrow\ \{\sigma\cdot \iota_{\us}(X)\ |\ \sigma\in SL(N_k+1)\}.
\ee

As described above, we would like rather an orbit of vectors $V$
rather than an orbit of projective varieties. In Chow-Mumford stability,
this is achieved by associating a Chow point ${\rm Chow}(\hat X)$
to each n-dimensional subvariety $\hat X\subset {\bf CP}^N$ of degree $d$
as follows: let the Chow variety $Z$ of $\hat X$ be the
hypersurface in the Grassmannian $Gr(N-n-1,{\bf CP}^N)$ defined by
\be
\label{chowvariety}
Z=\{w\in Gr(N-n-1,{\bf CP}^N);\ w\cap \hat X\not=0\ \}
\ee
Let $H$ be the vector space
$H=H^0(Gr(N-n-1), O(d))$ and let
$[V_{\hat X}]$ be the set of
elements of $H$ which vanish on $Z$.
Then $[V_{\hat X}]$ is a line in $H$
(known as the ``Chow line") that
is,
$[V_{\hat X}]\in \P H$. We say that
$\hat X$ is Chow-Mumford stable if
the orbit $SL(N+1)\cdot [V_{\hat X}]\sub \P H$ is  GIT stable. We say that a polarized
variety $(X,L)$ is Chow-Mumford r-stable
if $L^r$ is very ample and if, for
some (and hence every) basis $\us$ of $H^0(X,L^r)$, the subvariety  $\iota_\us(X)\sub \P^{N_r}$ is 
Chow-Mumford stable in the sense defined
above.
\v
To define Hilbert-Mumford stability,
one associates a different line to each subvariety $\hat X$ of ${\bf CP}^N$
as follows. For $k>>0$, the inclusion
\be
\label{hilbertmumford}
\{s\in H^0({\bf CP}^N,O(k)) :\ s\big\vert_{\hat X}=0\ \}
\subset H^0({\bf CP}^N,O(k))
\ee
defines an element of the Grassmannian of codimension $p(k)$ planes
in $H^0({\bf CP}^N,O(k))$, where
\be
p(k)\equiv {\rm dim}\, H^0(\hat X,O(k)|_{\hat X})
\ee
is the Hilbert polynomial. The maximum wedge product of this subspace is
by definition the Hilbert-Mumford k-line.
We say $\hat X$ is Hilbert-Mumford
k-stable if the orbit of this line
is GIT stable. We say $\hat X$ is
Hilbert-Mumford stable if it is
k-stable for $k$ sufficiently large.
We say that $(L,X)$ is asymptotically
Hilbert-Mumford stable if for every
$r>>0$, there is a basis $\us$ of
$H^0(X,L^r)$ such that $\iota_\us(X)\sub\P^{N_r}$ is Hilbert-Mumford stable.

\medskip
More generally, let the Hilbert scheme ${\cal H}={\cal H}(N,p)$ be the
parameter space of subschemes
of ${\bf CP}^N$ with the same Hilbert polynomial $p(k)$. The preceding procedure
for associating an orbit of vectors $V_{\hat X}$ corresponds to a choice
of an ample line bundle $\eta\to {\cal H}$, carrying a lift of the action
of $SL(N+1)$. The orbit of vectors is then obtained by selecting a vector $V_{\hat X}$
in the fiber of $\eta$ over $\hat X$, and applying the $SL(N+1)$ action.

\begin{definition}
{\rm (a)} Let $\eta\to {\cal H}$ be an ample line bundle over the Hilbert scheme
equipped with a lifting of the
$SL(N+1)$ action to $\eta$.
Then a projective variety $\hat X$ is said to be $\eta$-stable
if the orbit $\{\sigma\cdot V_{\hat X}\}$ is closed, and the stabilizer
of $V_{\hat X}$ is finite.

\smallskip
{\rm (b)}
In the context of positive line bundles:
Let $L\to X$ be an  ample line bundle with Hilbert polynomial $p$. For $r>>1$, let ${\cal H}_r$ be the Hilbert scheme of subchemes
of $\P^{N_r}$ with Hilbert polynomial $p_r(k)=p(rk)$ and let $\eta_r\ra {\cal H}_r$
be a family of ample line bundles with
$SL(N_r+1)$ action. We say $(X,L)$ is  $\eta_r$-stable
if  the image $\hat X=\iota_{\hat\us}(X)
\subset {\bf CP}^{N_r}$ of $X$ by the
Kodaira imbedding is $\eta_r$ stable. \
We say $(X,L)$ is asymptotically $\eta$ stable if it is $\eta_r$ stable for all
sufficiently large $r$.
\end{definition}

In general, the notion of stability in GIT depends on the choice of ample line bundle
$\eta\to {\cal H}$. However, for Chow-Mumford and Hilbert-Mumford stability, we have
the following recent theorem due to Mabuchi \cite{Ma06},

\begin{theorem}
\label{mabuchitheorem}
Let $L\ra X$ be a positive line bundle over a compact complex manifold
$X$. Then $L\ra X$ is asymptotically Hilbert-Mumford stable if and only if
$L\ra X$ is asymptotically Chow-Mumford stable.
\end{theorem}

\subsubsection{The Hilbert-Mumford criterion}

There is a remarkable condition which is equivalent to stability in the sense of
GIT, but which can be formulated without reference to the closures of orbits
of vectors in ${\bf C}^M$ and is thus susceptible to generalization. This is the
Hilbert-Mumford numerical criterion, which can be stated as follows.

\medskip
Let $\hat X\in {\cal H}$, and let $\lambda:{\bf C}^\times\to SL(N+1)$
be a homomorphism. Set
\be
\label{centralfiber}
\hat X_0={\rm lim}_{\tau\to 0}\ \lambda(\tau)\cdot\hat X.
\ee
Then $\hat X_0$ is fixed by $\lambda$, and, denoting  $\lambda(\tau)$ as the lif of $\lambda$ to $\eta$,
we define the numerical invariant $F$ by
\be
\label{numericalinvariant}
\lambda(\tau)\xi\ = \ \tau^{-F}\xi
\ee
where $\xi$ is any non-zero vector in the fiber of $\eta$ above $\hat X_0\in{\cal H}$.
Then the Hilbert-Mumford numerical criterion says that $\hat X$ is $\eta$-stable
if and only if $F>0$ for all homomorphisms $\lambda: {\bf C}^\times\to SL(N+1)$.

\subsubsection{Chow-Mumford and K-stability}

The Hilbert-Mumford criterion is a powerful tool, since it reduces matters to
one-dimensio\-nal complex subgroups of $SL(N+1)$, and to the study of the sign
of a single numerical invariant rather than the property of closure of orbits. Even
more important for the problem of constant scalar curvature metrics,
it can be used to define a notion of $\eta$-stability even when
$\eta\to{\cal H}$ is a line bundle over the Hilbert scheme which is not
necessarily ample.

\medskip
Thus given a line bundle $\eta\to{\cal H}$, not necessarily ample, we define
$L\to X$ to be $\eta$-stable if the corresponding numerical invariant $F$
defined by (\ref{numericalinvariant}) is strictly positive for all
homomorphisms $\lambda: {\bf C}^\times\to SL(N+1)$.

\medskip
We can now describe the notions of
stability which arise in the context of metrics of constant scalar curvature.
Each corresponds to a choice of line bundle $\eta\to{\cal H}$. They can be
defined in several ways, but we shall give for each a formulation
in terms of Deligne intersection pairings (a brief summary of basic
facts about Deligne intersection
pairings can be found in section \S 6.1.4).
The advantage of the Deligne pairing formulation is
that it provides a unified approach,
and also ties in naturally with the $K$-energy and the functional
$F_{\o_0}^0$ which are central to the analytic formulation of the problem:

\medskip
$\bullet$ {\it Chow-Mumford stability}: We have given earlier
a geometric construction of the Chow-Mumford
bundle. In view of a theorem of S. Zhang \cite{Z96} (see also section \S 9.1),
it can also be defined as the following Deligne intersection pairing,
\be
\label{chowzhang}
\eta_{Chow}=\<O(1),\cdots, O(1)\>
\ee

$\bullet$ {\it K-stability}: The line bundle $\eta_K$ defining the notion of
$K$-stability can be defined in two equivalent ways,
either in terms of the line bundles $\lambda_j=\lambda_j(L,X,{\cal H})$
of the Knudsen-Mumford expansion (see section \S 6.1.5  below),
or in terms of Deligne intersection pairings:

\be
\label{etaK}
\eta_K=\cases{\lambda_{n+1}^\mu\otimes ({\lambda_{n+1}^2\over\lambda_n^2})^{n+1}\cr
\<K_X,O(1),\cdots, O(1)\>^{n+1}
\<O(1),\cdots,O(1)\>^{nc_1(X)c_1(L)^{n-1}\over c_1(L)^n}\cr}
\ee

This algebraic notion of $K$-stability is due to Donaldson \cite{D02},
and is very similar to
an earlier analytic notion of $K$-stability due to Tian \cite{T97, T02}, who defined it in
terms of a Futaki invariant for special degenerations with $Q$-Fano
(and in particular, normal) central fibers (see \S 7.1).
Both notions of $K$-stability have been conjectured to be equivalent to
the existence of a constant scalar curvature K\"ahler metric in $c_1(L)$,
the analytic version in \cite{T97}, and the algebraic version in \cite{D02}.

\smallskip

The equivalence between Donaldson's original formulation
and the formulation in terms of Knudsen-Mumford expansions
is in \cite{PT}. The equivalence with the above
Deligne intersection pairing formulation is in \cite{PRS}.

\smallskip
A slightly stronger version of $K$-stability has been proposed by Szekelyhidi \cite{Sze}.
Another version in terms of slopes, and hence more similar to the notion of Mumford-Takemoto
stability for vector bundles, has been introduced by Ross and Thomas \cite{RT}.

\subsubsection{Deligne pairings and energy functionals}

We provide here a summary of Deligne pairings and explain how they
are related to energy functionals \cite{Z96, PS03, SW}.

\medskip
Let $L_0,\cdots,L_n$ be holomorphic line bundles over a complex manifold
$X$ of dimension ${\rm dim}\,X=n$. Then the Deligne pairing $\<L_0,L_1,\cdots,L_n\>$
is a one-dimensional space, generated by the symbol $\<l_0,l_1,\cdots,l_n\>$,
where $l_i$ are generic meromorphic sections of $L_i$ in general position,
transforming in the following way under changes of sections $l_i\to l_i'$,
\be
\<l_0',l_1,\cdots,l_n\>
=
\bigg(\prod_{z\in \cap_{i\not=0}{\rm div}\,l_i}{l_0'(z)\over l_0(z)}\bigg)
\<l_0,l_1,\cdots,l_n\>.
\ee

More generally, if $\pi: X\to B$ is a flat projective morphism of integral schemes
of pure relative dimension $n$, and $L_0,L_1,\cdots,L_n$ are line bundles over $X$,
then we can define the Deligne pairing
\be
\<L_0,L_1,\cdots,L_n\>(X/B)
\ee
which is a line bundle on $B$, locally generated by symbols $\<l_0,l_1,\cdots,l_n\>$
with $l_j$ rational sections of $L_j$
in general position and transforming as above under changes
of sections $l_j$. A useful property is the induction formula
\be
\<L_0,\cdots,L_{n-1},L_n\>(X/B)
=
\<L_0,\cdots,L_{n-1}\>({\rm div}\,l_n/B),
\ee
if all components of ${\rm div}\,l_n$ are flat over $B$.

A key feature of Deligne pairings is that they are equipped with a natural
metric $\<h_0,\cdots, h_n\>$, if each line bundle $L_j$ comes
equipped with a metric $h_j$.
This metric satisfies the following property
\be
\<h_0e^{-\phi},h_1,\cdots,h_n\>=e^{-\psi}\,\<h_0,h_1,\cdots,h_n\>
\ee
where $\psi:B\to{\bf C}$ is the function
\be
\psi=\int_{X/B} \phi\cdot \prod_{j=1}^n\o_j(L_j)
\ee
and $\o_j(L_j)\equiv -{i\over 2}\ddb\log h_j$ is the curvature of $h_i$.
Denoting by $O(f)$ the trivial bundle with metric $\|1\|=e^{-f}$,
we obtain the following induction formula with metrics
\be
\label{inductionwithmetrics}
\<L_0,\cdots,L_{n-1},L_n\>(X/B)
=
\<L_0,\cdots,L_{n-1}\>({\rm div}\,l_n/B)
\otimes
O(-\int_{X/B}\log\|l_n\|\wedge_{j=0}^{n-1}\o_j(L_j))
\ee

The lines bundles $\eta_{Chow}$ and $\eta_K$ defining Chow-Mumford and
K-stability can now be related to the functionals $F_{\o_0}^0$ and $K_{\o_0}$ by
\cite{PS03, PS04b}

\begin{theorem}
\label{delignepairingstheorem}
Let $L\to X$ be a positive holomorphic line bundle, and let
$\eta_{Chow}$ and $\eta_K$ be given by the Deligne pairings
{\rm(\ref{chowzhang})} and {\rm (\ref{etaK})}.
For each metric $h_0$ on $L$, let $\o_0=-{i\over 2}\ddb \log h_0$ be the corresponding
metric on $X$, and
equip $\eta_{Chow}$ and $\eta_K$ with
the corresponding metrics as Deligne pairings.
Then under a change of metric $h\to h e^{-\phi}$, we have
\bea
\eta_{Chow}(h_0e^{-\phi}) &=& \eta_{Chow}(h_0)\otimes O\bigg((n+1)c_1(L)^nF_{\o_0}^0(\phi)\bigg)
\nonumber\\
\eta_{K}(h_0e^{-\phi}) &=& \eta_{K}(h_0)\otimes O\bigg(-(n+1)c_1(L)^n\,K_{\o_0}^0(\phi)\bigg)
\eea
\end{theorem}

These basic relations provide yet additional evidence that
the functionals $F_{\o_0}^0(\phi)$ and $K_{\o_0}(\phi)$
are closely related to the algebraic-geometric notion of stability.
Similarly, the Futaki functional defined by Futaki in \cite{Fu2}
can also be realized in terms of Deligne pairings, more specifically
$\<K_X^{-1},K_X^{-1},\cdots,K_X^{-1}\>$.
For the relation of Deligne pairings to other functionals in K\"ahler
geometry, see \cite{SW}. At the 2002 Complex Geometry conference in Tokyo,
Professors T. Mabuchi and L. Weng have informed us that they have also been aware
of potential applications of Deligne pairings to K\"ahler geometry for some time.

\subsubsection{Deligne pairings and Knudsen-Mumford expansions}

We discuss now the two formulations of $\eta_K$ given in (\ref{etaK})
and their equivalence.

\medskip
Let $\pi:X\to B$ be again a morphism of flat integral schemes
with constant relative dimension, and let $L\to X$ be a
relatively ample line bundle. The theorem of Knudsen-Mumford
\cite{KM} says that there exist functorially defined line bundles
$\lambda_j=\lambda_j(X,L,B)\to B$ such that
the following isomorphism holds for all $k>>1$
\be
\label{knudsenmumford}
{\rm det}\,\pi_*(L^k)
\sim
\lambda_{n+1}^{k\choose n+1}\otimes\lambda_n^{k\choose
n}\otimes\cdots\otimes\lambda_0.
\ee
Deligne \cite{De} showed that,
in dimension $n=1$, we have
$\lambda_2=\<L,L\>_{{\cal X}/B}$ and $\lambda_1^2=\<LK_X^{-1},L\>_{X/B}$.
More generally, for general $n$,
Knudsen-Mumford showed that $\lambda_{n+1}$ is equal to the Chow bundle.
Combined with the theorem of Zhang \cite{Z96},
which says that the Chow bundle is in turn equal to $\<L,\cdots,L\>_{
X/B}$ (see Theorem \ref{zhangtheorem} below),
we have then
\be
\lambda_{n+1}(L,X,B)=\<L,\cdots,L\>_{X/B}.
\ee

The bundle $\eta_K$ in terms of the Knudsen-Mumford line bundles $\lambda_j$
in the first line of the definition (\ref{etaK}) was given by Paul-Tian \cite{PT},
who pointed out that
its numerical invariants coincided with the Donaldson-Futaki invariants
defined by Donaldson \cite{D02} in his definition of $K$-stability.
The bundle $\eta_K$ in terms of Deligne pairings in the second line
of the definition (\ref{etaK}) was first introduced in
\cite{PS03, PS04b}.
The following theorem due to \cite{PRS} identifies the second leading
term
$\lambda_n$ in the Knudsen-Mumford expansion,
and shows that the two bundles
in the definition of $\eta_K$ in
(\ref{etaK}) are equivalent:

\begin{theorem}
Let $\pi:X\to B$ be a proper flat morphism of integral
schemes of relative dimension $n\geq 0$, and let $L\to X$
which is very ample on the fibers. Assume that $X$ and $B$
are smooth. Let $K$
be the relative canonical line bundle of $X\to B$.
Then

\noindent
{\rm (i)} There is a canonical functorial isomorphism
\be
\lambda_n^2(L,X,B)
=
\<L^n K^{-1},\cdots,L\>_{X/B}.
\ee
{\rm (ii)} In particular,
\be
\lambda_{n+1}^\mu\otimes ({\lambda_{n+1}^2\over\lambda_n^2})^{n+1}
=
\<K_X,O(1),\cdots, O(1)\>^{n+1}
\<O(1),\cdots,O(1)\>^{n\,c_1(X)\,c_1(L)^{n-1}\over c_1(L)^n}
\ee
\end{theorem}

\subsubsection{Test configurations and Donaldson-Futaki invariants}

As we have seen in \S 6.1.3,
the definition of $\eta$-stability, when $\eta$ is not ample,
consists solely in the requirement that the corresponding
numerical invariant $F$ be strictly positive.
Thus it is useful to exhibit this numerical invariant more
concretely in the case of $K$-stability. For this as well as for
other subsequent uses, it is preferable to use the formalism of
``test configurations'' introduced by Donaldson \cite{D02}
rather than of one-parameter subgroups.
In view of its importance,
we describe this formalism in some detail in this section.

\begin{definition}
\label{testconfiguration}
Let $L\to X$ be a positive line bundle over a compact complex manifold $X$.
A test configuration ${\cal T}$ consists of

{\rm (1)} a scheme ${\cal X}$ with a ${\bf C}^\times$
action $\rho$,

{\rm (2)} an ${\bf C}^\times$
equivariant line bundle ${\cal L}\to{\cal X}$, ample on all fibers,

{\rm (3)} and a flat ${\bf C}^\times$ equivariant map $\pi:{\cal X}\to{\bf
C}$
where ${\bf C}^\times$ acts on ${\bf C}$ by multiplication,

\noindent
with the following
properties: the pair $(X_1,L_1)$ is isomorphic to $(X,L^r)$
for some $r>0$
 and $L^r$ is
very ample.
Here $X_w=\pi^{-1}(w)$ and $L_w={\cal L}_{\vert_{X_w}}$.
\end{definition}

It is convenient to denote the test configuration
${\cal T}$ by
\be
{\cal T}=\bigg(\rho:{\bf C}^\times\to {\rm Aut}({\cal L}\to{\cal X}\to{\bf C})\bigg).
\ee
We also say that a test configuration is trivial
if ${\cal L}=L\times {\bf C}$ with the trivial
action $\rho(\tau)(l,w)=(l,\tau w)$, for $(l,w)\in L\times {\bf C}$, $\tau\in{\bf C}^\times$.

\medskip

An important property of a test configuration $\cT$ is that it
defines for each $k$ an $(N_k+1)\times(N_k+1)$ diagonal matrix $B_k$ and
its traceless part $A_k$. Indeed,
in a test configuration
$\cT$, the fiber $L_0\to X_0$ is fixed by the ${\bf C}^\times$ action
$\rho$. Thus we obtain
an induced automorphism $\r_k(\tau)$ of $H^0(X_0,L_0^k)$
\be
\label{rhok}
\r_k(\tau)\ :
\  H^0(X_0,L_0^k)\ \longrightarrow\ H^0(X_0,L_0^k)
\ee
and hence, if
we choose a basis $\us$ of $H^0(X_0,L^k_0)$, we obtain a
one-parameter subgroup
$\r:\C^\times\ra GL(N_k+1)$. Since every one-parameter subgroup
is a direct sum of one dimensional subgroups, there exists a
 basis $\us$
with the following property. The matrix $\r(\tau)$ is  diagonal  with
entries $\tau^{\l_\al}$ where $\l_0^{(k)}\leq \cdots \leq \l_{N_k}^{(k)}$
are integers that are independent of $\tau$. We let $B_k$ be
the diagonal matrix with entries $\l_\al^{(k)}$ and we shall often
write
\be
\label{Bk}
\rho(\tau)\  = \ \tau^{B_k}
\ee
The matrix $A_k$ is then defined to be the traceless part of $B_k$.

\medskip
We shall see shortly (c.f. Lemma \ref{tc} below)
that a test configuration can be imbedded equivariantly into
${\bf CP}^{N_k}$, with $\tau^{B_k}$
defining a one-parameter subgroup completed by the central
fiber $X_0$ (\ref{centralfiber}). Thus,
for each choice of bundle $\eta$ over the
Hilbert scheme, as in the case of one-parameter subgroups
described in section \S 6.1.2,
we can define a numerical invariant.
Since we had not given any detail there of the construction
of the numerical invariant, we do so now,
in the equivalent context of test configurations.

\medskip
Let as before $p(k)\equiv {\rm dim}\,H^0(X,L^k)$
be the Hilbert polynomial of $L\to X$,
and let ${\cal H}$ be the Hilbert scheme of all subvarieties
of ${\bf CP}^{p(k)-1}$ with Hilbert polynomial $p(k)$.
Let $\eta\to {\cal H}$ be a line bundle on ${\cal H}$,
not necessarily ample, which carries a linearization of $GL(p(k),{\bf C})$,
that is, an action of $GL(p(k),{\bf C})$ which is a lift of the natural
action of $GL(p(k),{\bf C})$ on ${\cal H}$.

\smallskip

Let $\cT=(\rho:{\bf C}^\times\to {\rm Aut}({\cal L}\to{\cal X}\to{\bf C})$
be a test configuration. Then we have seen that $\cT$ induces an
endomorphism $\rho_k(\tau)$ on $H^0(X_0,L_0^k)$ defined by
\be
(\r_k(\tau)(s))(x)\ = \ \r(\tau)^{-1}(s(\r(\tau)x))\ \ \hbox{for all $x\in X_0$}
\ee
so if we fix a basis $\us=(s_0,...,s_{p(k)-1})$ then
$\r(\tau)(s_\al)=\s\sigma(\tau)_{\al\b}s_\b$ for some invertible
matrix $\sigma(\tau)=(\sigma(\tau)_{\al\b})\in GL(p(k))$.
Moreover, $\sigma:\C^\times\ra GL(p(k))$ is a one parameter subgroup.
For $k>>1$ we let $\iota_\us:X_0\hookrightarrow\C\P^{p(k)-1}$ be the imbedding
$x\mapsto [s_0(x),...,s_{p(k)-1}(x)]$.
Then one easily sees that $\iota_\us(X_0)$ is invariant under
the action of $\sigma(\tau)$. Indeed,  if $x\in X_0$ then
$\r(\tau)x\in X_0$ and
\be
\sigma(\tau)(\iota_\us(x))\ = \ (\iota_\us(\r(\tau)x))
\ee
To see this we observe
\bea
\sigma(\tau)(\iota_\us(x))\
&=& \ (\r_k(\tau)(s_\al))(x)\ = \
(\r(\tau)^{-1}(s_\al(\r(\tau)x))
\nonumber\\
&=&\ {s_\al(\r(\tau)x)\over s_{\al_0}(\r(\tau)x) }
\r(\tau)^{-1}(s_{\al_0}(\r(\tau)x))\ = \ (s_\al(\r(\tau)x))
\eea

Now let $[\iota_\us]\in\cH$ be the point in the Hilbert scheme
corresponding to the imbedding $\iota_\us:X_0\hookrightarrow\C\P^{p(k)-1}$.
The discussion above shows that $[\iota_\us]$ is a fixed point for the action
of $\sigma(\tau)\in GL(p(k))$. Since the $GL(p(k),\C)$ action on the Hilbert
scheme $\cH$ lifts to an action on $\eta$, it follows that there
exists $F\in \Z$ such that
\be
\label{donaldsonfutakiinvariant}
\sigma(\tau)n_\us\ = \ \tau^{-F}n_\us
\ee
for all $n_\us$ in the fiber of $\eta$ above $[\iota_\us]\in\cH$.
\v
We claim that  $F$ is independent of the basis $\us$. To see this, let $\ti\us$
be another basis and $\g\in GL(p(k))$ chosen so that
$\ti\us=\g\us$. Let  $\ti\sigma(\tau)$ be the matrix representing
$\r_k(\tau)$ with respect to the basis $\ti\us$.
Then $\ti\sigma(\tau)=\g\sigma(\tau)\g^{-1}$. Moreover
$\g\iota_\us(x)=\iota_{\ti\us}(x)$ so $\g[\iota_\us]=[\iota_{\ti\us}]$.
Fix $n_{\us}$ as above and define $n_{\ti\us}=\g n_\us$. Then
\be
\ti\sigma(\tau)n_{\ti\us}\ = \ \g\sigma(\tau)\g^{-1}\g n_\us\ =
\ \g\tau^{-F}n_\us\ = \ \tau^{-F}n_{\ti\us}
\ee
establishing the desired invariance.

\medskip
The bundle $L\to X$ is consequently defined to be $\eta$-stable if
$F>0$ for all test configurations. In view of the correspondence between
test-configurations and one-parameter subgroups given below,
this definition is just a rephrasing
in terms of test configurations
of the earlier definition of
$\eta$-stability. However, it allows us to make now contact,
in the case of $\eta=\eta_K$, with the original definition of K-stability
by Donaldson \cite{D02}:

\begin{definition}
Let $X_0$ be a projective scheme and $\al:\C^\times\ra\Aut(L_0\ra X_0)$
an algebraic homomorphism. Then  $F(\al)$ is defined by the asymptotic expansion
\footnote{Our conventions for $F(\al)$ differ from those of \cite{D02} by a minus sign.}

\be
{{\rm Trace}\,B_k\over k\, (N_k+1)}=F_0-F({\al})k^{-1}+ O(k^{-2})
\ee
where $B_k: H^0(X_0,L_0^k)\to H^0(X_0,L_0^k)$
is the infinitesimal generator  of the one-parameter group of endomorphisms
induced by the action of $\al$.
\v
\noindent
Now let $L\to X$ be an ample line bundle.
Let ${\cal T}$ be a test configuration
for $L\to X$.
Then the Donaldson-Futaki invariant $F({\cal T})$
is defined by $F({\cal T})=F(\tilde\r)$ where $\ti\r$ is
the restriction of $\r$ to the central fiber.
\end{definition}

In \cite{D02}, Donaldson defined $L\to X$ to be $K$-stable if
$F({\cal T})>0$ for any non-trivial test configuration ${\cal T}$.
As we had mentioned earlier, it was noted by Paul-Tian \cite{PT}
that $F({\cal T})$ coincided with the numerical invariant $F$ of the
bundle $\eta_K$. Thus the definition of $K$-stability of \cite{D02}
coincides with $\eta_K$-stability as defined in section \S 6.1.3.

\subsubsection{Equivariant imbeddings of test configurations}

The correspondence between one-parameter subgroups
and test configurations is described by the following lemma,
which provides an equivariant imbedding of test configurations
into projective space (see
the original statement in \cite{D05}, and \cite{PS07} for the version
presented here):

\begin{lemma}
\label{tc}
{\rm (i)} Let $X$ be a subvariety of ${\bf CP}^{N_k}$ and let
$\l:\C^\times\ra GL(N_k+1)$ be a one-parameter subgroup of $GL(N_k+1)$.
Let
\be
\cX^\times=\{\,(\l(\tau)x,\tau)\,: \tau\in\C^\times, x\in X\,\}\sub {\bf CP}^{N_k}\times\C^\times.
\ee
Then the closure (or, more precisely, the flat limit) of $O(1)\times\C^\times\to \cX^\times\to {\bf
C}^\times$ inside $O(1)\times \C\ra \P^{N_k}\times \C\ra \C$, is a test configuration.

\smallskip

{\rm (ii)} Conversely, let $L\to X$ be a positive line bundle,
and  $\rho:{\bf C}^\times\to {\rm Aut}({\cal L}\to{\cal X}\to{\bf C})$
be a test configuration.  Let $k$ be an integer such that $L^k$
is very ample. Let $B_k$ be the matrix defined in {\rm (\ref{Bk})}.
Then there exists
a basis $\underline s$ of $H^0(X,L^k)$,
and an imbedding
\be
I_{\underline s}:
({\cal L}^k\to{\cal X}\to{\bf C})
\to
(O(1)\times{\bf C}\to {\bf CP}^{N_k}\times{\bf C}\to{\bf C})
\ee
which restricts to $\iota_{\underline s}$ on the fiber $X_1$
and intertwines $\rho(\tau)$ and $B_k$, i.e.,
\be
\label{intertwining}
I_{\underline s}(\rho(\tau)l_w)
=
(\tau^{B_k}\cdot I_{\underline s}(l_w),\tau w)\ \ \
{\rm for\ each}\ l_w\in L^k
\ee

{\rm (iii)} If $h$ is a metric on $L$
with positive curvature $\o=-{i\over 2}\ddb\log h$, then we may choose
$\us$ to be an orthonormal basis of $H^0(X,L^k)$
with respect to the $L^2$ metric on $H^0(X,L^k)$ induced
by $h$ and $\o$. Moreover, the basis $\us$
is unique up to $U(N_k+1)$ elements commuting with $B_k$.
\end{lemma}

\subsubsection{Bott-Chern secondary characteristic classes}

We have seen how the energy functionals $F_\o^0(\phi)$ and $K_\o(\phi)$
are readily related to Deligne intersection pairings, and whence to
Chow-Mumford and $K$-stability. There is an alternative construction
of energy functionals, based on Bott-Chern secondary characteristic classes,
which is also of interest. It began with Donaldson's \cite{D87} construction of
the functional for Hermitian-Einstein metrics, and has been developed further
by Tian \cite{T00} in a more general context. We give a brief exposition
following \cite{T00}.

\medskip

Let $E\to X$ be a holomorphic vector bundle of rank $r$ over a compact
K\"ahler manifold $X$. For each function
$\Phi(A)$ on $gl(r,{\bf C})$ which is a symmetric polynomial
of degree $k$ in the eigenvalues of $A$,
we have a Chern-Weil $(k,k)$-form
$\Phi({i\over 2\pi}F(H))$, defined
for each Hermitian metric $H_{\bar\al\beta}$ on $E$, where
$F=-\bar\pl(H^{-1}\pl H)$
is the curvature of $H$. These forms are
cohomologically equivalent for different metrics $H_{\bar\al\beta}$.
In fact, we can write explicitly
\be
\Phi({i\over 2\pi}F(H))-\Phi({i\over 2\pi}F(H'))
=-\ddb {\bf BC}(\Phi;H,H'),
\ee
where ${\bf BC}$ is a form
given explicitly modulo ${\rm Im}\,\pl+{\rm Im}\,\bar\pl$ by
\be
{\bf BC}(\Phi;H,H')
=\int_0^1\Phi'({i\over 2\pi}F(H_t),H_t^{-1}\dot H_t)dt.
\ee
Here $H_t$ is a smooth one-parameter family of metrics joining $H=H_0$ to
$H'=H_1$, and the derivative $\Phi'$ of $\Phi$ is defined by
$\Phi'(A;B)={d\over dt}\vert{_{t=0}}(Ae^{tB})$. The form ${\bf BC}$ is called
the Bott-Chern secondary characteristic class associated to $\Phi$.
The corresponding Donaldson functional is then defined by
\be
\label{donaldsonfunctional}
D(\Phi;H,H')=\int_X {\bf BC}(\Phi;H,H').
\ee
To make connection with the functionals $F_{\o_0}^0(\phi)$ and
$K_{\o_0}(\phi)$ introduced in the previous sections,
let $L$ be now a positive line bundle over $X$
with metric $h_0$ and positive curvature $\o_0$.
If we choose $\Phi(A)={\rm Trace}(A^{n+1})$,
then for any other metric $e^{-\phi}h_0$ with positive curvature
$\o_\phi=\o_0+{i\over 2}\ddb\phi$,
we have
\be
F_{\o_0}^0(\phi)=D(\Phi;\o_0,\o_\phi).
\ee
Similarly, $K_{\o_0}(\phi)$ coincides with $D(\Phi;H,H')$,
when $E$ is the virtual bundle defined by
\be
E=(n+1)(K_X^{-1}-K_X)\otimes (L-L^{-1})^n-\mu(L-L^{-1})^{n+1}.
\ee

\subsection{Donaldson's infinite-dimensional GIT}

Beginning in the late 1990's, Donaldson \cite{D99} developed an approach
to the problem of constant scalar curvature K\"ahler metrics, motivated
partly
by the interpretation of the scalar curvature as the moment map of an
infinite-dimensional group of symplectic automorphisms.
Central to this approach is a natural interpretation of the space
${\cal K}$ of K\"ahler potentials\footnote{We denote
the space of K\"ahler potentials and the space of corresponding
K\"ahler metrics by the same letter ${\cal K}$, although
strictly speaking, the latter is equal to the former
modulo constants.} in $c_1(L)$,
\be
{\cal K}=\{\phi\in C^\infty(X); \ \o_0+{i\over 2}\ddb\phi>0\ \}
\ee
as an infinite-dimensional
symmetric space, whose geodesics would play a similar role
to one-parameter subgroups in the usual finite-dimensional
GIT.

\medskip
The starting point is that ${\cal K}$ admits a natural Riemannian metric,
given at
$\phi$ by
\be
\|\delta\phi\|_\phi^2=\int_X |\delta\phi|^2\o_\phi^n
\ee
where we have identified the tangent space to ${\cal K}$ at $\phi$
with $C^\infty(X)$. The following theorem is then due to Donaldson
\cite{D99},
Semmes \cite{Se}, and Mabuchi \cite{Ma87}:

\medskip
\begin{theorem}
The Riemannian manifold ${\cal K}$ is an infinite-dimensional symmetric
space
of non-positive curvature, in the sense that its curvature is covariant
constant
and non-positive. In fact, the Riemann curvature tensor $R_\phi$
and the sectional curvature $K_\phi$ at the point $\phi\in{\cal K}$
are given by the following expression
\be
R_\phi(\delta_1\phi,\delta_2\phi)\delta_3\phi
=
-{1\over 4}\{\{\delta_1\phi,\delta_2\phi\}_\phi,\delta_3\phi\}_\phi,
\quad
K_\phi(\delta_1\phi,\delta_2\phi)
=
-
{1\over 4}\|\{\delta_1\phi,\delta_2\phi\}_\phi\|_\phi^2,
\ee
where $\{\ ,\ \}_\phi$
is the Poisson bracket on $C^\infty(X)$ defined by the symplectic form
$\o_\phi$.
\end{theorem}

It is also shown in \cite{Calabichen} that the space ${\cal K}$ is non-positive
in the sense of Alexandrov, and essentially that the Calabi flow (see \S 10.1)
is distance decreasing in this space.

\subsubsection{Geodesic segments and geodesic rays}

The geodesics $(-T,0]\ni t\to\phi(t)$ in ${\cal K}$ are by definition
the solutions of the Euler-Lagrange equation of the length functional.
They are readily seen to be given by
\be
\label{geodesicequation}
\ddot\phi-g_\phi^{j\bar k}\pl_j\dot\phi\pl_{\bar k}\dot\phi=0,
\ee
where we have denoted by $(g_\phi)_{\bar kj}=g_{\bar kj}^0+\pl_j\pl_{\bar
k}\phi$
the metric corresponding to the K\"ahler form $\o_\phi$.
The geodesic equation admits an interpretation
as a completely degenerate complex Monge-Amp\`ere equation which will be
crucial in subsequent developments: let $A=\{w\in {\bf C}; e^{-T}\leq
|w|\leq 1\}$
and define the $S^1$ invariant function $\Phi(z,w)$ on $X\times A$
as in (\ref{Phi}) by $\Phi(z,w)=\phi(z,\log\,|w|)$.
Let $\Omega_0$ be again the K\"ahler form $\o_0$
viewed as a form on
$X\times A$.
Then in view of the equation
(\ref{MAandgeodesic}), we see that the above geodesic equation
is equivalent to the following
completely degenerate Monge-Amp\`ere equation \cite{Se, D99}
\be
\label{cdma}
(\Omega_0+{i\over 2}\ddb\Phi)^{n+1}=0\ \ {\rm on}\ X\times A.
\ee
Geodesic segments with given end points correspond to
$S^1$ invariant solutions of the degenerate complex Monge-Amp\`ere
equation
with given Dirichlet boundary conditions on $X\times \pl A$.
We shall also be particularly interested in geodesic rays,
which would be defined for $t\in (-\infty,0]$, and which correspond
to the same set-up with $A$ replaced by the punctured disk
$D^\times=\{w\in{\bf C};0<|w|\leq 1\}$. In the case of a
given starting point for the geodesic ray,
the Dirichlet boundary condition would be specified at $X\times
\{|w|=1\}$.
The formulation in terms of the complex degenerate Monge-Amp\`ere equation
also provides us with a natural notion of generalized geodesics: these will
be by definition $S^1$ invariant, locally bounded, generalized solutions
of the Monge-Amp\`ere equation (\ref{cdma}) in the sense of pluripotential
theory (see \S 12.3.1).

\subsubsection{Stability conditions in terms of geodesic rays}

We have seen that, by the Tian-Yau-Zelditch theorem,
the space ${\cal K}$ of K\"ahler potentials can be viewed as the limit
of the spaces ${\cal K}_k$ as $k\to\infty$.
Now stability conditions are conditions on the
numerical invariants of one-parameter orbits of ${\cal K}_k$.
Since ${\cal K}_k$ is a symmetric space, these are the same as geodesics,
and it is then natural to construct a GIT theory based directly on geodesic rays
in the infinite-dimensional symmetric space ${\cal K}$.

\medskip
The cost of working in an infinite-dimensional setting
is partly offset by the intuition we gain from the theory of symmetric
spaces
and moment maps \cite{D99}. There is also a very powerful incentive
for the use of geodesic rays in ${\cal K}$: it is the fact that
the $K$-energy is convex along these rays. In fact,
if $\phi(t)$ is any smooth path in ${\cal K}$, then we have
\be
\label{hessianofK}
\ddot K_\o(\phi(t))
=
{1\over V}\int_X |\bar\nabla\bar\nabla\dot\phi|^2\o_\phi^n-
{1\over V}\int(\ddot\phi-g_\phi^{j\bar k}\pl_j\dot\phi\pl_{\bar
k}\dot\phi)(R-\mu n)\o_\phi^n.
\ee
This implies that $K_\o(\phi(t))$ is convex when $\phi(t)$ is a geodesic,
and in fact, strictly convex if ${\rm Aut}^0(X)=0$. This strict convexity
has some immediate consequences: if geodesics always exist,
the critical points, and thus the metrics of constant scalar curvature,
would be unique.

\smallskip
The analogy with stability in GIT leads to the following
early conjecture/question of Donaldson \cite{D99}: the non-existence of constant scalar
curvature
K\"ahler metrics in $c_1(L)$ should be equivalent to the existence of an
infinite geodesic ray $\phi_t$, $t\in (-\infty,0]$ such that
\be
{1\over V}\int_X \dot\phi\,(R-\mu n)\,\o_\phi^n <0
\ee
for all $t\in (-\infty,0]$. Since the left-hand side is just $-\dot
K_\o(\phi)$
at $\phi_t$, this condition is very close to the formulation of stability
in terms of the asymptotic behavior of energy functionals, and can be
viewed
as an infinite-dimensional version of $K$-stability.

\subsection{Stability conditions on $Diff(X)$ orbits}

There is another concept of stability which arises naturally if we try to construct
the moduli space of complex structures on a smooth manifold $X$ as
the space of equivalence classes of integrable almost-complex structures $J$
on $X$, modulo diffeomorphisms
\be
\{J\ {\rm integrable\ almost-complex\ structure}\,\}/Diff(X).
\ee
A point in this moduli space would then be an orbit $Diff(X)\cdot J$
for some $J$, and the orbit would have to be closed in order for
the point to be a Hausdorff point. From this point of view,
the closedness of orbits is a condition similar to the ones in stability in GIT,
except that the orbits are now that of the infinite-dimensional group $Diff(X)$
rather than the finite-dimensional groups $SL(N+1)$. We describe now some specific
versions of such stability conditions, motivated as we shall see by the convergence
of the K\"ahler-Ricci flow (\cite{PS06, PSSW2}, and section \S 8.2).

\subsubsection{Condition (B)}

Let $X$ be a compact K\"ahler manifold and
let $J$ be
the almost-complex structure of
$X$, viewed as a tensor.
Condition (B) is the following condition \cite{PS06}:

\medskip

(B): There exists no almost-complex structure $\tilde J$ in the $C^\infty$ closure
of the orbit $Diff(X)\cdot J$ with
\be
\label{B}
{\rm dim}\,H^0(X,T_{\tilde J}^{1,0})
>
{\rm dim}\,H^0(X,T_J^{1,0})
\ee
Here $H^0(X,T_{\tilde J}^{1,0})$ denotes the space of holomorphic $(1,0)$
vector fields with respect to the integrable almost-complex structure $\tilde J$.

\medskip
Clearly, if Condition (B) is violated, the orbit $Diff(X)\cdot J$ cannot be closed,
and thus does not define a Hausdorff point.

\medskip
For our purposes (see section \S 8.2.4), we actually do not need the full strength of
condition (B), but actually only the following weaker version which may be closer
in spirit to the other stability conditions formulated in this chapter:

\medskip
${\rm (B^*)}$: Let
$g_m$ be any sequence of metrics in $c_1(X)$
with $K$-energy bounded from above, and which converges in $C^\infty$ in the sense
of Cheeger-Gromov, i.e., there is a family of diffeomorphisms
$F_m:X\to X$ so that $F_m^*(g_m)$ converges in $C^\infty$.
Then the sequence of almost-complex structures
$F_m^*(J)$ does not admit any $C^\infty$ limit point $J_\infty$
satisfying (\ref{B}).

\subsubsection{Condition (S)}

A second type of stability condition that arises in the study of the
K\"ahler-Ricci flow is a condition
called there Condition (S). Let $X$ be a
compact Fano manifold and let $J$ be again
the almost-complex structure of
$X$, viewed as a tensor.
Condition (S) is the following condition  \cite{PSSW2, PSSW3}:

\medskip

(S): There exists some solution $\o(t)$ to the K\"ahler-Ricci flow
$\dot g_{\bar kj}(t)=-(R_{\bar kj}-g_{\bar kj})$ so that
\be
\lambda\equiv {\rm inf}_{t\geq 0}\,\lambda_{\omega(t)} \,>\,0,
\ee
\noindent
where $\lambda_\o$ is the lowest strictly positive
eigenvalue of the Laplacian $\bar\partial^\dagger\bar\partial=-g^{j\bar
k}\nabla_j\nabla_{\bar k}$ acting on smooth $T^{1,0}(X)$ vector fields,
\be
\label{lambda}
\lambda_\o={\rm inf}_{V\perp H^0(X,T^{1,0})}{\|\bar\partial V\|^2\over \|V\|^2}
\ee
\v

Condition (S)  can be interpreted as a stability
condition in the following sense: Assume that there exist
diffeomorphisms $F_t:X\ra X$
so that $(F_t)^*(g(t))$ converges in $C^\infty$ to a metric
$\tilde g(\infty)$. Then if $J$ is the
complex structure of $X$, the pull-backs $(F_t)_*(J)$ converge also to a
complex structure $J(\infty)$ (see \cite{PS06}, \S 4). Clearly, the eigenvalues
$\lambda_{\o(t)}$ are unchanged under $F_t$.
If they do not remain bounded away
from $0$ as $t\to\infty$, then the complex structure $J(\infty)$ would
have a strictly higher
number of independent vector fields than $J$. Thus the $C^\infty$ closure of the orbit
of $J$ under the diffeomorphism group contains a complex structure different from
$J$, and $J$ cannot be included in a Hausdorff moduli space of complex structures.

\v
Just as the condition (B) can be weakened to condition ${\rm (B^*)}$,
we can strengthen condition (S) to a condition ${\rm (S^*)}$ which requires
more than we need, but which has now the advantage of not referring
specifically to the K\"ahler-Ricci flow:

\v
${\rm (S^*)}$ Let $g_m$ be any sequence of metrics in $c_1(X)$
with $K$-energy bounded from above. Then the corresponding eigenvalues
$\lambda_{g_m}$ are bounded uniformly from below by a strictly positive constant.

\newpage
\section{The Necessity of Stability}
\setcounter{equation}{0}

After describing at one end the analytic problem and at the other end
the expected algebraic-geometric answer, we turn now to describing
some of the partial results linking the two which have been obtained so far.
The more complete results go in the direction of necessity, that is,
the existence of constant scalar curvature metrics implies various
forms of stability or semi-stability. We describe some of these
results in this section.

\subsection{The Moser-Trudinger inequality and analytic K-stability}

Consider first the case of $L=K_X^{-1}$ and thus K\"ahler-Einstein metrics.
In this case, we have the following theorem of Tian \cite{T97},
which shows the equivalence between the existence of a K\"ahler-Einstein metric
and a properness property of the $F_{\o_0}$ functional, and which led to the
notion of $K$-stability (see Theorem \ref{Kstabilitynecessary} below):

\begin{theorem}
\label{mosertrudingertheorem}
Let $(X,\omega_0)$ be a compact manifold with $[\o_0]=c_1(X)>0$. If ${\rm Aut}^0(X)=0$,
then $X$ admits a K\"ahler-Einstein
metric if and only if there exists $\gamma>0$ and constants $A_\gamma>0$, $B_\gamma$ so that
\be
\label{mosertrudinger}
F_{\omega_0}(\phi)\geq A_\gamma\,\,J_{\omega_0}(\phi)^\gamma- B_\gamma
\ee
for all $\phi$ which are $\o_0$-plurisubharmonic. 
\end{theorem}

It was conjectured in \cite{T97} that the exponent $\gamma$ in
the inequality (\ref{mosertrudinger}) can be taken to be $\gamma=1$.
This conjecture, as well as the following extension to the case
${\rm Aut}^0(X)\not=0$,
was established in \cite{PSSW1}:

\begin{theorem}
Let $(X,\omega_0)$ be a compact manifold with $[\o_0]=c_1(X)>0$.
Let $G\subset {\rm Aut}^0(X)$
be any closed subgroup whose centralizer in the stabilizer ${\rm
Stab}(\omega_{KE})$
is finite. Then $X$ admits a $G$-invariant K\"ahler-Einstein metric
if and only if there exists constants $A>0$, $B$,
so that the following inequality holds for $G$-invariant
$\o_0$-plurisubharmonic potentials $\phi$,
\be
\label{mosertrudinger1}
F_{\omega_0}(\phi)\geq A\,\,J_{\omega_0}(\phi)- B
\ee
\end{theorem}

\medskip
As a consequence, Tian \cite{T97} deduced the following theorem,
which implies in particular that the Iskovskih manifolds
are smooth manifolds with $c_1(X)>0$, no holomorphic vector fields
(and hence automatically vanishing Futaki invariant), and yet no
K\"ahler-Einstein metrics. Define a special degeneration to be a
one-parameter
test configuration where the central fiber $X_0$ is a Q-Fano variety.
This means that $X_0$ is
a normal variety, and if
$K^{-1}$ is the anti-canonical line bundle on $X_0^{\rm reg}\sub X_0$ (the smooth
points of $X_0$) then there exists  $L$,
an ample line bundle on $X_0$ whose restriction to $X_0^{\rm reg}$ is
$K^{-1}$. In this set-up, by using a resolution of singularities,
Ding and Tian \cite{DT} showed that one can define a generalized
Futaki invariant $f_{X_0}$. This invariant can also be constructed explicitly
as follows.

\smallskip

Let $\phi_{L^k}:X_0\ra\P^N$ be a Kodaira imbedding corresponding
to some basis $s_j$ of $H^0(X_0,L^k)$, let  $h=\phi^*_{L^k}h_{FS}^{1/k}$,
which is a metric on the line bundle $K^{-1}$, and let $\o =
{1\over k}\phi^*_{L^k}\o_{FS}$, which is a K\"ahler metric on
$X_0^{\rm reg}$. Then $\o^n$ defines a second metric on $K^{-1}$. Two
metrics on the same holomorphic line bundle differ by a positive
function: $\o^n = he^{-f/k}$. Thus $Ric(\o)-\o={i\over 2}\ddb f$. In fact, we
have an explicit formula for $f$:
\be
f \ = \ -\log{ \s |s_j|^2_\o}
\ee
Let  $V$ be the holomorphic vector  field induced by the
${\bf C}^\times$ action on $X_0^{\rm reg}$. Then
 $f_{X_0}=\I_{M_{\rm reg}} V(f)\o^n$ is the generalized
Futaki invariant defined previously by \cite{DT}.

\begin{theorem}
\label{Kstabilitynecessary}
{\rm \cite{T97}}
Let $X$ be a compact K\"ahler manifold with $c_1(X)>0$ and no holomorphic
vector field. If $X$ admits a K\"ahler-Einstein metric,
then $X$ is analytically K-stable in the sense that, for any non-trivial special
degeneration, $f_{X_0}>0$.
\end{theorem}

\medskip
A very rough sketch of the arguments in \cite{T97} is as follows. Let $X_t\subset{\bf CP}^N$ be a special
degeneration of a smooth projective variety $X$,
and let $\phi_t$ be the K\"ahler potential of the Fubini-Study metric on $X_t$,
pulled-back to $X$.

\smallskip
- If the degeneration is non-trivial, then ${\rm
sup}_t\|\phi_t\|_{C^0}=\infty$.

- There exists
constants $C_1,C_2$ so that
$\|\phi_t\|_{C^0}\leq C_1\,J_{\o_0}(\phi_t)+C_2$. This is a Harnack-type
inequality,
and follows from the uniformity of Sobolev constants for the family
$X_t$ of subvarieties of ${\bf CP}^N$.

\smallskip
Assume next that $X$ has a K\"ahler-Einstein metric $\o_{KE}$ and that
${\rm Aut}^0(X)=0$.
By Theorem \ref{mosertrudingertheorem}, it follows that $K_{\o_0}(\phi_t)\to+\infty$ as
$t\to 0$.
On the other hand, we always have an expansion of the form
\be
K_{\o_0}(\phi_t)=f_{X_0}\,\log{1\over |t|}+O(1)
\ee
as $t\to 0$ for a degeneration. It follows that $f_{X_0}>0$, as claimed.
In the Iskovskih examples, $f_{X_0}=0$ always, and thus Theorem
\ref{Kstabilitynecessary}
implies that they do not admit K\"ahler-Einstein metrics.

\subsection{Necessity of Chow-Mumford stability}

Donaldson \cite{D01} has proved that the existence of
a constant scalar curvature metric implies Chow-Mumford
stability. More precisely, we have:

\begin{theorem}
\label{sc1}
Let $X$ be a compact complex manifold and $L\ra X$ a positive
holomorphic line bundle.
Assume that $X$ admits a metric $\o\in c_1(L)$ of constant scalar
curvature,
and that  $\Aut(X,L)/\C^\times$ is discrete.
Then $(X,L^k)$ is Hilbert-Mumford stable
for $k$ sufficiently large.
\end{theorem}
\noindent

\v
The discreteness
assumption was later weakened by Mabuchi \cite{Ma05}.
\v

The key intermediate notion for the
proof of Theorem \ref{sc1} is the notion of balanced imbedding.
By the result of Zhang \cite{Z96} (see also \S 9.1 below), the bundle $L\to X$
is $k$ Chow-Mumford stable if and only if
$L^k\to X$ admits a balanced imbedding into $O(1)\to {\bf CP}^{N_k}$.
Both notions of constant scalar curvature and of balanced imbedding are
related to properties of the density of states $\rho_k(z)$.
On one hand, by the Tian-Yau-Zelditch theorem and Lu's formula,
we have the following expansion for the density of states $\rho_k(z)$,
\be
\rho_k(z)=k^n+{1\over 2}R(\o)k^{n-1}+O(k^{n-2}),
\ee
and thus, the condition that $R(\o)$ is constant means that
$\rho_k(z)$ is constant up to errors of order $O(k^{n-2})$.
On the other hand, if $h$ is a metric on $L$, with respect to which
the density of states $\rho_k(z)$ is constant, then the orthonormal bases
$\us=\{s_\al(z)\}_{\al=0}^{N_k}$ of $H^0(X,L^k)$ with respect to $h$ provide a balanced imbedding.
Indeed,
in view of the formula (\ref{relation}), we have then
\be
\o-{1\over k}\iota_{\us}^*(\o_{FS})=-{i\over 2k}\ddb\log\,\rho_k(z)=0,
\ee
and hence
\be
{1\over V_{FS}}\int_X{\overline{s_\al}(z)s_\beta(z)\over
\sum_{\g=0}^{N_k}|s_\g(z)|^2}\iota_{\us}^*(\o_{FS})^n
=
{1\over V}\int_X{\<s_\beta,s_\al\>\over\rho_k}\o^n
=
{1\over \rho_k\,V}\int_X\<s_\beta,s_\al\>\o^n
=\lambda\,\delta_{\beta\al}
\ee
showing that the imbedding is balanced. Thus, heuristically,
one would like to construct a metric $h$ with $\rho_k(z)$ constant.
In practice, one proceeds in two steps. The first step uses
the Tian-Yau-Zelditch theorem to produce
metrics where $\rho_k(z)$ is constant to arbitrarily high order in $k^{-1}$,
and whence a basis of $H^0(X,L^k)$ which is ``approximately balanced''.
The second step shows how to construct a balanced metric from one which is
approximately balanced.

\v
In this section we give an outline
of the two steps in the proof. The
first is based directly on
\cite{D01}. The second step
is based on \cite{D01} together
with the estimates proved in \cite{PS04}.
This step is motivated by the methods
of \cite{D01} concerning moment maps and infinite
dimesnsional symplectic geometry, but it
does not make explicit use of these notions. Moreover, Deligne's formula
for the curvature of the Deligne pairing,
which plays a key role in the second step below, has an elementary proof (which
is given in \cite{PS03}) and one does not
need to invoke Deligne's theory \cite{De}
to implement the method (thus simplifying the presentation
in \cite{PS04}).

\v
\subsubsection{Approximately balanced imbeddings}

Donaldson \cite{D01} uses the Tian-Yau-Zelditch theorem together with the
formula of Lu \cite{Lu} to prove
the following two lemmas:
\noindent

\begin{lemma}
\label{lemma one} Assume that $\o_\i\in c_1(L)$ has constant scalar curvature and fix $q>0$.
Then for $k>>0$ there exists
a basis $\us=\us(k)$ of $H^0(X,L^k)$ with the following properties.
\be
\label{D basis}
M_{\bar\al\b}\ = \  {1\over V}\I_{\iota_{\us(k)}(X)} {\bar x_{\al}
x_\b\over |x|^2}\
\o_{FS}^n
\ =
\ {k^n\over N+1 }\delta_{\al\b}
\ + \ (E_k)_{\bar\al\b}
\ee
where $\iota_\us:X\rightarrow\P^{N_k}$ is the imbedding
of $X$ via the basis $\us$ and,
for $k$ sufficiently large,
\be\label{Ek}   \|E_k\|_{op} \ \leq \ k^{-q-1}
\ee
where $\|\cdot\|_{op}$ denotes the operator norm.
\end{lemma}
Remark: The matrix $E_k$ vanishes
precisely when the basis $\us(k)$ is
balanced. Thus we may think of
$\us(k)$ as ``approximately balanced''.

\begin{lemma}
\label{lemma two}
There
exists
$C>0$ with the  following properties.
Let $\xi\in  su(N_k+1))$, with $\|\xi\|_{\rm op}=1$ and, for $t\in (-{1\over
10},{1\over 10})$ and $x\in\P^{N_k}$ with $k>>0$, let
$\si_t(x)=\exp(it\xi x)$. Let $\ti\o_\i=k\o_\i$.
Then

\be\label{pullback} \|\iota_{\us(k)}^*\sigma_t^*\o_{FS}-\ti\o_\i\|_{C^4(\ti\o_0)}\ \leq
\  Ct+ O({1\over k})
\ee
\end{lemma}
Remark: Lemma \ref{lemma two} is a rewording of Proposition 27 in \cite{D01}.

\v

The idea of the proof of Lemma \ref{lemma one} goes as follows: Theorem \ref{tyz}
says that the density of states functions has an asymptotic expansion
\be
\r_k(\o)\ = \ k^n+A_1(\o)k^{n-1}+A_2(\o)k^{n-2}+\cdots
\ee
where the $A_p$ are polynomials in the curvature tensor of $\o$ and its
covariant derivatives. If $\o=\o_\i$ has constant scalar curvature,
then Lu's theorem implies that $A_1(\o_\i)$ is constant. The key step
in proving Lemma \ref{lemma one} is to construct a metric $\ti\o$ with
the property
\be \r_k(\ti\o)\ = \ k^n+\ti A_1k^{n-1}+\ti A_2k^{n-2}+\cdots
\ee
where the
$\ti A_j$ are constant for $1\leq j\leq q$. To do this,
we let $\ti\o=\o_\i+{i\over 2k}\pl\bar\pl\eta$ where $\eta$ is
an unknown smooth function. Now the variation $\delta R$
of the scalar curvature under a variation $\delta\phi$
of the potential is given by
\be
\label{scalarcurvaturevariation}
\delta R=-{\cal D}^*{\cal D}\delta\phi+g^{j\bar k}\pl_jR\,\pl_{\bar k}\delta\phi,
\ee
where ${\cal D}$ is the Lichnerowicz operator, mapping scalar functions
to symmetric two-tensors
\be
\label{lichnerowicz}
{\cal D}\delta\phi=\nabla_{\bar k}\nabla_{\bar j}\delta\phi.
\ee
Since $\o$ has constant scalar curvature, we have then
$A_1(\ti\o)=A_1(\o)+{1\over k}\cL(\eta)+ O({1\over k^2})$
and $A_2(\ti\o)=A_2(\o)+O({1\over k})$, where $\cL(\eta)\equiv {\cal D}^*{\cal D}$.
Under
the assumption that $\Aut(L\ra X)$ is discrete, $\cL$ is a self-adjoint elliptic
operator whose kernel consists only of the constants. Thus
\be
\r_k(\ti\o)\ = \ k^n+ A_1(\o_\i)k^{n-1}+
(A_2(\o_\i) + \cL(\eta))k^{n-2}+\cdots
\ee
Now we choose $\ti A_2$ to be the average of $A_2(\o_\i)$ and we choose $\eta$,
the unique smooth function (up to additions of scalars) satisfying
$\cL(\eta)=\ti A_2-A_2(\o_\i)$. This proves the result for $q=2$, and
the statement for arbitrary $q$ is proved in a similar manner.

\subsubsection{Estimates for the normal projection}

This step consists of two lemmas which
concern estimates for the normal projection
operator $\pi_N$ under the assumption of $R-$bounded geometry. We first recall
the definitions of  ``$R$-bounded
geometry" (introduced in \cite{D01})
and of $\pi_N$ and then we state the results.

\v
\noindent
{\it Definition of $R$-bounded geometry}.
Let $k$ be a positive integer, $R$ a
positive real number and let $\ti\o_0=k\o_0$. We say that a K\"ahler metric
$\ti\o\in kc_1(L)$ is $R$-bounded if
\be
\label{R bounded}
\ti\o>R^{-1}\ti\o_0
\ \ \ {\rm and} \ \ \
\|\ti\o-\ti\o_0\|_{C^4(\ti\o_0)}<R
\ee
We say that a basis $\us$ of $H^0(X,L^k)$ is $R$-bounded if the
corrseponding Bergman metric $\ti\o=\iota_\us^*\o_{FS}$ is $R$-bounded.
Note that
(\ref{pullback}) implies that  for $|t|\leq \r$ and $k$ sufficiently large,
the metric
$\iota_{\us}^*\sigma_t^*\o_{FS}$  is $R$ bounded with $R={1\over 2}$.
In other words, the basis $\sigma_t(\us)$ is $1\over 2$-bounded.
\v
\noindent
{\it Definition of the
operator $\pi_N$}.
Let $\us$ be a basis of $H^0(X,L^k)$ and consider the
exact sequence of holomorphic vector bundles
\be
0\ \ra \ TX\ \ \ra \ \iota_{\us}^*T\P^{N_k}\ \ra \ Q \ \ra \ 0
\ee
where $Q\ra X$ is, by defintion, the quotient bundle
$\iota_{\us}^*T\P^{N_k}/TX$. Let ${\cal N}\sub \iota_{\us}^*T\P^{N_k}$
be the orthogonal complement of $TX$ with respect to the Fubini-Study
metric. Then we define $\pi_T:\iota_{\us}^*T\P^{N_k}\ra TX$ and
$\pi_{N}:\iota_{\us}^*T\P^{N_k}\ra {\cal N}$ to be the corresponding
projections. If $V\in \iota_{\us}^*T\P^{N_k}$ then we shall
write $|V|_{FS}$ to be its norm with respect to the Fubini-Study
metric on $T\P^{N_k}$ and we define
\be
\|V\|^2_{L^2(\iota^*_\us\o_{FS})}\ = \ \I_X |V|^2_{FS}\ (\iota^*_\us\o_{FS})^n
\ee

\v
With these notions in place, we can state the two lemmas.
Fix $\xi\in
su(N_k+1)$, with $\|\xi\|=1$.
Let
$X_\xi$ be the holomorphic vector field on $\P^{N_k}$
generated by $\xi$. More explicitly, let us recall that
the holomorphic tangent bundle
$T^{1,0}(\P^{N_k})$ can be describe as follows.
$$
T^{1,0}(\P^{N_k})\ = \ \{(z,v): z,v\in \C^{N_k+1}, z\not=0\}/\sim
$$
where $(z,v)\sim(z',v')$ if and only if $z'=\l v$ and $v'=\l v+\m z$ for
some $\l,\m\in\C$. Then $X_\xi: \P^{N_k}\ra T^{1,0}(\P^{N_k})$
is the map $X_\xi(z)=(z,\xi z)$.
\v
To simplify
the exposition, we
assume throughout that $\Aut(X)$ is discrete (this assumption
can be removed - see \cite{PS04} for details).
\begin{lemma}\label{lemma three}. Let $\us$ be a basis of
$H^0(X,L^k)$ with $1\over 2$-bounded geometry. Then there
exists $C_1>1$ such that
\be\label{79} 1\leq C_1k\|X_\xi\|^2_{L^2(\iota_{\us}^*\o_{FS})}
\ee
\end{lemma}
\begin{lemma}
\label{lemma four}
Let $\us$ be a basis of
$H^0(X,L^k)$ with $1\over 2$-bounded geometry. Then
Then there
exists $C_2>1$ such that
\be\label{712} \|\pi_T X_\xi\|^2_{L^2(\iota_{\us}^*\o_{FS})}\ \leq \
C_2k\|\pi_NX_\xi\|^2_{L^2(\iota_{\us}^*\o_{FS})}
\ee
\end{lemma}
Remark: The assumption that $\Aut(X)$ is discrete implies
that $X_\xi$ is not tangent to $\iota_\us(X)$ and thus,
$\|\pi_NX_\xi\|^2_{L^2(\iota_{\us}^*\o_{FS})}>0$. Lemma
\ref{lemma four} shows that $\|\pi_NX_\xi\|^2_{L^2(\iota_{\us}^*\o_{FS})}$
cannot be too close to zero.

\v
\noindent

\subsubsection{Synthesis}

We show how steps one and two can be combined to
prove Theorem \ref{sc1}. We must show that the
existence of the almost balanced basis $\us(k)$ implies
the existence of a balanced basis.
For this, it is simplest to reformulate the problem
in terms of the functional $-F_{\iota_{\us}^*(\o_{FS})}^0$:
our task is to find a critical point (in fact, a minimum) for this
functional. This is now easy, since
section \S 7.2.1 provides upper bounds for the absolute value of
its first derivatives, while section \S 7.2.2 provides a positive lower bound for
its second derivative.

\medskip
We first introduce some necessary notation: Let $\ti\o(k)=\iota_{\us(k)}^*\o_{FS}$.
If
$\sigma\in SL(N_k+1,\C)$ we let $\phi_\sigma(x)=\log{|\sigma(x)|^2\over |x|^2}$.
Then to prove existence of a balanced metric,
in view of Theorem \ref{balancedcritical}, it suffices to show
that $-F^0_{\ti\o(k)}(\phi_\sigma)$ achieves its minimum at some
point $\sigma\in SL(N_k+1)$. Fix $\xi\in su(N_k+1)$
and let
$\si_t(x)=\exp(it\xi x)$ and $\phi_t(x)=\phi_{\sigma(t)}$ . Consider
\be
-(n+1)V_{FS}F^0(t)\ = \
-(n+1)V_{FS}F^0_{\ti\o(k)}(\phi_t)\ = \
\I_{\iota_\us(X)}\phi_t\cdot\s_{j=0}^n(\si_t^*\o_{FS}^j
\wedge\o_{FS}^{n-j})
\ee
Then (\ref{derivativesF0}) combined with Lemma \ref{lemma one} yields
\bea
\label{first der}
  |(n+1)\dot F^0(0)|\ &=& \ \left|{1\over V}\I_{\iota_\us(X)}{x^*\xi x\over x^*x}
\cdot\o_{\rm FS}^n\right|
\ =
\  |{\rm Tr}(\xi E_k)| \nonumber\\
&\leq& \ {(N_k+1)}^{1\over 2}\|\xi\|\cdot \|E_k\|_{\rm
op}
\
\leq
\ k^{n/2-q-1 }
\eea
On the other hand we have the following formula of Deligne \cite{De} for $\ddot F^0$
(see (\ref{derivativesF0}) and \cite{PS03} for an elementary proof):
\be
\label{second derivative}
-\ddot F^0(t)\ = \ \|\pi_NX_\xi\|^2_{L^2(\iota_{\us}^*\sigma_t^*\o_{FS})}\ = \
\I_X |\pi_NX_\xi|_{FS}^2\cdot (\iota_{\us}^*\sigma_t^*\o_{FS})^n
\ee
To estimate $ \|\pi_NX_\xi\|^2$ we wish to use Lemma \ref{lemma three} and
Lemma \ref{lemma four}. First we apply Lemma \ref{lemma two} to deduce
that the basis $\sigma_t\cdot\us(k)$ is $1\over 2$-bounded for $|t|\leq {1\over 4C}$
and for $k>>0$. Now Lemma \ref{lemma three} and Lemma \ref{lemma four} imply
that for such $t$,
\be
\label{second der}
-\ddot F^0(t)\ \geq \ \ {1\over
C_2(k+1)}(C_2k\|\pi_NX_\xi\|^2+\|\pi_NX_\xi\|^2)
\ \geq \ {1\over C_2(k+1)}\|X_\xi\|^2\ \geq \ {1\over C_1C_2k(k+1)}
\ee
Thus, if we take $q > n/2+1$, comparing (\ref{first der}) and
(\ref{second der}) we see that, for $k$ sufficiently large,
\be
-F^0(t)>-F^0(0)\ \ \ {\rm if} \ \ \ |t|\geq \ {1\over 4C_1}\ .
\ee
This shows that $-F^0_{\ti\o(k)}(\phi_\sigma)$ achieves its minimum at some
point $\sigma\in SL(N_k+1)$ which is of the form
$\sigma=\exp(i\xi)$ where $\xi\in su(N_k+1)$ and $\|\xi\|\leq {1\over 4C_1}$.
\v

\subsection{Necessity of semi $K$-stability}

The following theorem due to Donaldson \cite{D05} provides
an attractive lower bound for $\|R(\o)-\mu n\|_{L^2}^2$ in terms
of Futaki invariants of test degenerations:

\begin{theorem}
\label{lbc}
Let $L\to X$ be a positive line bundle. Then we have
\be
\label{calabilowerbound}
{\rm inf}_{\o\in c_1(L)}\|R(\o)-\mu n\|_{L^2}
\geq
{\rm sup}_{\cal T}(-{{\rm F}({\cal T})\over D({\cal T})})
\ee
where the supremum on the right hand side runs over all test configurations
${\cal T}$,
and the invariant $D({\cal T})$ of a test configuration ${\cal T}$ is
defined by,
\be
D({\cal T})={\rm lim}_{k\to+\infty}{\|A_k\|\over k^{{n\over 2}+1}}.
\ee
where $\|A\|=(\s_{i,j}|a_{ij}|^2)^{1/2}$ is the Hilbert-Schmidt norm.
\end{theorem}

\medskip
This implies immediately the following corollary:

\begin{corollary}
\label{Ksemistability}
If $L\to X$ is a positive line bundle over a compact complex manifold,
and if $X$ admits a metric of constant scalar curvature in $c_1(L)$,
then $L\to X$ is $K$-semistable.
\end{corollary}

\v

Another lower bound for the Calabi energy involving instead a version of the
Futaki invariant along geodesic rays can be found in \cite{Ch06}.

\smallskip
In the particular case when $L=K_X^{-1}$,
it is known that the boundedness of
the Mabuchi K-energy from below on $c_1(L)$ implies that
${\rm inf}_{\o\in c_1(L)}\|R(\o)-\mu n\|_{L^2}=0$
(see e.g. \cite{B, PS05, PSSW2}).
Thus we also have:

\begin{corollary}
Let $X$ be a compact complex manifold with $c_1(X)>0$.
If the Mabuchi K-energy is bounded from
below on $c_1(X)$, then $X$ is K-semistable.
\end{corollary}

\smallskip
Lower bounds for the K-energy have been obtained in certain cases by Weinkove \cite{W}
using the $J$-flow.

\medskip

The proof of Theorem \ref{lbc} is so basic and elegant that we can include it
here. In fact, it relies on the properties of $F_{\o_0}^0(\phi)$ on the spaces
${\cal K}_k$ of Bergman metrics which have already been discussed in section \S 5.
The version that we present here is slightly shorter, thanks to the use of
Lemma \ref{FTmu} below.

\smallskip

Let $\cT$ be a test configuration, and let $A_k, B_k$ be the endomorphisms
on $H^0(X_0,L_0^k)$ as defined in section \S 6.1.6. Let
$\sigma_t^{(k)}=e^{tB_k}:{\bf R}\rightarrow GL(N_k+1)$,
$\underline{\sigma}_t^{(k)}=e^{tA_k}:{\bf R}\rightarrow SL(N_k+1)$
be the one-parameter subgroup defined by $B_k$ and $A_k$.
For $k>>0 $ embed $X\sub{\bf CP}^{N_k}$ by a basis $\us$ of $H^0(X,L^k)$.
Let $\phi_t(x)=\log{|\sigma_t^{(k))}(x)|^2\over |x|^2}$ and
define
\be
\m(A_k)\ = \ -\lim_{t\to-\i}V_{FS}\,\dot F^0_{\iota_{\us}^*(\o_{FS})}(\phi_t)
\ee
Then we have the following alternative description of the Donaldson-Futaki
invariant of the test configuration $\cT$ \cite{PS07}

\begin{lemma}
\label{FTmu}
Let $\cT$ be a test configuration. Then
\be
F(\cT)\ = \ \lim_{k\to\i} {\m(A_k)\over k^n}
\ee
\end{lemma}

\medskip
Assuming the lemma for the moment, we can give the proof of the theorem.
Fix $h$ a metric on $L$, and let $\us$ be the Kodaira imbedding
defined by an orthonormal basis $\us$ of $H^0(X,L^k)$ with respect
to the $L^2$ norm defined by $h$ and $\o=-{i\over 2}\ddb\log\,h$.
By the convexity of $F^0_{\iota_{\us}^*(\o_{FS})}$ along one-parameter subgroups,
we have
\be
-{\rm lim}_{t\to-\infty}V_{FS}\,\dot F^0_{\iota_{\us}^*(\o_{FS})}(t)
\ \leq\
-V_{FS}\,\dot F^0_{\iota_{\us}^*(\o_{FS})}(0)
\ee
We can now evaluate both sides of the equation, using Lemma \ref{FTmu}
for the numerical invariant $\mu(A_k)$ and (\ref{derivativesF0}) for
$-\dot F^0_{\iota_{\us}^*(\o_{FS})}$ (with $\delta\sigma=\delta\sigma^*=A_k$).
This gives
\be
\label{muAk}
{\mu (A_k)\over\|A_k\|}
\leq
V_{FS}\,{{\rm Tr}(A_k\cdot {\underline M})\over \|A_k\|}\leq \,V_{FS}\,\|\underline{M}\|,
\ee
where $\underline M$ is the traceless part of the matrix $M$ defined by (\ref{M})
and (\ref{M2}),
and thus given by
\be
{\underline M}_{\al\bar\beta}
=
-{1\over kV_{FS}}\int_X\<s_\al,s_\beta\>(R-\bar R)\o^n +O({1\over k^{n+2}}).
\ee
We estimate the Hilbert-Schmidt norm of $\underline{M}$. By an orthonormal
transformation, we may assume that ${\underline M}_{\al\bar\beta}$ is diagonal.
The diagonal entries can be estimated by the Cauchy-Schwarz inequality
and the fact that $\|s_\al\|_{L^2}^2=1$,
\bea
V_{FS}\underline{M}_{\al\bar\al}
&=& -{1\over k}\int_X |s|_{h^k}^2(R(\o)-\bar R)\o^n+O({1\over k^{2}})
\nonumber\\
&\leq&
{1\over k}(\int_X |s|_{h^k}^2|R(\o)-\bar R|^2\o^n)^{1\over 2}+O({1\over k^{2}}).
\eea
Applying the Tian-Yau-Zelditch theorem, we obtain
\bea
V_{FS}^2\sum_{\al=0}^{N_k}\underline{M}_{\al\bar\al}^2
&\leq& {1\over k^{2}}\int_X\rho_k(z)|R(\o)-\bar R|^2\o^n
+O({1\over k^{-n+3}})
\nonumber\\
&\leq & {1\over k^{-n+2}}\|R(\o)-\bar R\|_{L^2}^2+O({1\over k^{-n+3}}).
\eea
Substituting this bound into (\ref{muAk}) and letting $k\to\infty$ gives the desired statement. Q.E.D.

\medskip
\noindent
{\it Proof of Lemma \ref{FTmu}}: We begin by recalling some basic facts about the Mumford numerical
invariant for projective varieties.
Let
$Z\sub \P^N$ be a smooth subvariety, and
$B$ an  $(N+1)\times (N+1)$ matrix. Let
$\o_{FS}$ be the Fubini-Study metric on $\P^N$. We shall also denote by
$\o_{FS}$ the restriction of the Fubini-Study metric to $Z$. For $t\in\R$
let $\sigma_t\in GL(N+1,\C)$ be the matrix $\sigma_t=e^{tB}$ and
let $\psi_t:\P^N\ra\R$ be the function
\be
\label{psi}
\psi_t(z)\ = \ \log{|\sigma_t z|^2\over |z|^2}
\ee
Here we view $z$ as an element in $\P^N$ and, when there is no
fear of confusion, a column vector in $\C^{N+1}$.
\v
Then $\psi_t$ is a smooth path in $\cH$: In fact,
$\sigma_t^*\o_{FS}=\o_{FS}+{i\over 2}\ddb\psi_t$.
Define
\be
\label{mu}
\m(Z,B)\ = \ -\lim_{t\ra -\i}V_{FS}\,\dot F^0_{\o_{FS}}(\psi_t)\ = \ - \dot
V_{FS}\,F^0(-\i)
\ee
where $V_{FS}$ denotes the volume of $Z$ with respect to $\o_{FS}$.
The
function $F^0(t)=F^0_{\o_{FS}}(\psi_t):\R\ra\R$
is  convex  so the above limit exists.
Suppose $\l:\C^\times\ra GL(N+1,\C)$
an algebraic homomorphism, let $B\in gl(N+1,\C)$ be such that
$\l(e^t)=e^{tB}$ for all $t\in\R$ and $Z^{(0)}=\lim_{\tau\to
0}\l(\tau)(Z)$ the flat limit of~$Z$. Thus $Z^{(0)}\sub\P^N$ is a subscheme
of $\P^N$ with the same Hilbert  polynomial as $Z$.
Let $M_0=O(1)|_{Z^{(0)}}$ and $M=O(1)|_Z$. We assume that $Z\sub\P^N$ is
an imbedding of $Z$ by a basis of $H^0(Z,M)$.
Then $\l(\tau)$ defines an automorphism of $H^0(Z^{(0)},M_0^p)$
and, for $p>0$, we let $\tilde w(Z, B, p)$ be the weight of this action
on $\det(H^0(Z_0,M_0^p))$. One easily sees that
\be
\ti w(Z,B,1)\ = \ Tr(B)
\ee
and, if $I$ is the identity matrix,
\be
\ti w(Z, I,p)\ = \ {p\dim H^0(Z^{(0)},M^p)}
\ee
It is then known by the work of Mumford \cite{Mu}
and Zhang \cite{Z96} that $\tilde w(Z,B,p)$ is a
polynomial in $p$ for $p$ large such that for $B=A\in sl(N+1,\C)$ we have
\be\label{tildew} \tilde w(Z,A, p) \ = \ {\m(Z,A)\over (n+1)!}\cdot p^{n+1}\ + \ O(p^n)
\ \ {\rm and} \ \ \tilde w(Z^{(0)},A,1)=0
\ee

We apply these general remarks to the case considered in
Lemma \ref{FTmu}, that is, to the case
$B=B_k$.
Thus we let
$L\to X$ be an ample line bundle,
and  $\rho:{\bf C}^\times\to {\rm Aut}({\cal L}\to{\cal X}\to{\bf C})$
 a test configuration $\cT$.  Let $k$ be an integer such that $L^k$
is very ample. By Lemma \ref{tc}, there exists
a basis $\underline s$ of $H^0(X,L^k)$,
and an imbedding
\be
I_{\underline s}:
({\cal L}^k\to{\cal X}\to{\bf C})
\to
(O(1)\times{\bf C}\to {\bf CP}^{N_k}\times{\bf C}\to{\bf C})
\ee
which restricts to $\iota_{\underline s}$ on the fiber $X_1$
and intertwines $\rho(\tau)$ and $B_k$, i.e.,
\be
\label{intertwining}
I_{\underline s}(\rho(\tau)l_w)
=
(\tau^{B_k}\cdot I_{\underline s}(l_w),\tau w)\ \ \
\hbox{for each $l\in \cL^k$}
\ee
Let $Z_k\sub\P^{N_k}$ be the image of $X$ by the map $\iota_\us$.
Then, $(L^r)^p=L^{rp}$ implies (with $k=rp$)
\be
\ti w(Z_r,B_r,p)\ = \ \ti w(Z_k,B_k,1)\ = \ Tr(B_k)
\ee
In particular we see that $Tr(B_k)$ is a polynomial
in $k$ of degree $n+1$ (for $k>>1$).
\v
On the other hand, $A_k=B_k-{1\over d_k}Tr(B_k)$ (where $d_k=N_k+1$)
implies
\be
rd_r Tr(B_k)= rd_r\ti w(Z_r,B_r,p)=\ti w(Z_r,rd_rB_r,p)\ = \
\ti w(Z_r, rd_rA_r,p)\ + \ \ti w(Z_r, rTr(B_r),p)
\ee
Thus
\be
rd_r Tr(B_k) - kd_kTr(B_r)\ = \ rd_r\ti w(Z_r,A_r,p) =
\ e_T(r)k^{n+1}\ + \ O(k^n) \ \ {\rm for} \ \ k>>1
\ee
where $e_T$ is a polynomial in $r$ (for $r>>1$) of degree at most $n$.
If follows from the definition of $F(T)$ that
$F(T)$ is the leading coefficient of $e_T(r)$. Comparing with
(\ref{tildew}) we get
\be
\lim_{r\to\i}\ {\m(Z_r,rN_rA_r)\over r^nr^{n+1} (n+1)!}\ = \ F(T)
\ee
Since $r^{-n}N_r = {1\over n!}\I\o^n + O(r^{-1})$,
we obtain the desired formula.

\newpage
\section{Sufficient Conditions: the K\"ahler-Einstein Case}
\setcounter{equation}{0}

We describe next some of what is known in the direction of sufficiency
of stability conditions. The case of K\"ahler-Einstein is better understood,
so we begin with this case. As noted earlier, the constant scalar curvature
equation reduces to a complex Monge-Amp\`ere equation, which can be
investigated either as an elliptic equation by the method of continuity,
or as a parabolic flow, giving rise to the K\"ahler-Ricci flow.

\subsection{The $\alpha$-invariant}

An approach to the problem of K\"ahler-Einstein metrics is through the
$\alpha$-invariant, defined on a compact K\"ahler manifold $(X,\o_0)$ as follows
\cite{T87},
\be
\label{alpha}
\alpha(X)={\rm sup}\{\kappa>0; {\rm sup}_{\phi}\int_X e^{-\kappa \phi}\o_0^n<\infty\},
\ee
where the supremum on the right hand side is taken over all $\phi$ which
are $\o_0$-plurisubharmonic and normalized by ${\rm sup}_X\phi=0$. It is shown in Tian \cite{T87},
using H\"ormander's
estimate for subharmonic functions \cite{Ho}, that $\alpha(X)>0$, and
that the lower bound $\alpha(X)>{n\over n+1}$ would imply the existence
of a K\"ahler-Einstein metric on $X$. Other applications of lower bounds for
$\alpha(X)$ are in \cite{TY87}. For toric manifolds, the $\alpha$-invariant has been
completely determined in Song \cite{S03}, generalizing earlier special cases established
in Batyrev and Selivanova \cite{BS} and \cite{S02}. It has very recently been evaluated for
the Mukai-Umemura 3-fold in Donaldson \cite{D07}.
More general exponential estimates for plurisubharmonic functions
with respect to certain probability measures are in Dinh, Nguyen, and Sibony \cite{DNS}.
The importance role of an $\alpha$-invariant for tame symplectic manifolds,
if it is strictly positive,
has been brought to light by Tosatti, Weinkove, and Yau \cite{TWY}.

\subsection{Nadel's multiplier ideal sheaves criterion}

In the case of vector bundles, as proved by Donaldson \cite{D87} and Uhlenbeck-Yau \cite{UY},
the obstruction to the existence of a Hermitian-Einstein metric is precisely the
presence of destabilizing subsheaves in the sense of Mumford-Takemoto.
It is natural to expect that the obstruction to K\"ahler-Einstein metrics,
and more generally, to K\"ahler metrics of constant scalar curvature,
can ultimately be also encoded in suitable notions of destabilizing sheaves.
An early important result is the following theorem of Nadel \cite{Na},
which we quote here in the simpler version of Demailly-Koll\'ar \cite{DK}. Let
$(X,\o_0)$ be a compact K\"ahler manifold with $c_1(X)>0$.
For each $\o_0$-plurisubharmonic $\psi$
(that is, upper semi-continuous and satisfying $\o_0+{i\over 2}\ddb\psi\geq 0$),
define the multiplier ideal sheaf ${\cal I}(\psi)$ by
\be
\label{multiplierideal}
{\cal I}_z(\psi)
=\{f; \exists U\ni z,\ f\in {\cal O}(U),
\ \int_X e^{-\psi}|f|^2\o_0^n<\infty\ \}
\ee

\begin{theorem}
\label{nadel}
If $(X,\o_0)$ does not admit a K\"ahler-Einstein
metric, then for any $p\in ({n\over n+1},1]$, there exists
a $\o_0$-plurisubharmonic function $\psi$
so that the multiplier ideal sheaf ${\cal I}(p\psi)$
defines a proper, coherent analytic sheaf on $X$ with acyclic
cohomology, i.e.,
\be
\label{acyclic}
H^q(X,{\cal J}(p\psi))=0,\qquad q\geq 1.
\ee
If $X$ admits a compact subgroup $G$ of holomorphic automorphisms,
and $\o_0$ is $G$-invariant, then ${\cal I}(p\psi)$ and the corresponding
subscheme are also $G$-invariant.
\end{theorem}

Multiplier ideal sheaves were introduced by Kohn \cite{K79} in the context
of subelliptic estimates for the $\bar\pl$-Neumann problem. Their applications
to complex geometry have been pioneered by Siu \cite{Siu2, Siu3}.

\smallskip

The original proof of Nadel's theorem \cite{Na} is based on the method of continuity
for the complex Monge-Amp\`ere equation (\ref{volumeequation}),
\be
\label{MAcontinuity}
(\o_0+{i\over 2}\ddb\phi)^n=e^{f_0-t\phi}\o_0^n,
\qquad 0\leq t\leq 1.
\ee
By Yau's estimates \cite{Y78}, a $C^\infty$ solution of (\ref{volumeequation}) exists
if the equation (\ref{MAcontinuity}) admits a $C^0$ a priori estimate
\be
\label{C0}
{\rm sup}_t\|\phi\|_{C^0}\leq C<\infty.
\ee
On the other hand, the solutions of the equation (\ref{MAcontinuity}) satisfy the following
\cite{Siu1, T87}
\bea
\label{MAharnack}
{1\over V}\int_X (-\phi)\o_\phi^n
&\leq& n\,{1\over V}\int_X\phi\o_0^n
\nonumber\\
{\rm osc}\,\phi
&\leq& A{1\over V}\int_X \phi\o_0^n +B
\eea
The second inequality is a Harnack-type inequality which can be proved by Moser iteration,
since the equation (\ref{MAcontinuity}) implies that the Ricci curvature
of $\o_\phi$ is bounded from below for $t\geq \epsilon>0$, and the Sobolev
constants with respect to $\o_\phi$ are then uniformly bounded from below.
Another important observation is the following lemma:

\begin{lemma}
\label{nadellemma}
If there exists a constant $p\in ({n\over n+1},1]$ so that
\be
\label{nadelcondition}
{\rm sup}_t\,{1\over V}\int_X {\rm exp}\bigg(-p(\phi-{1\over
V}\int_X\phi\o_0^n)\bigg)\o_0^n <\infty,
\ee
then it follows that
\be
\label{technical}
{\rm sup}_t\,{1\over V}\int_X\phi\,\o_0^n\leq C<\infty.
\ee
\end{lemma}

\medskip
\noindent
{\it Proof.} Since the logarithm is concave, the given inequality implies
\be
C>{1\over V}\log\big(\int_X e^{(t-p)\phi-f_0+p{1\over V}\int_X\phi\o_0^n}\o_\phi^n \big)
\geq
{1\over V}
\int_X ((t-p)\phi-f_0)\o_\phi^n+p {1\over
V}\int_X\phi\o_0^n,
\ee
and thus, in view of (\ref{MAharnack}),
\be
p {1\over V}\int_X\phi\o_0^n
\leq (t-p){1\over V}\int_X(-\phi)\o_\phi^n+C
\leq
n(1-p){1\over V}\int_X\phi\o_0^n+C.
\ee
where the last inequality applies if ${1\over V}\int_X(-\phi)\o_\phi^n\geq 0$.
The desired bound follows in this case if $p>n(1-p)$. When ${1\over V}\int_X(-\phi)\o_\phi^n<0$,
we just observe that it suffices to prove the upper bound on ${1\over V}\int_X\phi\o_0^n$
for $t\geq t_0$, for some fixed small time $t_0>0$. Since the given bound (\ref{nadelcondition})
implies a similar bound
for $p_0={1\over 2}t_0$, we can just apply the first inequality in (\ref{technical})
with $p\to p_0$ and deduce that ${1\over V}\int_X\phi\o_0^n$ is uniformly bounded from above.
Q.E.D.

\medskip
The proof of Nadel's theorem can now be completed as follows. If $X$ does
not admit a K\"ahler-Einstein metric, then the $C^0$ estimate (\ref{C0}) must fail.
Now, if we normalize $f_0$ by
\be
{1\over V}\int_X e^{f_0}\o_0^n=1
\ee
then the equation (\ref{MAcontinuity}) implies that $\phi$ must vanish somewhere,
and thus $\|\phi\|_{C^0}\leq {\rm osc}\,\phi$. It follows that ${\rm osc}\,\phi$ is unbounded,
and in view of the Harnack inequality (\ref{MAharnack}) that ${1\over V}\int_X \phi\o_0^n$
is unbounded. By Lemma \ref{nadellemma}, for any $p\in ({n\over n+1},1]$,
the left hand side of the inequality (\ref{nadelcondition}) is unbounded
for some sequence $\phi=\phi(t_j)$.
Let $\psi$ be a limit point of $\phi(t_i)-{1\over V}\int_X\phi(t_i)\o_0^n$.
By the semi-continuity theorem of Demailly-Koll\'ar \cite{DK},
$e^{-p\psi}$ is not in $L^1$, and thus ${\cal I}(p\psi)$ is a proper sheaf.
It is also coherent and acyclic as a consequence of the general theory of multiplier
ideal sheaves. This completes the proof of Theorem \ref{nadel}.

\subsection{The K\"ahler-Ricci flow}

Although there has been important developments since the late 1980's, particularly
concerning complex surfaces \cite{TY87, T90}, we jump now to the more
recent progresses, based on the K\"ahler-Ricci flow.

\medskip
Let $(X,\o_0)$ be a compact K\"ahler manifold, with $\mu\o_0\in c_1(X)$.
Recall that the K\"ahler-Ricci flow is the flow of metrics
defined by (\ref{KRmetric}). Since it preserves the K\"ahler class,
we may write as usual $g_{\bar kj}=(g_0)_{\bar kj}+\pl_j\pl_{\bar k}\phi$,
$\o_\phi=\o_0+{i\over 2}\ddb\phi$, and the flow for the metrics $g_{\bar kj}$
is equivalent to
the flow (\ref{KRpotential}) for the potentials $\phi$.
It follows readily from the maximum principle that
\be
\|\phi\|_{C^M(X\times [0,T))}\leq C_{MT},
\ee
and thus the flow exists for all $t\in [0,\infty)$ \cite{Ca}. The main issue
is the convergence of the flow. It has been shown to converge when $\mu<0$ and
$\mu=0$ \cite{Ca}, and thus we restrict to the case $\mu>0$, which is the case of
positive curvature. Clearly $\mu$ can be normalized to be $\mu=1$, and we shall
do so henceforth. For the convenience of the reader, we reproduce here the
equation
\be
\label{KRpotentialbis}
\dot\phi =\log {\o_\phi^n\over\o_0^n}+\phi-f_0,
\qquad \phi(0)=c_0
\ee
where $c_0$ is a constant, and $f_0$ is the Ricci potential of the initial
K\"ahler form, defined with normalization by
\be
\label{riccipotentialnormalized}
Ric(\o_0)-\o_0={i\over 2}\ddb f_0,
\qquad {1\over V}\int_X e^f\o_0^n=1.
\ee

.

\subsubsection{Perelman's estimates}

The following are key estimates for the K\"ahler-Ricci flow. They are all
consequences
of Perelman's crucial monotonicity formulas
\cite{Pe, CH, ST}, with the first three statements
directly due to him, and the last one to Ye \cite{Ye} and Zhang \cite{Zq} :

\begin{theorem}
\label{perelmanestimates}
{\rm (i)} The Ricci potential $f$ satisfies the following estimates along
the
K\"ahler-Ricci flow,
\be
{\rm sup}_{t\geq 0}(\|f\|_{C^0}+\|\nabla f\|_{C^0}
+
\|\Delta f\|_{C^0}) <\infty.
\ee
{\rm (ii)} The diameters of $X$ with respect to the metrics $g_{\bar kj}(t)$
are uniformly bounded for $t\geq 0$.

\noindent
{\rm (iii)} Let $\rho>0$ be fixed. Then there exists a constant $c>0$ so
that for all
$x\in X$, $t\geq 0$, and $r$ with $0<r\leq\rho$, we have
\be
\label{noncollapse}
\int_{B_r(x)}\o_\phi^n > c\, r^{2n},
\ee
where $B_r(x)$ is the geodesic ball centered at $x$ of radius $r$ with
respect to the
metric $g_{\bar kj}(t)$.

\noindent
{\rm (iv)} There exists a constant $C$, independent of $t\geq 0$, so that
the
Sobolev inequality
\be
\label{uniformsobolev}
\|u\|_{L^{2n\over n-1}}\leq C\, (\|\nabla u\|_{L^2}+\|u\|_{L^2}),
\qquad u\in C^\infty(X),
\ee
holds, with all norms taken with respect to the metric $g_{\bar kj}(t)$.
\end{theorem}

\medskip
All these estimates are for quantities which depend only on the metrics
$g_{\bar kj}$,
and not on the potentials $\phi$ and the normalization $c_0$
for the initial potential. To translate them into
estimates for $\phi$, we need to pick a precise normalization for $\phi$.
First, we note that the quantity $\int_0^\infty\|\na\dot\phi\|^2e^{-t} dt$
(first written down in \cite{CT}; see also \cite{Li})
is finite and independent of the choice of initial condition $c_0$.
This is because two solutions of (\ref{KRpotentialbis}) with different initial
values $c_0$ and $\tilde c_0$ differ by the expression
$(\tilde c_0-c_0)\,e^t$, which cancels out in the norm $\|\na\dot\phi\|^2$.
Furthermore,
\be
\dot\phi=-f+\al(t),
\ee
for some $\al(t)$ independent of $z$,
and by Perelman's estimate (i) in Theorem \ref{perelmanestimates}
above, the quantity $\|\na\dot\phi\|^2$ is bounded, and hence the integral
in $t$ converges.
It is then shown in \cite{PSS} that,
if the initial value $c_0$ is chosen to be
\be
\label{initialc}
c_0=\int_0^\infty\|\na\dot\phi\|^2e^{-t} dt+{1\over V}\int_X f_0\o_0^n,
\ee
then the constant $\al(t)$ is uniformly bounded, and
(i) in Theorem \ref{perelmanestimates} implies
\be
\label{perelman}
\|\dot\phi\|_{C^0}\leq C.
\ee

\subsubsection{Energy functionals and the K\"ahler-Ricci flow}

Recall the energy functionals $F_\o^0(\phi)$, $F_\o(\phi)$, and $K_\o(\phi)$
introduced in section \S 3.
It is well-known that both energy functionals $F_\o(\phi)$ and $K_\o(\phi)$
decrease along the K\"ahler-Ricci flow. For the $K$-energy, this is an
immediate consequence of its variational definition, and the fact
that $\dot\phi=-f+\al$, where $f$ is the Ricci potential,
\be
{d\over dt}K_\o(\phi)=-{1\over V}\int_X\dot\phi(R-n)\o_\phi^n
=
{1\over V}\int_X \dot\phi\Delta \dot\phi\o_\phi^n=-{1\over
V}\int_X|\nabla\dot\phi|^2\o_\phi^n.
\ee
As for the functional $F_\o(\phi)$, we have, with $W\equiv \int_X
e^{-\dot\phi}\o_\phi^n$,
\bea
{d\over dt}F_\o(\phi)
=\int_X ({1\over W}e^{-\dot\phi}-{1\over V})\dot\phi\o_\phi^n
=\int_X ({1\over W}e^{-\dot\phi}-{1\over V})
(\dot\phi-\log{V\over W})\o_\phi^n,
\eea
which is negative since the integrand is of the form $-(x-y)(e^x-e^y)\leq
0$.

\medskip
Henceforth, we choose the initial value $c_0$ for
the flow (\ref{KRpotentialbis}) to be (\ref{initialc}),
so that the inequality (\ref{perelman}) holds.
The following identity and bounds will also be very useful:

\begin{lemma}
\label{harnack1}
{\rm (i)} The following identity holds along the K\"ahler-Ricci flow
\be
K_\o(\phi)-F_\o^0(\phi)-{1\over V}\int_X\dot\phi\,\o_\phi^n=C.
\ee
{\rm (ii)} There exists a constant $C$ so that
\be
|F_\o(\phi)-K_\o(\phi)|+|F_\o^0(\phi)-K_\o(\phi)|\leq C
\ee
along the K\"ahler-Ricci flow.
\end{lemma}

\noindent
{\it Proof}. From the derivative of $F_\o^0(\phi)$ and the definition of the
flow, we have
\be
-{d\over dt}F_\o^0(\phi)
={1\over V}\int_X\dot\phi\o_\phi^n
=
{1\over V}\int_X\ddot\phi\o_\phi^n,
\ee
since
$0={d\over dt}({1\over V}\int_X\o_\phi^n)
=
{d\over dt}\bigg({1\over V}\int_Xe^{f_0-\phi+\dot\phi}\o_0^n\bigg)
=
-{1\over V}\int_X\dot\phi\o_\phi^n
+
{1\over V}\int_X\ddot\phi\o_\phi^n$.
On the other hand,
\bea
{1\over V}\int_X\ddot\phi\o_\phi^n
&=&
{d\over dt}\bigg({1\over V}\int_X\dot\phi\o_\phi^n\bigg)
-
{1\over V}\int_X\dot\phi\Delta\dot\phi\,\o_\phi^n
=
{d\over dt}\bigg({1\over V}\int_X\dot\phi\o_\phi^n\bigg)
+
{1\over V}\int_X|\nabla\dot\phi|^2\o_\phi^n
\nonumber\\
&=&
{d\over dt}\bigg({1\over V}\int_X\dot\phi\o_\phi^n\bigg)
-{d\over dt}K_\o(\phi).
\eea
This establishes (i). Since $|\dot\phi|$ is bounded by Perelman's results,
$|F_\o^0(\phi)-K_\o(\phi)|$ is bounded.
Finally, consider the difference between $F_\o^0(\phi)$ and $F_\o(\phi)$,
\be
|\log\,({1\over V}\int_X e^{f_0-\phi}\o_0^n)|
=
|\log({1\over V}\int_X e^{-\dot\phi}\o_\phi^n)|
\leq
C{1\over V}\int_X\o_\phi^n=C,
\ee
and thus $|F_\o(\phi)-K_\o(\phi)|$ is bounded as well. Q.E.D.

\bigskip
It is convenient to group together the essential inequalities between
the quantities ${1\over V}\int_X \phi\,\o_0^n$,
${1\over V}\int_X(-\phi)\o_\phi^n$, and $J_\o(\phi)$
in the following lemma. Note that the first inequality
is the analogue for the K\"ahler-Ricci flow of
the first inequality in (\ref{MAharnack}) for the complex
Monge-Amp\`ere equation:

\begin{lemma}
\label{harnack2}
There exists constants $C$ so that
\bea
&&
{1\over n}{1\over V}\int_X (-\phi)\,\o_\phi^n-C\ \leq\ J_\o(\phi)
\
\leq {1\over V}\int_X \phi\,\o_0^n+C
\nonumber\\
&&
{1\over V}\int_X \phi\,\o_0^n
\ \leq \ n\,{1\over V}\int_X(- \phi)\,\o_\phi^n
-(n+1)K_\o(\phi)+C.
\eea
uniformly along the K\"ahler-Ricci flow.
\end{lemma}

\medskip
\noindent
{\it Proof.}
Since $F_\o(\phi)$ is monotone decreasing along the K\"ahler-Ricci flow, we
have
$F_\o(\phi)\leq 0$. As noted in the proof of Lemma \ref{harnack1},
this implies that $F_\o^0(\phi)\leq C$.
Now let $I_{\o}(\phi)$ be the following functional,
\be
\label{I}
I_\o(\phi)={1\over V}\int_X\phi(\o^n-\o_\phi^n)
={i\over 2}\sum_{k=0}^{n-1}\int_X \pl\phi\wedge\bar\pl\phi\wedge\o^{n-1-k}\wedge\o_\phi^k.
\ee
Comparing with the expression for $J_\o(\phi)$ in \S 3.3.2,
we readily see that
\be
0\leq {1\over n}J_\o(\phi)\leq {1\over n+1}I_\o(\phi)
\leq J_\o(\phi).
\ee
Furthermore, the functional $F_\o^0(\phi)$ can be written in terms of $I_\o(\phi)$
and $J_\o(\phi)$ in two different ways,
\be
F_\o^0(\phi)
=
J_\o(\phi)-{1\over V}\int_X\phi\o^n=
-
[(I_\o-J_\o)(\phi)+{1\over V}\int_X\phi\o_\phi^n].
\ee
Using the first way of writing $F_\o^0(\phi)$,
we obtain the inequality on the right of the first statement
of Lemma \ref{harnack2}.
Using the second way of writing $F_\o^0(\phi)$,
we obtain
\be
{1\over V}\int_X(-\phi)\o_\phi^n
\leq
(I_\o-J_\o)(\phi)+C
\leq n\,J_\o(\phi)+C.
\ee
and the inequality on the left of the first statement of Lemma \ref{harnack2}
also follows.
The second statement of Lemma \ref{harnack2}
is an easy consequence of Lemma \ref{harnack1}. The
inequality $F_\o^0(\phi)-K_\o(\phi)\geq -C$ can be rewritten as
\be
{1\over V}\int_X\phi\o_0^n
\leq J_\o(\phi)-K_\o(\phi)+C
\leq
{n\over n+1}I_\o(\phi)-K_\o(\phi)+C,
\ee
and expressing $I_\o(\phi)$ by its definition gives the desired result.
Q.E.D.

\medskip
The following lemma can be found in Rubinstein \cite{Ru}.
It is the analogue of the second inequality in (\ref{MAharnack})
for the complex Monge-Amp\`ere equation. As in that case,
it is Moser iteration, applied originally by Yau \cite{Y78}
to the Monge-Amp\`ere equation in his proof of the Calabi conjecture, and subsequently
used in the method of continuity for the case of positive $c_1(X)$
by \cite{T97, TZ}.
Its key ingredient is a uniform Sobolev constant,
which has now become available for the K\"ahler-Ricci flow
thanks to Theorem \ref{perelmanestimates}:

\begin{lemma}
\label{harnack3}
We have the following estimate along the K\"ahler-Ricci flow,
\be
{\rm osc}\,\phi \leq A\,{1\over V}\int_X\phi\,\o_0^n +B.
\ee
\end{lemma}

\noindent
{\it Proof.} Let $\psi={\rm max}_X\phi-\phi+1\geq 1$. Then for any $\al\geq
0$, we have
\bea
\int_X \psi^{\al+1}\o_\phi^n
&\geq&
\int_X \psi^{\al+1}(\o_\phi^n-\o_0\wedge\o_\phi^{n-1})
=
{i\over 2}(\al+1)\int_X \psi^\al
\,\pl\psi\wedge\bar\pl\psi\wedge\o_\phi^{n-1}
\nonumber\\
&=&
{i(\al+1)\over 2
({\al\over 2}+1)^2}
\int_X\pl(\psi^{{\al\over 2}+1})\wedge
\bar\pl(\psi^{{\al\over 2}+1})\wedge\o_\phi^{n-1}
\eea
Thus we obtain
\be
\label
{gradientalpha}
\|\na(\psi^{{\al\over 2}+1})\|
\leq {n({\al\over 2}+1)^2\over \al+1}\int_X\psi^{\al+1}\o_\phi^n,
\ee
and hence, in view of Theorem \ref{perelmanestimates}, (iv)
and setting $\beta\equiv {n\over n-1}>1$, $p=\al+2$,
\be
\bigg[\int_X \psi^{p\beta}\bigg]^{1\over \beta}\ \leq\ C\,p\int_X
\psi^p\o_\phi^n,\qquad p\geq 2.
\ee
Taking $p=2$ and iterating, $p\to p\beta\to\cdots p\beta^k$, it follows that
\bea
\log ||\psi||_{L^\infty(\o_\phi)}
\leq \sum_{k=1}^\infty {\log\, (2C\beta^k)\over 2\beta^k}+\log
||\psi||_{L^2(\o_\phi)}
= C_1+\log ||\psi||_{L^2(\o_\phi)}.
\eea
On the other hand, a simple Bochner-Kodaira argument shows that
\be
{1\over V}\int_X \psi^2 e^{f}\o_\phi^n
\leq {1\over V}\int_X |\na \psi|^2 e^{f}\o_\phi^n
+
({1\over V}\int_X \psi e^{f}\o_\phi^n)^2
\ee
since $R_{\bar kj}-g_{\bar kj}=\pl_j\pl_{\bar k}f$
(see \cite{Fu3}, or \cite{TZ, PSSW2}). By Theorem \ref{perelmanestimates}, (i),
the measures $e^{f}\o_\phi^n$ and $\o_\phi^n$ are equivalent. Together with
(\ref{gradientalpha})
with $\al=0$, we find that
\be
\bigg[{1\over V}\int_X \psi^2 \o_\phi^n\bigg]^{1\over 2}
\leq\ C\,(1+
{1\over V}\int_X \psi \o_\phi^n).
\ee
Finally, up to an additive constant, the expression on the right hand side
can clearly be bounded by ${\rm sup}_X\phi$
and
${1\over V}\int_X(-\phi)\o_\phi^n$. Both of these expressions are bounded,
up to an additive constant,
by ${1\over V}\int_X \phi\o_0^n$, the first because of the
plurisubharmonicity property $\o_0+{i\over 2}\ddb\phi>0$,
and the second in view of Lemma \ref{harnack2}. Q.E.D.

\medskip
\begin{lemma}
\label{convergencelemma}
Let $(X,\o_0)$ be a compact K\"ahler manifold, $\o_0\in c_1(X)$,
and consider the K\"ahler-Ricci
flow
(\ref{KRpotentialbis}), with initial value $c_0$
given by {\rm (\ref{initialc})}. If there exists a constant $C$ with
\be
{\rm sup}_{t\in[0,\infty)}{1\over V}\int_X\phi\,\o_0^n\ \leq \ C\ <\infty,
\ee
then the K\"ahler-Ricci flow converges exponentially fast
in $C^\infty$ to a K\"ahler-Einstein
metric.
\end{lemma}

\medskip
\noindent
{\it Proof.} By Lemma \ref{harnack3}, the hypothesis implies that the
oscillation ${\rm osc}\,\phi$
is uniformly bounded. But since $\int_X\o_\phi^n=1$ and $\|\dot\phi\|_{C^0}$
is bounded, we also
have
\be
0\ <\ C_1\ \leq\ {1\over V}\int_X e^{-\phi}\o_0^n\ \leq\ C_2
\ee
which implies that
\be
\label{normalizationestimates}
{\rm inf}_X\,\phi\leq -\log\,C_1,\quad {\rm sup}_X\,\phi\geq -\log \,C_2.
\ee
Combined with the bound for ${\rm osc}\,\phi$, this implies that
$\|\phi\|_{C^0}$ is bounded.
As in Yau's proof of the Calabi conjecture \cite{Y78}, the $C^0$ bound for
$\phi$ implies a uniform bound for $||\phi||_{C^k}$, for each $k\in {\bf
N}$.
Detailed derivations in the case of the K\"ahler-Ricci flow can be found in
\cite{PSS, Pa, Ca}.
This already implies that there exists subsequences of times $t_m\to+\infty$
with $g_{\bar kj}(t_m)$
converging to a K\"ahler-Einstein metric. The proof that the full flow
$g_{\bar kj}(t)$
converges exponentially fast to a K\"ahler-Einstein metric
is more involved, and actually makes use of bounds for the lowest eigenvalue
$\lambda_{\o(t)}$ for the operator $\bar\partial^\dagger\bar\partial$
introduced earlier in section \S 2.6.2.
It can be found
in
the second part of the proof of Lemma 6 in \cite{PSSW2}.

\subsubsection{Perelman's convergence theorem}

We are now in position to give a proof, different from
the earlier one in Tian-Zhu \cite{TZ}, of the following version of a result
announced by Perelman in private communications:

\begin{theorem}
\label{perelmanconvergence}
If $X$ admits a K\"ahler-Einstein metric $\o_{KE}$
and ${\rm Aut}^0(X)=0$, then for any initial metric $(g_0)_{\bar kj}$,
the K\"ahler-Ricci flow
converges to a K\"ahler-Einstein metric.
More generally, if ${\rm Aut}^0(X)\not=0$ and $G\subset {\rm Stab}(\o_{KE})$
is a closed subgroup whose centralizer in
the stabilizer ${\rm Stab}(\o_{KE})$ of $\o_{KE}$
is finite,
then the K\"ahler-Ricci flow converges to
a K\"ahler-Einstein metric for all $G$-invariant
initial metrics $(g_0)_{\bar kj}$.
\end{theorem}

\smallskip
\noindent
{\it Proof.}
Under the hypotheses of the theorem, the Moser-Trudinger inequality
(\ref{mosertrudinger})
holds. Since $F_\o(\phi)$ is decreasing under the flow,
it follows that $J_\o(\phi)$ is uniformly bounded.
The bound (\ref{mosertrudinger}) also shows that $F_\o(\phi)$ is bounded
from below.
By Lemma \ref{harnack1}, $K_\o(\phi)$ is then bounded.
By the second statement in Lemma \ref{harnack2},
${1\over V}\int_X\phi\,\o_0^n$ is bounded from above.
By Lemma \ref{convergencelemma}, the K\"ahler-Ricci flow
converges. Q.E.D.

\medskip

We note that by
an early result of Bando-Mabuchi \cite{BM},
K\"ahler-Einstein metrics are unique up to diffeomorphisms.
An extension of Theorem \ref{perelmanconvergence}
of the case of K\"ahler-Ricci solitons has been given
recently in \cite{TZ}, using an extension of works of Kolodziej \cite{K98}.
When restricted to the K\"ahler-Einstein case, \cite{TZ} also yields a proof
of Theorem \ref{perelmanconvergence}.

\subsubsection{Condition (B)}

There has been so far relatively few results on the convergence of the
K\"ahler-Ricci
flow when $c_1(X)>0$. The first such result is due to Hamilton \cite{H88},
who showed convergence for $X={\bf CP}^1$, when the initial metric
$(g_0)_{\bar kj}$
has positive curvature everywhere. This assumption was removed later by B.
Chow \cite{CH1},
who showed that, for any initial metric $(g_0)_{\bar kj}$, the curvature
eventually
becomes positive everywhere. Convergence in higher dimensions under the assumption of positive
bisectional curvature is treated in Chen and Tian \cite{CT1,CT2}.
This assumption is preserved under the K\"ahler-Ricci flow \cite{B87, Mok}, but it is
restrictive, since it implies that $X$ is holomorphically equivalent to
${\bf CP}^N$ \cite{Mori, SY}.
In fact, the arguments in \cite{CT1, CT2} rely on the existence a priori
of a K\"ahler-Einstein metric. More recently, convergence in the case of
toric
varieties with vanishing Futaki invariant has been established by Zhu
\cite{Zh},
and, as we have seen above, for general manifolds manifolds $X$ under the
assumption
that $X$ admits a K\"ahler-Einstein metric (see Theorem \ref{perelmanconvergence}
above) or
K\"ahler-Ricci soliton
\cite{TZ}. Since toric manifolds with vanishing Futaki
invariant
are known to admit K\"ahler-Einstein metrics \cite{WZ}, all these results
turn out to require manifolds for which the existence of such metrics were known
in advance.

\smallskip
In view of the conjecture of Yau, the convergence of the
K\"ahler-Ricci flow
should be tied with the stability of $X$ in the sense of GIT. We describe
next
some recent results in this direction, where the stability conditions
involved are the conditions (B) and (S) described in section \S 2.
The following theorem was proved in \cite{PS06}:

\begin{theorem}
\label{conditionBtheorem}
Let $(X,\omega)$ be a compact K\"ahler manifold with $c_1(X)>0$.
Then the K\"ahler-Ricci flow converges if and only if the Riemann
curvature tensor is bounded along the flow, the K-energy
is bounded from below,
and Condition (B) holds.
\end{theorem}

We have stated the theorem as it was stated in \cite{PS06}.
But clearly, the condition (B) can be replaced by the weaker condition
${\rm (B)^*}$ defined in section \S 2.6.1, since in the proof of Theorem
\ref{conditionBtheorem}, only metrics along the K\"ahler-Ricci flow are
considered, and their $K$-energy decreases along the flow and hence is bounded above.

\medskip

We discuss briefly the role of Condition (B) in the proof of Theorem
\ref{conditionBtheorem}.
Note that the assumption of uniform boundedness of the curvature, together
with
the fixed volume and the boundedness of the diameter by Perelman's results,
implies already that there exist
subsequences $g_{\bar kj}(t_m)$ which converge in the sense of
Cheeger-Gromov,
that is, after suitable reparametrizations depending on $t_m$.
The issue is full convergence, and more important, the convergence of
$g_{\bar kj}(t)$
as a sequence of tensors, point by point on $X$.

\smallskip

Now the assumptions that the K-energy is bounded from below and the
Riemann curvature tensor is bounded imply that
$\|\dot g_{\bar kj}\|_{C^k}=\|R_{\bar kj}-g_{\bar kj}\|_{C^k}$
tends to 0 for all $k$, but where norms are taken with respect to the
evolving
metric $g_{\bar kj}(t)$. The main step is then to show that the metrics
$g_{\bar kj}(t)$ are uniformly equivalent. For this, according to a lemma of
Hamilton \cite{H82}, it suffices to establish an exponential
decay for $\|\dot g_{\bar kj}\|_t$.
The key starting point is the following differential inequality
(\cite{PS06}, eq. (3.7)) for $Y(t)$
defined by
\be
Y(t)=\int_X |\na f|^2\o_\phi^n
\ee
where $f$ is the K\"ahler-Ricci potential defined by (\ref{riccipotential}),
\bea
\label{Y}
\dot Y(t)
&\leq &
-2\lambda_t\,Y(t)+2\lambda_t\, {\rm Fut}(\pi_t(\nabla^jf))
-
\int_X|\nabla f|^2(R-\mu n)\o^n\nonumber\\
&& \quad - \int_X\na^j f\overline{\na ^kf}(R_{\bar kj}-g_{\bar
kj})\o_\phi^n.
\eea
Here $\pi_t$ is the orthogonal projection of
$(1,0)$-vector fields onto the space of holomorphic vector fields,
${\rm Fut}(V)$ is the Futaki invariant acting on $V\in
H^0(X,T^{1,0}(X))$, and $\lambda_t$ is the lowest strictly positive
eigenvalue of the $\bar\partial^\dagger\bar\partial$ operator on vector fields. Under our assumptions,
all the terms on the right
hand side of (\ref{Y}) tend to $0$ except for the term $-2\lambda_t
Y(t)$. To obtain exponential decay for $Y(t)$, we shall produce a strictly
positive lower bound for $\lambda_t$.

\smallskip

This follows from Condition
(B): assume otherwise. Since the curvatures, volume, diameter of the
metrics $g(t)\equiv g_{\bar kj}(t)$ are bounded above and their
injectivity radii bounded from below, by the Cheeger-Gromov-Hamilton
compactness theorem \cite{C,G, H94}, a subsequence $F_{t_j}^*(g(t_j))$
converges in $C^\infty$, after suitable reparametrizations
$F_{t_j}$. The eigenvalues of the $\bar\partial^\dagger\bar\partial$ operator with
respect to the metric $F_{t_j}^*(g(t_j))$ and the almost-complex
structure $F_{t_j}(J)$ are the same as $\lambda_{t_j}$. By going
to a subsequence if necessary, we may assume that $F_{t_j}(J)$ converges
to an almost-complex structure $J_\infty\in \overline{Diff(X)\cdot J}$.
But $\lambda_t\to 0$, and thus ${\rm dim}\,H^0(X,T_{J_\infty}^{1,0}(X))
> {\rm dim}\,H^0(X,T_J^{1,0}(X))$, contradicting Condition (B).

\smallskip

Once the exponential decay of $Y(t)=\|\nabla f\|_{(0)}^2$ is
established, repeated applications of Bochner-Kodaira formulas
show that $\|\nabla f\|_{(s)}$ converge exponentially to
$0$ for all $s$. This implies the exponential convergence to $0$
of $\|\dot g_{\bar kj}(t)\|_{C^0}$, in the $g_{\bar kj}(t)$ norm,
which implies the uniform equivalence of all metrics $g_{\bar kj}(t)$,
by the lemma of Hamilton \cite{H82}. Once the uniform equivalence
of $g_{\bar kj}(t)$ is established, the proof of $C^\infty$ convergence is
easy.

\smallskip
The method of integral estimates has also been applied successfully by Hou and Li
\cite{HouLi} to boundary value problems for real Monge-Amp\`ere equations.

\subsubsection{Condition (S)}

One advantage of Condition (B) is that it is clearly
a necessary condition for the orbit
of $J$ to be a Hausdorff point in the moduli space
of orbits of almost-complex structures.
However, the assumption of uniform boundedness of the Riemann
curvature tensor along the K\"ahler-Ricci flow in
Theorem \ref{conditionBtheorem} is very restrictive.
The following theorem due to \cite{PSSW2}
eliminates the curvature assumption
completely, by replacing
Condition (B) by the closely related Condition (S):

\begin{theorem}
\label{conditionStheorem}
{\rm (i)} If ${\rm inf}_{\omega\in \pi c_1(X)}K_{\o_0}(\omega)>-\infty$,
and Condition (S) holds, then the K\"ahler-Ricci flow $g_{\bar kj}(t)$
converges exponentially fast to a K\"ahler-Einstein metric.

{\rm (ii)} Conversely, if the metrics $g_{\bar kj}(t)$ converge in
$C^\infty$
to a K\"ahler-Einstein metric, then the above two conditions
are satisfied.

{\rm (iii)} In particular, if $g_{\bar kj}(t)$ converge in $C^\infty$,
then the convergence is exponential.

\end{theorem}

\medskip
We sketch the proof of the convergence of the flow
under the assumptions in the theorem. It depends on the following criterion
for the convergence
of the flow \cite{PSSW2}:

\medskip
\begin{lemma}
\label{convergenceintegralcriterion}
If the following inequality is satisfied,
\be
\int_0^\infty \|R-n\|_{C^0}dt \ <\ \infty,
\ee
then the K\"ahler-Ricci flow converges exponentially fast in $C^\infty$
to a K\"ahler-Einstein metric.
\end{lemma}

Indeed, by definition of the K\"ahler-Ricci flow,
$\pl_t (\o_\phi^n/\o_0^n)=-(R-n)$, and the condition
in the lemma implies that the volume forms $\o_\phi^n$ are all equivalent in size to
the volume form $\o_0^n$.
But then $\phi=-\log(\o_\phi^n/\o_0^n)+\dot\phi-f_0$ is bounded in $C^0$,
and we can apply Lemma \ref{convergencelemma}.

\smallskip

Returning to the proof of the theorem proper, the first important step
is to show that, under the sole assumption that the K-energy is bounded
from below, we have
\be
\|R(t)-n\|_{C^0}\rightarrow 0\quad {\rm as}\ t\to+\infty.
\ee
This is the analogue for the K\"ahler-Ricci flow of an estimate
established by Bando \cite{B} for the method of continuity.
Combined with the inequality (\ref{Y}),
this leads to the following
key difference-differential
inequality for $Y(t)=\|\na f\|_{L^2}^2$,
\be
\label{differencedifferential}
\dot Y(t)
\leq -2\lambda_t\, Y+\epsilon
\prod_{j=0}^N Y(t-j)^{\delta_j\over 2},
\qquad
\sum_{j=1}^N\delta_j=2,\quad \delta_j>0,
\ee
for a fixed integer $N$,
any given $\epsilon>0$, and $t\in [T_\epsilon,\infty)$,
for a suitable $T_\epsilon$ large enough.
The point is that a difference-differential inequality
of this form can still guarantee the exponential decay
of $Y(t)$. In fact, let $Z(t)=Z_0\,e^{-\mu (t-T_\epsilon)}$,
for some $0<\mu<1$,
and take $\epsilon={1\over 2}{\rm inf}_t\lambda_t>0$.
We claim that, for $Z_0$ large enough and $\mu$ small enough,
$Z(t)$ is a barrier for $Y(t)$. Indeed, write $T_\e=0$ for simplicity,
and let $Z_0$ be chosen so that $Y(t)<Z(t)$ for $t\in [0,N]$.
We claim that $Y(t)<Z(t)$ for all $t$. Otherwise, let $T$ be the first
time with $Y(T)=Z(T)$. Then
\bea
\dot Y(t)-\dot Z(t)
&\leq& -2\e Y(t)+\mu Z(t)+\e \prod_{j=0}^NY(t-j)^{\delta_j}
\nonumber\\
&\leq&
Z_0\bigg\{(\mu-2\epsilon)e^{-\mu t}+\e\,e^{-\mu
t+\mu\sum_{j=0}^Nj\delta_j}\bigg\}
<0
\eea
for $t\in [0,T]$. But then $Y(T)<Z(T)$, which is a contradiction.
Thus we do have $Y(t)\leq Z_0e^{-\mu t}$ for all $t\geq 0$.
Combined with Futaki's $L^2(e^f\o^n)$ Poincare inequality
\cite{Fu3} and Perelman's non-collapse theorem,
the exponential decay of $Y(t)=\|\na f\|_{L^2}^2$
can be shown to imply an exponential decay for $\|f\|_{C^0}$,
and ultimately for $\|\Delta f-n\|_{C^0}=\|R-n\|_{C^0}$,
and the theorem follows.

\medskip
Several extensions and applications of Theorems \ref{conditionBtheorem}
and \ref{conditionStheorem} have been obtained in \cite{PSSW3}.
In particular, the condition that the K-energy be bounded from below
in Theorem \ref{conditionBtheorem} can be weakened to just the vanishing
of the Futaki invariant. It is also shown in \cite{PSSW3}
that, under the assumption that the K energy is bounded from
below (or just the vanishing of the Futaki invariant
when ${\rm dim}\,X \leq 2$), if the initial metric has positive bisectional curvature,
then the K\"ahler-Ricci flow converges to a K\"ahler-Einstein metric. We stress
that the convergence of the K\"ahler-Ricci flow has been treated in \cite{CT1, CT2},
but the arguments there rely on the existence of a K\"ahler-Einstein metric,
and thus on the solution of the Frankel conjecture. Here it is essential that
the arguments do not assume the a priori existence of such a metric,
in order to be viewed as a progress in the problem of giving an independent
proof of the Frankel conjecture by flow methods.

\subsubsection{Multiplier ideal sheaves}

We have seen in section \S how multiplier ideal sheaves arise from the continuity
method for the complex Monge-Amp\`ere equation. Here we discuss the adaptation of such
ideas to the context of the K\"ahler-Ricci flow.

\begin{theorem}
\label{multipliertheorem}
Consider the K\"ahler-Ricci flow {\rm (\ref{KRpotentialbis})} on a compact K\"ahler manifold
$(X,\o_0)$, with $\o_0\in c_1(X)$, and initial value $c_0$ given by {\rm (\ref{initialc})}.
Then the K\"ahler-Ricci flow converges if and only if there exists $p>1$
so that
\be
{\rm sup}_{t\geq 0}{1\over V}\int_X e^{-p\phi}\o_0^n <\infty.
\ee
The convergence is then in $C^\infty$ and exponentially fast.
\end{theorem}

This theorem is proved in \cite{PSS}, using a priori
estimates of Kolodziej \cite{K98, K03}
for Monge-Amp\`ere equations with $L^p$ right hand sides,
$p>1$. (In \cite{PSS}, the proof of the full
convergence of the flow, by opposition to convergence only along a
subsequence
of times $t_m\to+\infty$, was given only under the assumption that
${\rm Aut}^0(X)=0$.
But we now have at our disposal Lemma \ref{convergencelemma},
based on the recent improvements in \cite{PSSW2} for convergence arguments,
and this technical assumption can be removed.)
Clearly, Theorem \ref{multipliertheorem}
can be restated in terms of a version of multiplier
ideal sheaves: let ${\cal J}^p$ be the sheaf whose stalk ${\cal J}_z^p$
at $z\in X$ is defined by
\be
\label{stalk1}
{\cal J}^p_z
=
\{f;\quad \exists U\ni z,\ f\in{\cal O}(U)\quad
{\rm sup}_{t\geq 0}\int_Ue^{-p\phi}|f|^2\o_0^n<\infty
\ \}
\ee
Then the necessary and sufficient condition for the convergence of the
K\"ahler-Ricci flow is that there exists $p>1$ so that ${\cal J}^p$ admits
the global section $1$.

\smallskip
An application of Theorem \ref{multipliertheorem} to del Pezzo
surfaces can be found in \cite{Heier}.

\smallskip

Theorem \ref{multipliertheorem} implies the following weaker, but simpler
statement. First note that for any $p>1$ we
have
\be
{\rm sup}_{t\geq 0}{1\over V}\int_X e^{-p\phi}\o_0^n =\infty
\ \Rightarrow\
{\rm sup}_{t\geq 0}{1\over V}\int_X e^{-p(\phi-{1\over
V}\int_X\phi\o_0^n)}\o_0^n =\infty.
\ee
In fact, replacing $\phi\to \phi-{\rm sup}_X\phi$
in the integrals on the left hand side only increases their sizes,
in view of the estimate (\ref{normalizationestimates}). But ${\rm sup}_X\phi\leq {1\over
V}\int_X\phi\o_0^n+C$ by $\o_0$-plurisubharmonicity,
hence the assertion.
Assume now that a K\"ahler-Einstein metric does not exist, so that
for all $p>1$, there exists a subsequence of times $t_m\to+\infty$.
Let $\psi$ be a weak limit of $\phi-{1\over V}\int_X\phi\o_0^n$.
Then the multiplier ideal sheaf ${\cal I}(p\psi)$ defined as in
(\ref{multiplierideal})
defines a proper, coherent sheaf, with
$H^q(X,K_X^{-[p]}\otimes {\cal I}^p)=0$,
$q\geq 1$.

\medskip
This last statement has been recently strengthened by Rubinstein \cite{Ru}
to the same range $p\in ({n\over n+1},\infty)$ that Nadel \cite{Na}
established for the multiplier ideal sheaves obtained from the method of
continuity.
The point is that, together with Lemmas \ref{harnack2},
\ref{harnack3} and the following Lemma \ref{nadellemma1},
we have now the analogues in the case of the K\"ahler-Ricci flow
of all the ingredients required for Nadel's arguments in the method of
continuity, namely (\ref{MAharnack}) and Lemma \ref{nadellemma}:

\begin{lemma}
\label{nadellemma1}
The exact same statement as in Lemma \ref{nadellemma} holds,
with $\phi$ the solution of the K\"ahler-Ricci flow
{\rm (\ref{KRpotentialbis})}.
\end{lemma}

\medskip
\noindent
{\it Proof.} As in the proof of Lemma \ref{nadellemma},
the concavity of the logarithm implies
\be
C>{1\over V}\int_X e^{(1-p)\phi-f_0+\dot\phi}\o_\phi^n
\geq
{1\over V}
\int_X ((1-p)\phi-f_0+\dot\phi)\o_\phi^n+p {1\over
V}\int_X\phi\o_0^n,
\ee
and thus, by Perelman's uniform bound for $|\dot\phi|$ and Lemma \ref{harnack2},
\be
p {1\over V}\int_X\phi\o_0^n
\leq (1-p){1\over V}\int_X(-\phi)\o_\phi^n+C
\leq n(1-p){1\over V}\int_X\phi\o_0^n+C.
\ee
and Lemma \ref{nadellemma1} follows as before \cite{Ru}. Q.E.D.

\medskip
Assume now that $X$ does not admit a K\"ahler-Einstein metric.
Then ${1\over V}\int_X\phi\o_0^n$ must diverge to $+\infty$
for some subsequence of times $t_m\to+\infty$,
for otherwise Lemma \ref{convergencelemma} would imply that
the K\"ahler-Ricci flow converges to a K\"ahler-Einstein metric.
Thus the integrals in (\ref{nadelcondition}) must diverge to $\infty$
for some $p\in ({n\over n+1},1]$. If we let $\psi$ be an $L^1$ limit
point of $\phi-{1\over V}\int_X\phi\o_0^n$, then ${\cal I}(p\psi)$
provides the desired coherent, acyclic multiplier ideal sheaf.

\newpage
\section{General $L$: Energy Functionals and Chow Points}
\setcounter{equation}{0}

With this section, we begin the description of results which address, at least partially
in some way, the eventual sufficiency of stability conditions for constant scalar curvature
metrics in a general K\"ahler class $L$, which is not necessarily $K_X^{-1}$.
A first class of results links directly the energy functionals $K_\o(\phi)$ and
$F_\o^0(\phi)$ to the Chow point. That there should be some relation is to some extent
already built into the notion of stability: we have seen that $K_\o(\phi)$ and
$F_\o^0(\phi)$ are just the conformal changes of metrics in the very Deligne pairings
which define the line bundles $\eta_K$ and $\eta_{Chow}$ over the Hilbert scheme
giving the notions of K and Chow-Mumford stability. But the direct relations which we
describe below are much more precise, and one can hope that, combined with some suitable
$k\to\infty$ limiting process, they may eventually allow to deduce the asymptotic growth of
$K_\o(\phi)$ and $F_\o^0(\phi)$ on ${\cal K}$ from stability conditions.

\subsection{$F_\o^0$ and Chow points}

We start with the case of $F_\o^0$. Here the basic result is the following theorem due
to Zhang \cite{Z96}, which links all three concepts of critical points for $F_\o^0(\phi)$,
balanced imbeddings, and Chow-Mumford stability:

\begin{theorem}
\label{zhangtheorem}
Let $\hat X\subset{\bf CP}^N$ be a smooth projective variety. Let $\o=\o_{FS}$,
$\phi_\sigma(x)=\log{|\sigma x|^2\over |x|^2}$, and view $F_\o^0(\phi_\sigma)$ as
a function of $\sigma\in SL(N+1)/SU(N+1)$. Then

\smallskip
{\rm (a)} $-F_\o^0$ is convex along one-parameter subgroups;

{\rm (b)} A point $\sigma_0$ is a critical point for $F_\o^0$ if and only if $\sigma_0(\hat X)$ is
balanced, in the sense c.f. (\ref{balanced}) that
\be
\int_{\sigma(\hat X)}{\bar x_\al x_\beta\over |x|^2}\sigma_{FS}^n \sim \delta_{\bar\al\beta}.
\ee

{\rm (c)} Define a norm $\| f\|$ on $f\in H^0(Gr, O(d))$ by
\be\label{znorm}
\log \|f\|^2
=
{1\over D}\int_{Gr}\log {|f(z)|^2\over |Pl(z)|^{2d}}\,\o_{Gr}^m
\ee
where $Gr=Gr(N-n-1,{\bf CP}^N)$ is the Grassmannian of $N-n-1$ planes in ${\bf CP}^N$,
$Pl$ is the Pl\"ucker imbedding, $\o_{Gr}$ is the Fubini-Study metric restricted to
$Gr$, and $D$ is its volume. Then we have
\be
\label{zhangformula}
-F_\o^0(\phi_\sigma)
=
\log {\|\sigma\cdot {\rm Chow} (\hat X)\|^2\over
\|{\rm Chow} (\hat X)\|^2}.
\ee

{\rm (d)} In particular, $\hat X$ is Chow-Mumford stable if and only if there exists
a unique $\sigma_0\in SL(N+1)/SU(N+1)$ with $\sigma_0(\hat X)$ balanced.
\end{theorem}

The statements (a) and (b) and their proofs have been given in section \S 5.
Zhang's original proof of (c) using Deligne pairings is in \cite{Z96}. A different
and perhaps simpler proof based on comparing the derivatives of both sides
along one-parameter subgroups can be found in \cite{PS03}.
We show now how (d) follows from (a), (b), and (c). The theorem
of Kempf and Ness \cite{KN} says that $\hat X$ is stable if and only if
$\|\sigma\cdot {\rm Chow}(\hat X)\|^2$ is a proper map from $SL(N+1)$
in to ${\bf R}$, in the sense that the inverse image of any compact subset
is compact.
Equivalently, $\hat X$ is stable if and only if
\be
\label{kempfness}
\log \|\sigma\cdot {\rm Chow}(\hat X)\|^2 \geq - C
\quad{\rm and}
\quad
{\rm lim}_{\sigma\to\infty}\log \|\sigma\cdot {\rm Chow}(\hat X)\|^2=\infty
\ee
(where the latter statement means that for any constant $M$, there exists
a compact set $K\subset SL(N+1)$ so that $\log \|\sigma\cdot {\rm Chow}(\hat X)\|^2\geq M$
for $\sigma\notin K$.)

\smallskip

Assume that $\hat X$ is Chow-Mumford stable.
In view of the above Kemp-Ness characterization of stability,
this implies that $\log \|\sigma\cdot {\rm Chow}(\hat X)\|^2$
must attain its mimimum, and $-F_\o^0$ must have a critical point $\sigma_0$.
This critical point is unique, since if there are two distinct critical points,
the convexity of $F_\o^0$ along the one-parameter subgroup joining them
would force $\log \|\sigma\cdot {\rm Chow}(\hat X)\|^2$ to be constant along
this geodesic, contradicting its properness. By (b),
a point $\sigma_0$ is a critical point if and only if $\sigma_0(\hat X)$
is balanced.

\smallskip

Conversely, assume that there exists a unique balanced, or equivalently,
a unique critical point $\sigma_0$. By the convexity of $-F_\o^0$ along
one-parameter subgroups, this point must be a minimum along any such
path. And since any point in $SL(N+1)$ can be joined to $\sigma_0$ by such
a path, $\sigma_0$ must be a minimum for $-F_\o^0$ on the whole of $SL(N+1)$.
It is a strict minimum since the existence of another minimum would contradict
the uniqueness assumption of balanced points. It follows that
$\log \|\sigma\cdot {\rm Chow}(\hat X)\|^2$ is bounded from below, and
must tend to $\infty$ when restricted to any one-parameter subgroup.
This implies that the numerical invariant along any one-parameter subgroup
must be strictly positive, and, in view of the Hilbert-Mumford criterion,
the variety $\hat X$ must be stable.

\subsection{$K_\o$ and Chow points}

The analogue for the K-energy of the formula of the previous section
is the following \cite{PS03}:

\begin{theorem}
Let $X\subset{\bf CP}^N$ be as in Theorem \ref{zhangtheorem}. Then
\bea
\label{GITmabuchi}
&&
K_{\omega}(\phi_{\sigma})
+{1\over V}
\langle [Y_s],\Phi_{\sigma}\sum_{i=0}^{m-1}
\omega_Z^i\sigma^*\omega_Z^{m-1-i}\rangle
-{D\over V}
{m\, {\rm deg}(Y_s)\over m+1}
{\rm log}{||\sigma\cdot {\rm Chow}(X)||^2\over ||{\rm Chow}(X)||^2}
\nonumber\\
&&
\qquad\qquad
={D(m+2)(d-1)\over V(m+1)}
{\rm log}{||\sigma\cdot {\rm Chow}(X)||_{\#}^2
\over ||{\rm Chow}(X)||_{\#}^2}
\eea
\end{theorem}

\medskip
Here $Z=\{w\in Gr(N-n-1,{\bf CP}^N);\ w\cap \hat X\not=0\ \}$ is the
Chow variety of $X$. In general, $Z$ is a singular variety, and we denote
by $Y_s=\{w\in Gr(N-n-1,{\bf CP}^N);\ \# (w\cap \hat X) >1\ \}$ the
subvariety where it is singular. We denote by $\o_Z$ the restriction of
the Fubini-Study K\"ahler form $\o_{Gr}$ on the Grassmannian to the
regular part
$Z\setminus Y_s$ of $Z$, and by
$[Y_s]$ the current of integration on $Y_s$.
Let $P\ell:Gr(N-n-1,{\bf CP}^N)\to {\bf P}(\wedge ^{N-n}{\bf CP}^{N+1})$
be the Pl\"ucker imbedding, and set
for each
$\sigma\in GL(N+1)$, $\phi_\sigma(z)=\log {|\sigma z|^2\over |z|^2}$,
$\Phi_\sigma(z)=\log{|P\ell(\sigma z)|^2\over |P\ell(z)|^2}$. The norm
$\|\cdot\|$ is the norm defined in Theorem \ref{zhangtheorem}, while
$\|\cdot\|_{\#}$
is a degenerate semi-norm defined for $f\in H^0(Gr, O(d))$ as follows
\bea
\label{sharpnorm}
\log||f||_{\#}^2
&=&{m+1\over (m+2)(d-1)}{1\over D}
\int_Z\log \bigg({\o_{Gr}^m\wedge \pl\bar\pl{|f(z)|^2\over
|P\ell(z)|^{2d}}\over \o_{Gr}^{m+1}}\bigg)\o_{Gr}^m
\nonumber\\
&&
+{d-m-2\over (m+2)(d-1)}{1\over D}\int_{Gr}\log {|f(z)|^2\over
|P\ell(z)|^{2d}}\o_{Gr}^{m+1},
\eea
with $m=(N-n)(n+1)-1$ and $D$ the dimension and the volume of the
Grassmannian.
The various ingredients in the formula (\ref{GITmabuchi}) have an
interesting interpretation.
We can define a K-energy $K_{\o_Z}(\Phi_\sigma)$ associated to the regular
part $Z\setminus Y_s$
of the variety $Z\subset Gr$ just like the K-energy $K_\o(\phi_\sigma)$
for the variety $X\subset{\bf CP}^N$. Remarkably, an exact Radon transform
argument shows that $K_{\o_Z}(\Phi_\sigma)=K_\o(\phi_\sigma)$.
The left hand side of the equation (\ref{GITmabuchi}) can then be viewed as
a notion of K-energy associated to the full variety $Z$, and it is this
modified K-energy
associated to a singular variety which satisfies an identity analogous to
the identity for $F_\o^0(\phi_\sigma)$ stated in Theorem \ref{zhangtheorem}.

\smallskip
A formula for the K-energy in terms of a Quillen norm, up to a bounded
error,
has also been obtained by
Tian \cite{T94}. The Futaki ivariant and $K$-energy for hypersurfaces have been
evaluated by Lu \cite{Lu99, Lu04}, Yotov \cite{Yo}, and in \cite{PS04b}.

\smallskip
We note that (\ref{GITmabuchi}) can be interpreted as
a degenerate norm on a suitable line bundle over the Hilbert scheme.
More precisely, let $\eta_{Chow}$ be the
Chow
line bundle defined before, and define another line bundle
$\eta_{Chow_s}$  by associating to each variety $X$ the Chow line
of
$Y_s \sub {\bf P}(\wedge ^{N-n}{\bf CP}^{N+1})$,
equipped with the corresponding norm $\|\cdot\|_s$ defined as in
(\ref{znorm}). Then K-energy restricted to $SL(N+1)$ orbits
is the change in
the norm $(\|\cdot\|\otimes \|\cdot\|_{\#})\otimes
\|\cdot\|_s$
on
the line bundle
\be
\eta_{Chow}\otimes \eta_{Chow_s}
\ee
Besides this formal application, it would be very
interesting to explore the full consequences of
the formula (\ref{GITmabuchi}), since it relates several basic objects,
namely the K-energy, the Chow point, the singular locus $Y_s$,
and a degenerate semi-norm.

\newpage
\section{General $L$: The Calabi Energy and the Calabi Flow}
\setcounter{equation}{0}

We have seen that in the case $L=K_X^{-1}$, the metrics of constant scalar curvature
metrics are K\"ahler-Einstein metrics, and from the point of view of geometric flows,
the problem can be reduced to the issue of convergence of the K\"ahler-Ricci flow.
For general $L$, we need to deal with the full 4th-order equation, the natural
parabolic version of which is the Calabi flow, which we discuss briefly in this
section.

\subsection{The Calabi flow}
Let $L\to X$ be a positive line bundle over a compact complex manifold $X$,
and let $\o_0={i\over 2}g_{\bar kj}^0 dz^j\wedge d\bar z^k$
be a K\"ahler form in $c_1(L)$.
The Calabi flow is the 4th-order flow defined by the following equation
\be
\dot g_{\bar kj}=\pl_{\bar k}\pl_jR,\qquad g_{\bar kj}(0)=g_{\bar kj}^0.
\ee
Clearly, it preserves the K\"ahler class of $\o_0$. If we write then $g_{\bar kj}=
g_{\bar kj}^0+\pl_{\bar k}\pl_j\phi$,
the Calabi flow is equivalent to the following flow for $\phi$,
\be
\dot\phi=R-\bar R.
\ee
In view of the fact that the variational derivative of the $K$-energy is $-(R-\bar R)\o_\phi^n$,
we see that the Calabi flow is just the gradient flow of the K-energy,
and that
\be
{d\over dt}K_{\o_0}(\phi)=-{1\over V}C(\phi),
\ee
where $C(\phi)$ is the Calabi functional, defined by
\be
\label{calabifunctional}
C(\phi)=\int_X |R(\o_\phi)-\bar R|^2\o_\phi^n.
\ee

Clearly, metrics of constant scalar curvature are the minima of the
Calabi functional $C(\phi)$. However, the other critical points of
$C(\phi)$, called ``extremal metrics'', are also of considerable
interest. Now a straightforward calculation gives the following
variational formula for $C(\phi)$
\be
\delta C=-\int_X
\delta\phi\, (\Delta^2R+|\nabla R|^2+R_{\bar kj} \nabla^{\bar
k}\nabla^jR)\,\o_\phi^n.
\ee
Applying the Bianchi inequality, we can
easily verify that
\be
\Delta^2R+|\nabla R|^2+R_{\bar kj}
\nabla^{\bar k}\nabla^jR=\nabla^{\bar l}\nabla^{\bar q}\nabla_{\bar
q}\nabla_{\bar l}R.
\ee
Thus the extremal metrics are given by the
following equation
\be
\label{extremalmetric}
\nabla_{\bar
q}\nabla_{\bar l}R=0
\ee
which means exactly that $\nabla^jR$ is a
holomorphic vector field. It follows also immediately from the
variational formula for $C(\phi)$ that it decreases along the Calabi
flow,
\be {d\over dt}C(\phi)\ = \ -2\int_X
|\na_{\bar q}\na_{\bar l}R|^2\o_\phi^n\ \leq\ 0.
\ee
The operator $f\rightarrow \na_{\bar q}\na_{\bar l}f$
mapping functions to symmetric two-tensors the Lichnerowicz operator
${\cal D}$ which we have already encountered in section \S 7.2.
Its central role in the problem of canonical metrics had been stressed by
Calabi \cite{Ca}.

\medskip
When ${\rm dim}\,X=1$, the Calabi flow has been shown to converge
to a metric of constant scalar curvature
by Chrusciel \cite{Chr}, Chen \cite{Ch01}, and Struwe \cite{St}.
In higher dimensions,
both the long-time existence and the convergence of the Calabi flow
are open problems. There has however been several recent progresses.
In \cite{ChHe}, the flow was shown to exist as long as the Ricci curvature
stays uniformly bounded. In \cite{TW}, the flow was shown to converge if
the Calabi functional is initially small enough, and either $c_1(X)=0$,
or $c_1(X)<0$ and $c_1(L)$ is sufficiently close to $c_1(X)$. In the case of
ruled surfaces, the long-time existence for initial metrics given by the momentum
construction has been established in \cite{GuD1} under the assumption of existence
of an extremal metric, and in general in \cite{Szea}.

\subsection{Extremal metrics and stability}

We have seen that the extremal metrics are the critical points of the Calabi energy,
while the metrics of constant scalar curvature are its minima.
Just as in the case of constant scalar curvature metrics,
the existence of extremal metrics is expected to be equivalent to
some suitable form of stability. In \cite{Sze},
a notion of $K$-stability of a line bundle $L\to X$,
relative to a maximal torus of automorphisms, is introduced
and conjectured by Szekelyhidi to be equivalent to the existence of an extremal
metric in $c_1(L)$. We give here a brief description of these ideas.

\v
It is well known that the map $J\mapsto R(g_J)-\overline R$ is a moment map
for the action of the symplecto-morphism group of $(X,\o)$ on the
space ${\cal J}$, consisting of integrable complex structures on $X$
compatible with the symplectic form $\o$ (here $g_J$ is the K\"ahler
form $g_J(v,w)=\o(v,Jw)$). Thus the Calabi function may be viewed as
the norm squared of the moment map and, in order to gain some
insight into the nature of its critical points, Szekelyhidi first examines
the finite-dimensional picture: Let $L\ra (X,\o)$ be as above and
let $K$ be a compact Lie group acting biholomorphically on $L\ra X$
and preserving $\o$. A moment map for the action of $K$ on $(X,\o)$
is a smooth $K$-equivariant map $\m: X\ra Lie(K)^*$ with the property
\be
W(\langle \m,\xi\rangle) \ = \ \o(V_\xi,W)
\ee
for every smooth vector field $W$ and every $\xi\in Lie(K)$ (here $V_\xi$ is the vector field on $X$
generated by $\xi$). Let $G$ be the complexification of $K$ and assume
$G$ acts as well on $L\ra X$ and that the action is compatible with that of $K$.
Recall that an orbit $Gx$ is stable if for every one-parameter
subgroup $\l(t)=\exp(t\al)$, the weight $F_x(\al)=F_x(\l)$ is positive. We say
that the orbit is polystable if $F_x(\l)\geq 0$ with equality only if $\l$
fixes $x$.
Then the Kempf-Ness theorem says that a $G$ orbit $Gx$ contains a zero
of $\m$ if and only if the orbit is polystable.
\v
The analogue of the Calabi functional in this finite-dimensional
setting is the function $c=\|\m\|^2:X\ra  [0,\i)$ (the norm is taken
with respect  to a fixed inner product $(\ , \ )$ on $Lie(K)$), and extremal
metrics correspond to critical points of the function $c(x)$. The
following generalization of the Kempf-Ness theorem is proved in \cite{Sze}:

\begin{theorem}
\label{Sz}
A point $x\in X$ is in the $G$-orbit of a critical point of $c$ if
and only if it is polystable relative to a maximal torus in $G_x$,
where $G_x$ is the stabilizer of $x$.
\end{theorem}
Let us explain the terminology of the theorem: Let $T\sub G_x$ be a
maximal torus and let $G_T$ be the connected component of the centralizer of $T$
(elements of $G$ which commute with all the elements of $T$).
Then there is a connected subgroup
$G_{T^\perp}\sub G_T$ which can be characterized as follows: it is isomorphic to $G/T$ under the map
$G_{T^\perp}\hookrightarrow G_T\ra G_T/T$, and it has  the property
that $( \al,\b)=0$ for all $\al\in Lie(T)$ and all  $\b\in Lie(G_{T^\perp})$.
We say that $x\in X$ is polystable relative to $T$ if it is polystable for the
action of $G_{T^\perp}$ on $(X,L)$, that is, if and only if
\be\label{extremal1} F_x(\al)\geq 0\ \  \hbox{for
all $\al\in Lie(G_{T^\perp})$ with equality if and only if
$\l_\al$ fixes $x$.}
\ee
One can rewrite condition (\ref{extremal1}) as a condition on the set of all
one-parameter subgroups $G_T$ as follows: consider the
linear map $Lie(T)\ra \R$ given by $\al\mapsto F_x(\al)$ where
$\al\in Lie(T)$, and let $\chi\in Lie(T)$ be
its dual.  Thus
$F_x(\al)=(\al,\chi)$
for all  $\al\in Lie(T)$. Let $F_{x,\chi}(\al)=F_x(\al)-(\al,\chi)$.
Then one shows that (\ref{extremal1}) can
be rewritten as follows
\be
\label{extremal2}
F_{x,\chi}(\al)\geq 0\ \ \ \hbox{for all
$\al\in Lie(K)\cap Lie(G_T)$ }
\ee
with equality iff and only if $\al$ fixes $x$.
\v

To generalize these notions to the infinite-dimensional setting, Szekelyhidi first constructs an
inner product on the space of test configurations: Let $V$ be a
a projective scheme and $L\ra V$ a very ample line bundle.
A $\C^\times$ action for $L\ra V$ is a homomorphism
$\al:\C^\times\ra \Aut(L\ra V)$. If $\al,\b$ are two such
$\C^\times$ actions, let $A_k,B_k$ be the infinitesimal generators
on $H^0(V,L^k)$ and define $\langle \al,\b\rangle$ by the equation
\be
\label{inner product}
{{\rm Tr} (A_kB_k)}\ - \ {{\rm Tr}(A_k){\rm Tr}(B_k)\over \dim\,H^0(V,L^k)}\ = \ \langle \al,\b\rangle k^{n+2}
\ + \ O(k^{n+1})
\ee
Next we construct the  {\it extremal} $\C^\times$ action $\chi$ as follows.
Let $T\sub\Aut(L\ra V)$ be a maximal torus. Define $\chi:\C^\times\ra T$
by requiring
\be
\label{chi}
F(\al)\ = \ \langle\chi,\al\rangle\ \ \ {\rm for\ all} \ \ \
\al:\C^\times\ra T\ .
\ee
\v
Now let $L\ra X$ be an ample line bundle over a smooth projective variety and
let $T\sub\Aut(L\ra X)$ be a maximal torus.
Let $\r:\C^\times\ra\Aut(\cL\ra\cX\ra\C)$ be a test configuration.
Let
$\cX^\times =\pi^{-1}(\C^\times)$
and define
a homomorphism $T\ra \Aut(\cL^\times\ra \cX^\times)$  by
\be
\label{non-compact action} t\cdot l_\tau\ = \ \r(\tau)t\r(\tau^{-1})l_\tau\ \ \ {\rm for\ all}\ \ \
l_\tau\in L_\tau
\ee
We say that $\r$ is compatible with $T$ if the action defined by
(\ref{non-compact action})
extends to an action $T\ra \Aut(\cL\ra\cX)$. If $\r$ is compatible with $T$,
we let $\ti\chi:\C^\times\ra\Aut(L_0\ra X_0)$ the be restriction of
$\chi:\C^\times\ra T\ra \Aut(\cL\ra\cX)$
to the central fiber, and we let $\ti\r:\C^\times\ra\Aut(L_0\ra X_0)$ be
the restriction of $\r$ to the central fiber.
\begin{definition}
Let $(X,L)$ be a polarized variety and $T\sub\Aut(L\ra X)$
a maximal torus. Let $\chi:\C^\times\ra T$ be the extremal
$\C^\times$ action. We say that $(X,L)$ is K-stable relative to $T$ if for
all test configurations $\r$ compatible with $T$, we have
\be F_{\ti\chi}(\ti\r)\ \equiv \ F(\ti\r)\ - \ \langle \ti\chi,\ti\r\rangle\ \geq\ 0
\ee
with equality if and only if the test configuration is a product.
\end{definition}

It is then conjectured in \cite{Sze}
that a polarized variety $L\to X$ admits an extremal metric if and only
if it is $K$-stable relative to a maximal torus.

\newpage
\section{General $L$: Toric Varieties}
\setcounter{equation}{0}

In the case of general $L$, perhaps the greatest advances in the direction of sufficiency
have taken place in the context of toric varieties, where the equations for the K\"ahler
potential can be re-expressed in terms of its Legendre transform, namely the symplectic
potential. They become then real equations, and the tools from convex analysis can
be brought to bear. There has been many remarkable developments in this direction, even a
perfunctory description
of which is beyond the scope of this paper. We shall limit ourselves to a few
words, mainly to provide references for further reading.

\subsection{Symplectic potentials}

The basic properties of K\"ahler forms and their scalar curvatures on toric varieties
have been worked out by Guillemin \cite{Gui} and Abreu \cite{Abreu}.
A brief summary is as follows.
Let $L\to X$ be a positive toric line bundle over a toric variety.
Then a dense orbit in $X$ of the $({\bf C}^\times)^n$ action
can be parametrized by $z=(z_1,\cdots,z_n)$, with
$z_i=\xi_i+i\eta_i$. If $\o$ is a K\"ahler form on the orbit which
is invariant under the $(S^1)^n$ subgroup of $({\bf C}^\times)^n$,
then $\o$ must be of the form
\be
\o={i\over 2}{\pl^2\phi\over\pl \xi_j\pl\xi_k}\,dz^j\wedge d\bar z^k
\ee
with $\phi$ a strictly convex
function of the variables $\xi_i$ alone, $\phi=\phi(\xi)$.
Associated to $\phi$ is its Legendre transform $u(x)$ defined by
\be
u(x)={\rm sup}_\zeta (\<x,\zeta\>-\phi(\zeta))=\<x,\xi\>-\phi(\xi),
\qquad x={\pl\phi\over\pl\xi},
\ee
which is a function on the polytope $P$ given by the image of the
moment map
\be
\mu\ :\ z\ \rightarrow\ x={\pl\phi\over\pl\xi}.
\ee
The polytope $P$ does not depend on the choice of K\"ahler form
$\o$ within $c_1(L)$.
The function $u(x)$ is called the symplectic potential. It is strictly convex.
The metric $\omega$ can be re-written in terms of the coordinates
$(x,\eta)$ as $u_{ij}dx^idx^j+u^{ij}d\eta^id\eta^j$, where $u_{ij}$
is the Hessian metrix of $u_{ij}$, and $u^{ij}$ is its inverse.
The scalar curvature becomes
\be
\label{abreu}
R(\o)=-{\pl^2 u^{ij}\over\pl x^i\pl x^j}.
\ee
The form $\o$ extends to a smooth K\"ahler form on $X$ if and only if
$u(x)$ satisfies
\be
\label{guillemin}
u(x)-\sum \delta_k(x)\log\,\delta_k(x)\,\in\, C^\infty(\bar P).
\ee
Here we have written $P$ as the intersection of half-spaces $\delta_k(x)>0$.

\subsection{K-stability on toric varieties}

The main concepts in $K$-stability can also be re-written completely explicitly
in terms of the symplectic potential $u$ and the polytope $P$. In particular,
test configurations correspond precisely to rational piecewise linear convex functions $f$ on $P$,
and the $K$-energy and the Futaki invariant are given by the following
functionals of $u$ \cite{D02},
\bea
F(f)&=&-\mu \int_P f\,dx+\int_{\partial P}f\,d\sigma\nonumber\\
K(u)&=&
-\int_P\log\,({\rm det}\,u_{ij})\,dx+F(u),
\eea
where $dx$ is the Lebesgue measure on $P$ and $d\sigma$ is the measure on $\partial P$
defined by the requirement that $d\sigma\wedge d\delta_k=\pm dx$
on the face defined by $\delta_k(x)=0$.
The problem reduces then to finding a smooth convex solution
in $P$ for the following 4-th order equation,
called Abreu's equation,
\be
-{\pl^2 u^{ij}\over\pl x^i\pl x^j}=\bar R
\ee
subject to the boundary condition (\ref{guillemin}), assuming that the Futaki-Donaldson
invariant is strictly positive for all test configurations. This remains a challenging
problem, but when ${\rm dim}\,X=2$, the following progresses have been made
by Donaldson:

\smallskip
(a) It has been shown in \cite{D02} that $K$-stability implies that the
$K$-energy is bounded from below;

(b) Interior a priori estimates for Abreu's equation have been obtained in \cite{D04a};

(c) A method of continuity has been developed in \cite{D06}.

\smallskip
\noindent
Generalized solutions, assuming the properness of the $K$-energy, are also studied in \cite{ZZ}.

\subsection{The $K$-unstable case}

In the $K$-unstable case, it has been shown by Szekelyhidi \cite{Szeb} that
there is a maximally destabilizing configuration, in the class of
functions which are in $L^2(P)$,
continuous and convex on the union of $P$ with its codimension $1$ faces,
and whose boundary values are in $L^1(\pl P)$. It is not yet known, and perhaps
not always the case that this configuration is piecewise linear.
If it is, he also shows how this maximally destabilizing
configuration would lead to a decomposition into semi-stable components
analogous to the Harder-Narasimhan filtration for an unstable vector bundle.
The Calabi functional is the analogue in this context of the Yang-Mills
functional.
In general, this is not known and not expected that this maximally destabilizing
should be piecewise linear. Thus it should be viewed as
a limit of test configurations \cite{Szeb}. An attractive scenario
for how to construct still a Harder-Narasimhan like filtration
is also given in \cite{Szeb}.

\newpage
\section{Geodesics in the Space ${\cal K}$ of K\"ahler Potentials}
\setcounter{equation}{0}

As we have seen in section \S 6.2,
Donaldson's infinite-dimensional GIT theory provides yet another possible approach
to the problem of constant scalar curvature metrics. In this approach,
geodesic rays in the space ${\cal K}$ of K\"ahler potentials play the role
of test configurations, and numerical invariants and stability conditions
are to be constructed from the asymptotic behavior of the $K$-energy
along geodesic rays.

\smallskip
The starting point for this approach is then the construction and regularity
of geodesics. We have seen earlier that this can be viewed as an existence and
regularity problem for a Dirichlet problem for a completely degenerate
Monge-Amp\`ere equation. There has been significant progress in this direction, some
of which we describe below.

\subsection{The Dirichlet problem for the complex Monge-Amp\`ere equation}

Let $L\to X$ be a positive line bundle over a compact complex manifold $X$
of dimension $n$. Let $h_0$ be a metric on $L$, with $\o_0=-{i\over 2}\ddb\log\,h_0>0$.
Let $A=\{w\in {\bf C}; e^{-T}<|w|<1\}$, and write $\Omega_0$ for $\o_0$,
viewed as a non-negative $(1,1)$-form
on $X\times A$. We consider the following Dirichlet problem
\be
\label{MADirichlet}
(\Omega_0+{i\over 2}\ddb \Phi)^{n+1}=0\ {\rm on}\ X\times A,
\qquad \Phi\big\vert_{X\times \pl A}=\Phi_b,
\ee
where $\Phi$ is continuous at $X\times\pl A$,
and the boundary $\pl A$ and the boundary data are specified in the following
two possible manners:

\bigskip
$\bullet$ {\it Geodesic segments}: Here another metric $h_1$
on $L$ with $\o_1=-{i\over 2}\ddb \log\,h_1>0$ is prescribed,
$T>0$ is finite, $\pl A\equiv\{w\in {\bf C}; |w|=1\ {\rm or}\ |w|=e^{-T}\}$,
and the boundary value $\Phi_b$ is given by
\bea
\Phi_b(z,w)&=&0\qquad {\rm for}\qquad |w|=1
\nonumber\\
\Phi_b(z,w)&=&\log {h_1(z)\over h_0(z)}\qquad {\rm for}\qquad |w|=e^{-T}.
\eea

\medskip
$\bullet$ {\it Geodesic rays}: Here $T=\infty$, and $\pl A\equiv \{w\in {\bf C};|w|=1\}$,
and the boundary value $\Phi_b$ is given by
\be
\Phi_b(z,w)=0\qquad {\rm for}\qquad |w|=1.
\ee

\subsection{Method of elliptic regularization and a priori estimates}

The theory of a priori estimates for elliptic Monge-Amp\`ere equations has been
developed by Calabi \cite{Ca}, Yau \cite{Y78}, and Aubin \cite{A}
in the compact case, and by Caffarelli,
Kohn, Nirenberg, and Spruck \cite{CKNS} in the case of domains in ${\bf C}^n$ with boundary.
The geodesics problem requires an extension to the case of degenerate equations on
complex manifolds with boundary. The following theorem was established in \cite{PS07a}.
The main interest lies in its formulation, otherwise it
is a generalization, requiring no new estimates, of results in Chen \cite{Ch00}.
As can be seen from the sketch of
the proof given below, it is a direct outgrowth of the estimates for Monge-Amp\`ere equations
obtained in \cite{CKNS, Y78, GuB, Ch00}:

\v

\begin{theorem}
\label{MADirichletMtheorem}
Let $\bar M$ be a compact complex manifold of dimension $m$,
with smooth boundary $\pl M$. Assume that $\bar M$ admits a K\"ahler form $\Omega$. Then
the Dirichlet problem for the completely degenerate Monge-Amp\`ere equation
\be
\label{MADirichletM}
(\Omega+{i\over 2}\ddb \Psi)^m=0\ {\rm on}\ M,
\qquad \Psi|_{\pl M}=0
\ee
for a function $\Psi$ which is $\Omega$-plurisubharmonic,
admits a unique $C^{1,1}$ solution.
\end{theorem}

\v
We sketch the proof of Theorem \ref{MADirichletMtheorem}. Consider the elliptic regularization
\be
\label{ellipticregularizationM}
(\Omega+{i\over 2}\ddb \Psi)^m=\epsilon\,\Omega^m\ {\rm on}\ M,
\qquad \Psi|_{\pl M}=0
\ee
for a constant $\epsilon>0$ which we shall ultimately let tend to $0$.
Since $\Omega$ is a K\"ahler form, and hence positive definite,
the function $\underline{\Psi}=0$ is a subsolution
of (\ref{ellipticregularizationM}) for any $0<\epsilon\leq 1$.
The elliptic regularization admits then a $C^\infty$ solution,
and it suffices to establish $C^2$ estimates for the solution
uniform in $\epsilon$ in order to obtain the theorem.

\smallskip

- The $C^0$ estimates follow as in the case of domains in ${\bf C}^m$
from the maximum principle \cite{CKNS},
once a subsolution has been constructed;

- In the case of domains in ${\bf C}^m$,
the estimates for the first-order derivatives can be reduced to the boundary,
by differentiating the equation along any constant, global vector
field, and applying the maximum principle. Since the solution
is bounded by the subsolution and a harmonic function
with same boundary condition, it follows that the boundary values
of its gradient must be bounded \cite{CKNS}.
This argument does not generalize to complex manifolds.
However, the blow-up arguments of Chen \cite{Ch00} still apply.

- The differential inequalities of Yau \cite{Y78} and Aubin \cite{A}
for the trace of the Hessian of the solution hold on general
K\"ahler manifolds.
Applying the maximum principle, they reduce the estimates
for the second derivatives
to the boundary. A difficult barrier argument gives next
these boundary estimates, in the case of strongly pseudoconvex domains
in ${\bf C}^n$ \cite{CKNS}. This barrier argument was subsequently extended
by Guan \cite{GuB} to the general case, assuming instead of strong
pseudoconvexity the existence of a subsolution with given boundary
values. Q.E.D.

\medskip
We observe that the $C^2$ estimates
give an upper bound for the eigenvalues of
$\Omega+{i\over 2}\ddb\Psi$. If the Monge-Amp\`ere equation is elliptic,
the determinant is bounded from below, and hence all the eigenvalues are
bounded from both above and below. The $C^3$ identity of Calabi and Yau
can then apply, reducing again the $C^3$ estimates to the boundary.
An ingenious argument is then provided in \cite{CKNS} for the logarithmic modulus
of continuity of the second derivatives, and hence the $C^3$ boundary
estimates. Alternative approaches have also been provided by Evans \cite{Evans}
and Krylov \cite{Kry}.

\smallskip
In the degenerate case, the $C^2$ estimates do not imply bounds from below for the eigenvalues
of the form $\Omega+{i\over 2}\ddb\Psi$, and we cannot go further and obtain estimates for
higher order derivatives. This is consistent with the optimal $C^{1,1}$
regularity of known examples of solutions of degenerate Monge-Amp\`ere equations \cite{Gam, Lempert}.
A partial regularity theory for completely degenerate complex Monge-Amp\`ere equations
has been put forth in \cite{CT}.

\bigskip
We return now to the problem proper of geodesics in the space ${\cal K}$
of K\"ahler potentials for a positive line bundle $L\to X$.
We consider first the case of geodesic segments, when $M=X\times A$
with $A$ an annulus,
the case of geodesic rays being somewhat different and treated separately
later. The following simple, but key lemma will allow
to reduce the original equation (\ref{MADirichlet}) with degenerate background
form $\Omega_0$
to the situation treated in Theorem \ref{MADirichletMtheorem}:

\begin{lemma}
There exists a function $\underline\Phi \in C^\infty(X\times \bar A)$ which
is a subsolution of the equation
(\ref{MADirichlet}), in the sense that
\be
\Omega_0+{i\over 2}\ddb\underline\Phi\,>\,0,
\qquad
\underline{\Phi}\big\vert_{X\times\pl A}=\Phi_b.
\ee
\end{lemma}

\medskip
\noindent
{\it Proof.} Let $\phi(z,t)=t\,\log\ {h_0(z)\over h_1(z)}$, and
$\underline{\tilde \Phi}(z,w)=\phi(z,t)$, $t=\log\,|w|$.
Then clearly $\underline{\tilde \Phi}\big\vert_{X\times\pl A}=\Phi_b$,
and we have
\be
\Omega_0+{i\over 2}\ddb
\underline{\tilde \Phi}(z,w)
=\pmatrix{(1-t)\o_0+t\o_1 & {i\over 4w}\pl_{\bar z}\log {h_0(z)\over h_1(z)}\cr
-{i\over 4\bar w}\pl_z\log {h_0(z)\over h_1(z)} & 0}
\ee
The upper left entry is uniformly strictly plurisubharmonic on $X$.
Thus, setting $\underline{\Phi}(z,w)=\underline{\tilde \Phi}(z,w)+f(w)$,
where $f(w)$ is the solution of the Dirichlet problem
$\Delta f=C$, $f\big\vert_{\pl A}=0$ for some constant $C>0$
large enough, we obtain the desired subsolution. Q.E.D.

\medskip
Let now $\Omega=\Omega_0+{i\over 2}\ddb\underline{\Phi}$,
$\Psi=\Phi-\underline{\Phi}$. The equation (\ref{MADirichlet}) for a function
$\Phi$ which is $\Omega_0$-plurisubharmonic, is then equivalent
to the equation in (\ref{MADirichletM}) on $\bar M=X\times\bar A$,
for a function $\Psi$ which is $\Omega$-plurisubharmonic.
Theorem \ref{MADirichletMtheorem} gives then the
following existence and $C^{1,1}$ regularity for geodesic segments which was
proved by X.X. Chen \cite{Ch00}, and which led to the more
general formulation provided in Theorem \ref{MADirichletMtheorem}:

\begin{theorem}
\label{chentheorem}
\cite{Ch00}
Let $L\to X$ be a positive line bundle over a compact complex manifold,
and $h_0, h_1$ $C^\infty$ metrics on $L$ with positive curvatures $\o_0,\o_1$.
Then the Dirichlet problem {\rm (\ref{MADirichlet})}
admits a unique $C^{1,1}$ solution. The solution is $C^1$ invariant,
and defines a $C^{1,1}$ geodesic joining $h_0$ to $h_1$ in the space ${\cal K}_k$
of K\"ahler potentials.
\end{theorem}

\bigskip
Next, we discuss the construction of geodesic rays. Recall that they are the analogues for ${\cal K}$
of the one-parameter subgroups in ${\cal K}_k$. Since one-parameter subgroups are essentially
the same as test configurations, the natural question which arises is whether one can
associate a geodesic ray to each test configuration.

\smallskip

In \cite{AT}, Arezzo and Tian showed that, given a test configuration
${\cal T}$ for $L\to X$ with smooth central fiber $X_0$,
then one can use the Cauchy-Kowalevska theorem to find families of local
analytic solutions to the geodesic equations near infinity.
The geodesic rays obtained in this manner are real-analytic, but
their origins cannot be prescribed. The condition that the central
fiber be smooth is also a severe restriction. In \cite{PS07a},
the following was established:

\begin{theorem}
Let $L\to X$ be a positive line bundle on a compact complex manifold,
and let $cT: \cL\to cX\to {\bf C}$ be a test configuration
for $L\to X$ in the sense of Definition \ref{}.
Let $h_0$ be any metric on $L$ with $\o_0=-{i\over 2}\ddb\log\,h_0>0$,
and consider the Dirichlet problem (\ref{MADirichlet})
in the geodesic ray case, that is, when $A=\{w\in {\bf C};0<|w|\leq 1\}=D^\times$,
and $\pl A\equiv \{w\in {\bf C};|w|=1\}$.

Let $p:\tilde{\cal X}\to{\cal X}\to\cal {\bf C}$ be any smooth,
$S^1$ equivariant resolution of ${\cal X}$.
Then the Dirichlet problem
admits a solution $\Psi:X\times D^\times\to {\bf R}$, with
$\Omega_0+{i\over 2}\ddb\Psi$ the restriction to $p^{-1}({\cal
X}_{\vert_{D^\times}})$
of a non-negative $(1,1)$ current $\Omega+{i\over 2}\ddb\Psi$
on $\tilde{\cal X}_D\equiv
p^{-1}({\cal X}_{\vert_{D}})$ satisfying
\be
(\Omega+{i\over 2}\ddb\Psi)^{n+1}=0 \quad {\rm on}\quad \tilde{\cal X}_{\vert_D}.
\ee
Here $\Omega$ is a smooth K\"ahler metric
on $\tilde{\cal X}_{\vert_D}$,
and $\Psi$ is a $C^{1,1}$ function.
\end{theorem}

\medskip
A key step in the proof of this theorem is the construction of the K\"ahler
form $\Omega$ on $\tilde{\cal X}_{\vert_D}$, the remaining part following
readily from
Theorem \ref{MADirichletMtheorem}.
This is accomplished by constructing a
line bundle ${\cal M}\to \tilde\cX$ with the properties that

\smallskip

(a) $p^*\cL^m\otimes {\cal M}\to\tilde \cX$ is a positive line bundle;

(b) ${\cal M}\big\vert_{\tilde\cX^\times}$ is trivial, in the sense
that ${\cal M}$ admits a meromorphic section which is holomorphic and nowhere vanishing on
$\tilde\cX^\times$.

\smallskip
The desired K\"ahler form can then be taken as the curvature of
$p^*\cL^m\otimes {\cal M}\to\tilde \cX$.

\subsection{Geodesics in ${\cal K}$ and geodesics in ${\cal K}_k$}

An important question in the problem of constant scalar
curvature metrics is to determine in what precise sense
Donaldson's infinite-dimensional GIT is the limit of GIT.
In this context, it would be particularly valuable to
realize geodesics in ${\cal K}$ as limits of geodesics in ${\cal K}_k$.
This can be viewed also as the natural next step in Yau's general
strategy of approximations by algebraic-geometric objects: the Tian-Yau-Zelditch
theorem says that ${\cal K}$ is the ``pointwise'' limit of ${\cal K}_k$,
and the natural next step is to understand the external geometry of ${\cal K}_k$
as $k\to\infty$. If geodesics in ${\cal K}$ can be approximated by geodesics
in ${\cal K}_k$, as is desirable from the point of view of GIT,
this would mean that the subspaces ${\cal K}_k$ satisfy a remarkable property,
namely that they become asymptotically geodesically flat as $k\to\infty$.
In this section, we describe in what precise sense the answer to this question
is indeed affirmative \cite{PS06a, PS07}.

\subsubsection{An Ansatz for geodesic approximations}

Our approach is based on the following general Ansatz \cite{PS06, PS06a, PS07a}.
Let $L\to X$ be a positive line bundle over a compact complex manifold $X$ of dimension $n$,
and let $h_0$ be a metric on $L$ with positive curvature $\o_0=-{i\over 2}\ddb\log\,h_0$.
Let $\us=\{s_\al\}_{\al_0}^{N_k}$ be an orthonormal basis for $H^0(X,L^k)$
with respect to the metric $h_0$ and the volume form $\o_0^n$.
For each $k>>1$, let $\lambda_\al^{(k)}$ be a sequence of real numbers,
$0\leq\al\leq N_k={\rm dim}\,H^0(X,L^k)$. Set
\be
\label{bergmanmetric}
\Phi_k(z,w)
=
{1\over
k}\log\,\sum_{\alpha=0}^{N_k}|w|^{2\lambda_\alpha^{(k)}}|s_\alpha(z)|^2h_0^k
-n{\log\,k\over k}.
\ee
The function $\Phi_k(z,w)$ is manifestly $\Omega_0$-plurisubharmonic
on $X\times D^\times$, where $D^\times=\{w\in {\bf C};0<|w|<1\}$, and
$\Omega_0$ is the form $\o_0$, viewed as a non-negative form on $X\times D^\times$.
Let
\be
\label{approximation}
\Phi(z,w)={\rm lim}_{\ell\to\infty}
\bigg[\,{\rm sup}_{k\geq \ell}\ \Phi_k(z,w)\,\bigg]^*,
\ee
where $u^*(z)\equiv {\rm lim}_{\epsilon\to 0}({\rm sup}_{|z-\zeta|<\epsilon}u(\zeta))$
denotes the upper semi-continuous envelope
of a function $u(z)$. Then $\Phi(z,w)$ is a $\Omega_0$-plurisubharmonic function
on $X\times D^\times$, and $\Phi_k(z,w)=\phi_k(z,\log |w|)-n{\log\,k\over k}$,
where
\be
\phi_k(z,t)
=
{1\over k}
\log\,
\sum_{\alpha=0}^{N_k} e^{2\lambda_\alpha^{(k)}t}|s_\alpha(z)|^2h_0^k
\ee
is a geodesic in the space ${\cal K}_k$ of Bergman potentials.
Geometrically, the expression (\ref{bergmanmetric}) has another nice
motivation. Let $\pi_*(L^k)$ be the direct image of $L^k$ over $A_T$,
that is, the vector bundle over $A_T$ whose fiber at each $w\in A_T$
is the vector space $H^0(X,L^k)$. Holomorphically,
the bundle is trivial. However, the choice of a metric $h_0$ on $L$
equips the fiber of $\pi_*(L^k)$
at, say, $w=1$ with the corresponding $L^2$ metric,
and then, by rotation $w\to e^{i\theta}w$, to the fibers along $|w|=1$.
Then, at each $w\in A_T$, the expression
\be
\<s_\al,s_\beta\>=|w|^{2\lambda_\al^{(k)}}\delta_{\al\bar \beta}
\ee
in (\ref{bergmanmetric}) defines a metric which restricts
on the boundary of $A_T$ to the $L^2$ metric induced by $h_0$,
and is flat. The flatness is a consequence of the fact that the bundle $\pi_*(L^k)$
admits locally an orthonormal basis of holomorphic sections, namely
$w^{-\lambda_{\al}^{(k)}}s_\al$.
Returning to our original problem,
the question is then under what circumstances does $\Phi$
satisfy the degenerate Monge-Amp\`ere equation. The following theorem can be extracted
from \cite{PS06a, PS07} and provides
an answer to this question:

\begin{theorem}
\label{ansatz}
Let the set-up be as described above.
Let $A_T=\{w\in {\bf C}; e^{-T}<|w|<1\}$,
where $T$ can be both finite or infinite.
If the following two conditions
are satisfied,

\smallskip
{\rm (a)} There exists a constant $C>0$, independent of both $k$ and $\al$,
so that
\be
|\lambda_\al^{(k)}|\,\leq \, C\,k
\ee

{\rm (b)} There exists a constant $C>0$ independent of $k$
so that
\be
\int\int_{X\times A_T}\Omega_k^{n+1} \,\leq\, C\,k^{-1}
\ee
\smallskip
\noindent
then $\Phi(z,w)$ is continuous near $|w|=1$,
and we have, in the sense of pluripotential theory,
\bea
\label{MADirichlet2}
(\Omega_0+{i\over 2}\ddb\Phi)^{n+1}=0\ {\rm on}\ X\times A_T
\qquad
\Phi\big\vert_{|w|=1}=0.
\eea
\end{theorem}

Here the notion of
$(\Omega_0+{i\over 2}\ddb\Phi)^{n+1}$ in the sense of pluripotential theory, for
$\Phi$ a $\Omega_0$-plurisubharmonic function, can be defined as follows.
It suffices to define it locally, so we consider the case of
a plurisubharmonic $u$ function on ${\bf C}^n$.
Let $C$ be any non-negative $(1,1)$ closed current. Since the coefficients of $C$ are
then non-negative measures, we may define
\be
(i\ddb u)\wedge C= i\ddb (uC).
\ee
Applying this to $C=i\ddb u$, and iterating, we obtain a definition
of $(i\ddb u)^n$ as an $(n,n)$-form with a positive measure as coefficient.

\medskip
Some steps in the proof of Theorem \ref{ansatz} are as follows.
The condition (a) guarantees some a priori estimates for $\Phi_k(z,w)$,
including a crucial uniform bound for the normal derivative
of $\Phi_k(z,w)$ at the component $|w|=1$ of the boundary of $A_T$,
which will guarantee that $\Phi$ has the desired boundary value.
The condition (b) suggests that a reasonable limit of $\Phi_k(z,w)$
should have mass $0$, and thus satisfy the degenerate Monge-Amp\`ere equation.
In the foundational work \cite{BT76, BT82} of Bedford and Taylor,
it is shown that this is indeed the case, if the $\Phi_k(z,w)$ converge
either uniformly or monotonically. This is not the situation in the setting
of Theorem \ref{ansatz}, but a suitable extension of the Bedford-Taylor
pluripotential theory can be established \cite{PS06a}, which does
give the desired convergence, when $T$ is finite and we have a standard
Dirichlet problem for a domain with smooth codimension 1 boundary.
When $T=\infty$, the argument has to be supplemented by a careful
limiting process, first with $T$ finite, and then letting $T$ tend to
$\infty$ \cite{PS07}.

\smallskip

Clearly, the assumption (b) can be weakened to $\int\int_{X\times A_T}\Omega_k^{n+1}\to 0$,
in which case the conclusion of the theorem would be that there exists a subsequence
$\Phi_{k_j}(z,w)$, whose corresponding limit $\Phi(z,w)$ in the sense of
(\ref{approximation}) satisfies the equations (\ref{MADirichlet2}).

\smallskip
Other extensions of the Bedford-Taylor theory to the manifold setting
can be found in \cite{GuZe} and \cite{BlK}.

\subsubsection{Construction of geodesic segments}

We apply Theorem \ref{ansatz} in the context of geodesic segments.
In this case, we are given a second metric $h_1$ on $L$ with positive curvature
$\o_1=-{i\over 2}\ddb\log\,h_1$, and the problem is to construct
the geodesic joining $h_0$ to $h_1$ (whose existence has been proved by
the Theorem \ref{chentheorem} of X.X. Chen, using a priori estimates and elliptic
regularization). In this case, let $\us^{(1)}=\{s_\al^{(1)}(z)\}_{\al=0}^{N_k}$
be a basis for $H^0(X,L^k)$, orthonormal this time with respect to
the metric $h_1$ and the volume form $\o_1^n$. Without loss of generality,
we may assume that the change of bases from $\us$ to $\us^{(1)}$
is given by a diagonal matrix,
\be
s_{\al}^{(1)}= e^{\lambda_\al^{(k)}}s_\al.
\ee
Let $\Phi_k(z,w)$ be defined as in (\ref{approximation}), with this
choice of weights $\lambda_\al^{(k)}$. Then it is shown in \cite{PS06a}
that the two conditions (a) and (b) of Theorem \ref{ansatz}
are satisfied. It may be instructive to see how to verify condition (b),
assuming condition (a), since
we have already at our disposal all the tools needed to
estimate the Monge-Amp\`ere masses in this context.
By (\ref{MAmasses}), this reduces to
computing $\dot F_\o^0$ at the two end points $t=0$ and $t=1$, and
since $\phi_k(z,t)$ is a path inside ${\cal K}_k$, the
formulas for $\dot F$ in section \S apply. Thus we have
\bea
\label{F0MA}
\int\int_{X\times A_T}
\Omega_k^{n+1}
=
\int_X \dot\phi_k(1)\o_{\phi(1)}^n-\int_X\dot\phi_k(0)\o_{\phi(0)}^n
\eea
and hence, more explicitly,
\bea
{2\over k^{n+1}}
\int_X \sum_{\al=0}^{N_k}
\lambda_\al^{(k)}|s_\al^{(1)}(z)|^2h_1(k)^k \o_1(k)^n
-
{2\over k^{n+1}}
\int_X \sum_{\al=0}^{N_k}
\lambda_\al^{(k)}|s_\al(z)|^2h_0(k)^k\o_0(k)^n.
\eea
Applying the Tian-Yau-Zelditch theorem, we obtain easily the following
asymptotics,
\be
\label{dotFasymptotics}
{2\over k^{n+1}}
\int_X \sum_{\al=0}^{N_k}
\lambda_\al^{(k)}|s_\al(z)|^2h_1(k)^k \o_1(k)^n
=
{2\over k^{n+1}}
\sum_{\al=0}^{N_k}\lambda_\al^{(k)}
+
O({1\over k^{n+2}})\cdot N_k\cdot {\rm max}_\al |\lambda_\al^{(k)}|
\ee
as well as an analogous expression with $h_1\leftrightarrow h_0$.
The leading term in this expression cancels out between the contributions
from $h_1$ and from $h_0$, and the second term is $O(k^{-1})$ in view
of condition (a). Thus condition (b) is satisfied. Theorem
\ref{ansatz} implies then that the corresponding
$\Phi(z,w)$ defines then a solution of the equation (\ref{MADirichlet2}).
Furthermore, the roles of $h_0$ and $h_1$ are clearly reversible, and
it follows that $\Phi(z,w)$ satisfies the desired boundary condition
also when $|w|=e^{-1}$, and is actually a solution
of the Dirichlet problem (\ref{MADirichlet}). By uniqueness, it
must coincide with the solution provided by Theorem \ref{chentheorem}. Thus we have the
following theorem:

\begin{theorem}
\label{bergmanapproximation}
Let $L\to X$ be a positive line bundle over a compact complex manifold $X$.
Let $h_0,h_1$ be two metrics on $L$ with positive curvatures $\o_i=-{i\over 2}\ddb\log\,h_i$.
Let $\us$, $\us^{(1)}$ be two bases of $H^0(X,L^k)$, orthonormal with
respect to the metrics and volume forms $h_0,\o_0$ and $h_1,\o_1$ respectively.
Assume without loss of generality that the matrix of change of bases is diagonal,
with eigenvalues $\lambda_\al^{(k)}$. Then the $C^{1,1}$ geodesic joining
$h_0$ and $h_1$ can be obtained by the construction {\rm (\ref{approximation})}.
\end{theorem}

\subsubsection{Construction of geodesic rays}

We turn next to the construction of geodesic rays. Let $L\to X$ be a
positive line bundle over a compact complex manifold as before, and
let $\cT: \cL\to\cX\to {\bf C}$ be a test configuration, in the
sense of Definition \ref{testconfiguration}. Let $A_k$ be the
traceless endomorphism of $H^0(X_0,L_0^k)$, as defined in \ref{B}
and subsequent line, and let $\lambda_\al^{(k)}$ be its eigenvalues. Then \cite{PS07}

\begin{theorem}
\label{bergmanrays}
With this set-up, and this choice of weights $\lambda_\al^{(k)}$,
the expressions {\rm (\ref{bergmanmetric})} and {\rm (\ref{approximation})}
produce a generalized geodesic ray starting at $h_0$, i.e.,
$\Phi(z,w)$ is continuous and equal to $0$ at $|w|=1$,
and we have, in the sense of pluripotential theory,
\be
(\Omega_0+{i\over 2}\ddb\Phi)^{n+1}=0
\quad{\rm on}\quad X\times \{w\in {\bf C};0<|w|<1\}.
\ee
The geodesic is non-constant when the test configuration $\cT$ is non-trivial.
\end{theorem}

To establish the theorem, we apply Theorem \ref{ansatz} and verify conditions (a)
and (b). Once again, in this exposition, we assume condition (a), and concentrate on condition (b).
Again, by (\ref{F0MA}), we can write
\be
\int\int_{X\times A_\infty}\Omega_k^{n+1}
=
{\rm lim}_{T\to\infty}\int_X\dot\phi(T)\o_k(T)^n
-
\int_X\dot\phi(0)\o_k(0)^n
\ee
We expand the second term on the right hand side, using the Tian-Yau-Zelditch
theorem as before. In view of (\ref{dotFasymptotics}) and the fact that the matrix $A_k$
is traceless, this term is $O(k^{-1})$. By Lemma \ref{FTmu}, we have
\be
{\rm lim}_{T\to\infty}\int_X\dot\phi(T)\o_k(T)^n
={1\over k}F
\ee
where $F$ is the Donaldson-Futaki invariant. In particular, this shows
that this term is also $O(k^{-1})$, establishing condition (b).
Thus Theorem \ref{bergmanrays} follows from Theorem \ref{ansatz}.

\medskip
The construction of Theorem \ref{ansatz} is arguably canonical,
and thus, given a test configuration $\cT$, we have a canonical way
of associating to any point $h_0\in{\cal K}$ a generalized vector,
namely the initial velocity vector at $h_0$ of the geodesic we
have just constructed. In this sense, a test configuration defines a
generalized vector field on the space ${\cal K}$ of K\"ahler potentials.
It would be very valuable to be able to write this vector field down more
explicitly in terms of data from the test configuration $\cT$. Indeed,
it would provide a valuable model for how to relate behavior at $\infty$
to behavior well inside ${\cal K}$, a ubiquitous underlying theme
in the problem of stability and constant scalar curvature metrics.

\subsubsection{Variations on the Ansatz and rates of convergence}

Theorems \ref{bergmanapproximation} and \ref{bergmanrays}
show that geodesics in ${\cal K}$ can be approximated by geodesics in ${\cal
K}_k$.
It is of great interest to understand
how good such approximations can be,
and in particular, to analyze more precisely the approximation given in
(\ref{approximation}). Replacing the potentials $\Phi_k$ in
(\ref{approximation}) with potentials $\tilde\Phi_k$ associated instead
to the bundle $L^k\otimes K_X$,
Berndtsson \cite{Be} has obtained
another version of (\ref{approximation}) with
$C^0$ convergence and with
precise error bounds $O(k^{-1}\log\,k)$.
If we stay instead with $\Phi_k$, the sharpest results to date have been
provided by
Song and Zelditch \cite{SZ1, SZ2}, in the context of toric varieties:

\medskip
\begin{theorem}
Let $L\to X$ be an ample toric line bundle over a compact toric variety.
Let $h_0$, $h_1$ be toric Hermitian metrics on $L$ with positive curvatures
$\o_0$, $\o_1$.
Then the approximation $\Phi_k$ described in {\rm (\ref{bergmanmetric})}
converges in $C^2(X)$ to the geodesic joining $h_0$ and $h_1$.
\end{theorem}

We note that, for toric varieties, the geodesic equation becomes a linear
equation
for the Legendre transform of the potentials \cite{GuD,SZ1}, so that the
geodesic joining $h_0$ and
$h_1$ is known to be smooth if $h_0$ and $h_1$ are smooth.

\medskip

A similar sharp analysis of the construction of geodesic rays associated to a test
configuration by Theorem \ref{bergmanrays} has appeared very recently
in \cite{SZ3}.
In particular, it is shown there that the regularity $C^{1,1}$ for such rays
is optimal. Some remarkable and potentially far-reaching relations
between the Ansatz (\ref{bergmanmetric})
for complex Monge-Amp\`ere equations
and ideas from the theory of large deviations
are also brought to light in this paper.
Other unexpected connections between approximations of the form
(\ref{bergmanmetric}) and classical topics such as Bernstein polynomials
and Dedekind-Riemann sums over
lattice points in polytopes can be found in \cite{Z07}.

\end{document}